\title{On geometrically defined extensions of the Temperley--Lieb category 
 in the Brauer category}
\author{
Z K\'ad\'ar,
P P Martin %\and
and 
S Yu \\ %\myleedsaddress
Department of Pure Mathematics, University of Leeds}
\date{}
\documentclass[12pt]{article} 
\usepackage{epic} 
\usepackage{eepic} 
\usepackage{latexsym}
\usepackage{amsthm}
\usepackage{amssymb}
\usepackage{graphicx}
\usepackage[all]{xy}
\providecommand{\noglossaryignore}[1]{}
\newcommand{\globalglossaryentry}[3]{\makebox[1.5in][l]{\tt $\backslash${#1}} 
\makebox[1.1in][l]{{$#2$}} \makebox[2.5in][l]{{#3}}\newline} 
\newcommand{\newcommandabbreviation}[3]{\newcommand{#1}{#2}%
\noglossaryignore{\globalglossaryentry{#3}{#2}{}}}
\newcommand{\renewcommandabbreviation}[3]{\renewcommand{#1}{#2}%
\noglossaryignore{\globalglossaryentry{#3}{#2}{}}}
\newcommand{\newcommandmacro}[4]{\newcommand{#1}{#2}%
\noglossaryignore{\globalglossaryentry{#3}{#2}{#4}}}
\newcommand{\gge}[3]{\noglossaryignore{\globalglossaryentry{#1}{#2}{#3}}}
\newcommand{\myaddress}%
{\parbox{3in}{\footnotesize \begin{center} 
Mathematics Department, City University, \\  
Northampton Square, London EC1V 0HB, UK.\end{center}}}
\newcounter{minidef}[section]
\renewcommand{\theminidef}{\thesection.\arabic{minidef}}
\newcommand{\mdef}{\refstepcounter{minidef} 
\medskip \noindent ({\bf \theminidef}) }

\newcommand{\mystufffont}{\textsc} %% this sets up font for 
\newtheoremstyle{pu}% name
{7pt}%
{7pt}%
{\it}% body font
{}% indent amount
{}% {\textsc}% headfont
{.}% puntuation
{ }% gap
{\thmnumber{({\bf #2}) }\thmname{\textsc{#1}}\thmnote{#3}}% headspec
\newtheoremstyle{puu}% name
{3pt}%
{3pt}%
{\rm}% body font
{}% indent amount
{}% {\textsc}% headfont
{.}% puntuation
{ }% gap
{\thmnumber{({\bf #2}) }\thmname{\textsc{#1}}\thmnote{#3}}% headspec

\theoremstyle{pu}
\newtheorem{mmpr}[minidef]{Proposition}
\newcommand{\mpr}[1]{\begin{mmpr} #1 \end{mmpr}}
\newtheorem{mmco}[minidef]{Corollary}
\newcommand{\mco}[1]{\begin{mmco} #1 \end{mmco}}

\theoremstyle{puu}

\newtheorem{mex}[minidef]{Example}
\newtheorem{mrem}[minidef]{Remark}

\newcommand{\mupr}{\smallskip\noindent{\mystufffont{Proposition.}} }

\newcounter{minicapt}

\newtheorem{de}[minidef]{Definition}     
\newtheorem{pr}[minidef]{Proposition} 
\noglossaryignore{GREEK ETC.\newline}
\newcommandabbreviation{\e}{\epsilon}{e}        
\newcommandabbreviation{\lam}{\lambda}{lam}  
\newcommandabbreviation{\la}{\langle}{la}        
\newcommandabbreviation{\ran}{\rangle}{ran}
\newcommandabbreviation{\ha}{\#}{ha}             
\newcommandabbreviation{\rmap}{\rightarrow}{rmap}
\newcommandabbreviation{\aaa}{\alpha}{aaa}        
\newcommandabbreviation{\ab}{\alpha,\beta}{ab}
\newcommandabbreviation{\aab}{a(\ab )}{aab}       
\noglossaryignore{\newline RINGS\newline}
\newcommandabbreviation{\HH}{H \!\!\! I}{HH}               % Hecke algebra
\newcommandabbreviation{\C}{\mathbb C}{C}
\newcommandabbreviation{\N}{\mathbb N}{N}   %new AMS versions (was \Bbb)!
\newcommandabbreviation{\Z}{\mathbb Z}{Z}      % AMS versions!!!!!!
\renewcommandabbreviation{\Re}{\mathbb R}{Re}
\newcommandabbreviation{\R}{{\mathbb R}}{R}
\newcommandabbreviation{\Q}{\mathbb Q }{Q}
\renewcommandabbreviation{\H}{\mathbb H }{H}
\noglossaryignore{\newline SYMMETRIC GROUP\newline}
\def\Sym(#1){\Sigma(#1)}                   % generic symbol for symmetric group
\gge{Sym(-)}{\Sym(-)}{symmetric group on - objects}
\def\Sy(#1){\Sigma_{#1}}                   % symmetric group irrep.
\gge{Sy(-)}{\Sy(-)}{symmetric group irreducible -}
\def\sym(#1){\mbox{\LARGE s}(#1)}        % another generic symbol for
\gge{sym(-)}{\sym(-)}{symmetric group on - objects (variant)}
\def\sy(#1){\mbox{\LARGE s}({#1})}        % another symm gp irrep. 
\gge{sy(-)}{\sy(-)}{symmetric group irreducible - (variant)}
\newcommandmacro{\cs}{\C \, \sy(n)}{cs}{symmetric group algebra over $\C$}
\noglossaryignore{\newline PARTITIONS/SETS\newline}
\newcommand{\Nset}[1]{\underline{#1}}
\gge{Nset\{-\}}{\Nset{-}}{set of natural numbers to -} %%% (default n)}
\def\nset(#1){ \{ #1 \}_{ \underline{n} }} % the set {#1}_{n}
\gge{nset(-)}{\nset(-)}{a set $-\times\Nset$}
\def\ul(#1){_{\underline{#1}}}             % _underline #1
\gge{ul(-)}{{}\ul(-)}{subscript underline -}
\def\Ee(#1){{\bf E}_{#1}}                  % set of equiv. relations
\gge{Ee(-)}{\Ee(-)}{set of equivalence relations on set -}
\def\Eee(#1){{\bf E}_{\{ #1 \}_{\underline{n}}}}   %ditto for nset
\gge{Eee(-)}{\Eee(-)}{ditto for nset}
\def\Een(#1,#2){{\bf E}_{\{ #1 \}_{\underline{#2}}}}   %ditto for n+1set
\def\Ssn(#1,#2){{\bf S}_{\{ #1 \}_{\underline{#2}}}}   %partitions for n+1set
\def\Ss(#1){{\bf S}_{#1}}                  % set of partitions
\def\Sss(#1){{\bf S}_{\{ #1 \}_{\underline{n}}}}   %ditto for nset
\def\bbc(#1){((\beta_1)(\beta_2)...(\beta_{#1}))}      % beta singletons 
\newcommandmacro{\Ln}{{\Gamma}^{n}}{Ln}{large index set}
\newcommandmacro{\LnQ}{{\Gamma}^{n}_Q}{LnQ}{index set}
\newcommandmacro{\Zz}{\zeta}{Zz}{`shape' function}
\noglossaryignore{\newline PARTITION ALGEBRA\newline}
\def\ka(#1){\kappa_{#1}}                   % maps AP_{n=#1}A to P_{n-1}
\def\Sm(#1){\Sigma_{#1}}                   % image of S_lamda in P_n
\newcommandmacro{\com}{\bullet}{com}{bullet composition}
\newcommandmacro{\enm}{\; e^n(\! m\! ) \;}{enm}{product of idempotents}
\def\Ai(#1){ A^{ #1 \cdot } }              % A_i
\def\Aij(#1,#2){ A^{ #1  #2 } }            % A_ij
\newcommandmacro{\One}{\mbox{\bf $1 \!\!\! 1$}}{One}{algebra unit 1}
\newcommandmacro{\Bp}{B_p}{Bp}{partition basis}
\def\Bb(#1){B_p[#1]}                       % partition basis
\def\Pp(#1){P_n[#1]}                       % left module
\def\Ps(#1){P_n[#1] \! /}                  % left module
\newcommandmacro{\Ph}{\hat{P}}{Ph}{P hat  algebra}
\def\Is(#1){\sim^{#1}}                     % is equivalent under a to
\noglossaryignore{\newline MODULES\newline}
\def\Wm(#1){{\cal S}_{#1}}                 % Weyl module S 
\gge{Wm(-)}{\Wm(-)}{Weyl module with index -}
\def\wm(#1,#2){{}_{#1}{\cal S}_{#2}}       % Weyl module nS
\gge{wm(-1,-)}{\wm(-1,-)}{Weyl module with index -}
\def\Ind(#1,#2,#3){\mbox{Ind}_{#1}^{#2}#3} % induction
\gge{Ind(-1,-2,-)}{\Ind(-1,-2,-)}{induction}
\def\Res(#1,#2,#3){\mbox{Res}_{#1}^{#2}#3} % restriction
\gge{Res(-1,-2,-)}{\Res(-1,-2,-)}{restriction}
\newcommandabbreviation{\weyl}{standard}{weyl}
\newcommandabbreviation{\mod}{\mbox{mod}}{mod}
\newcommandabbreviation{\head}{\mbox{head }}{head}
\newcommandabbreviation{\Weyl}{Weyl}{Weyl}
\def\SS(#1){{\cal S}_{#1}}                 % Specht/Weyl module
\gge{SS(-)}{\SS(-)}{Specht/Weyl module index -}
\def\LL(#1){{\cal L}_{#1}}                 % Simple module
\gge{LL(-)}{\LL(-)}{simple module index -}
\noglossaryignore{\newline FUNCTORS/MAPS\newline}
\newcommandmacro{\Gg}{{\cal G}}{Gg}{G Functor}
\newcommandmacro{\Fg}{{\cal F}}{Fg}{F Functor}
\newcommandmacro{\ra}{\rightarrow}{ra}{}
\def\ses(#1,#2,#3){0\ra #1 \ra #2 \ra #3 \ra 0}   %short exact sequence
\gge{ses(1,2,-)}{\ses(1,2,-)}{\hspace{.5in} short exact sequence}
\def\starr(#1){ \stackrel{ #1 }{\longrightarrow} }
\gge{starr(-)}{\starr(-)}{}
\newcommandmacro{\doublerightarrow}{\; -\!\!\! -\!\!\!\!\!\! \gg \;}
{doublerightarrow}{}%{ $---->>$ }
\noglossaryignore{\newline PARTITION ALGEBRA MAPS\newline}
\newcommandmacro{\smap}{s}{smap}{`inclusion' map}
\newcommandmacro{\tmap}{t}{tmap}{$ P_n -> S_n$}
\newcommandmacro{\pmap}{\psi}{pmap}{$ S_n -> P_n $}
\noglossaryignore{\newline MISC.\newline}
\def\Amap(#1){{\cal A}_{#1}}               % variant inclusion A_Gamma
\gge{Amap(-)}{\Amap(-)}{}
\def\Rr(#1){R_{#1}}                        % restriction of E
\gge{Rr(-)}{\Rr(-)}{restriction of E}
\def\Cr(#1){C_{#1}}                        % restriction of E to N
\gge{Cr(-)}{\Cr(-)}{restriction of E to N}
\newcommandmacro{\Tm}{{\cal T}}{Tm}{Transfer Matrix}
\def\On(#1){{\cal I}_{#1}}
\gge{On(-)}{\On(-)}{}
\newcommandmacro{\UU}{\underline{\sqcup}}{UU}{}  
\newcommandmacro{\UUU}{\sqcup}{UUU}{}  
\newcommandmacro{\Vq}{V_Q^{\otimes n}}{Vq}{Potts config. space}
\def\bs(#1,#2){\mbox{{\Large $\ast$}}^{#1}_{#2}}  % general plumbing multiplier
\gge{bs(-,-)}{\bs(-,-)}{general plumbing multiplier}
\newcommand{\ignore}[1]{}
\gge{ignore\{-\}}{\ignore{-}}{ignore argument!}
\def\choo(#1,#2){ \left( \begin{array}{c} #1 \\ #2 \end{array} \right) } %choose
\gge{choo(-1,-)}{\choo(-1,-)}{choose}
\newcommand{\Qed}{$\Box$}%{Qed}{QED}
\gge{Qed}{\mbox{\Qed}}{QED}
\def\staq(#1){\stackrel{#1}{=}}            % stack =
\gge{staq(-)}{\staq(-)}{}
\def\stam(#1){\stackrel{#1}{\rightarrow}}  % stack ->
\gge{stam(-)}{\stam(-)}{}
\def\mat{ \left( \begin{array} }    
\def\tam{ \end{array}  \right) }
\gge{mat/tam}{...}{matrix delimiters}
\newcommand{\beq}{\begin{equation} }
\def\eql(#1){ \begin{equation} \label{#1} 
}
\newcommand{\eq}{\end{equation} }
\def\eqal(#1){\begin{eqnarray} \label{#1} }
\def\eqa{\end{eqnarray} }
\def\lab(#1){\label{#1}
}
\def\prl(#1){ \begin{pr} \label{#1} 
}
\def\del(#1){ \begin{de} \label{#1} 
}
\gge{smeq\{-\}}{...}{small equation}
\gge{fneq\{-\}}{...}{very small equation}
\noglossaryignore{\newline HECKE/BLOB\newline}
\newcommandmacro{\Hnq}{H_n(q)}{Hnq}{ * freestanding symbol}
\newcommandmacro{\Hn}{H_n}{Hn}{      *-mod etc.}
\newcommandmacro{\A}{{\cal A}}{A}{}
\newcommandmacro{\Cwts}{C}{Cwts}{}
\newcommandmacro{\CA}{{\cal A}}{CA}{}

\newcommandmacro{\calA}{{\cal A}}{calA}{}
\newcommandmacro{\modi}{\mbox{Mod} }{modi}{was mod not modi!}
\newcommandmacro{\Wgen}{{\Bbb S}}{Wgen}{}
\def\ol(#1){\overline{#1}}
\newcommandmacro{\st}{\mbox{St}}{st}{}
\def\CMult(#1,#2){(#1:#2)}
\def\CM(#1,#2){( #1 : #2 )}
\def\FMult#1,#2{(#1:#2)}
\def\CF#1,#2{(#1:#2)}

\newcommandmacro{\Top}{\mbox{Top}}{Top}{}
\newcommandmacro{\Soc}{\mbox{Soc}}{Soc}{}
\newcommandmacro{\Head}{\mbox{Head}}{Head}{}
\newcommandmacro{\Filt}{{\cal F}}{Filt}{}
\newcommandmacro{\Mod}{\mbox{mod}}{Mod}{}
\newcommandmacro{\Resi}{\mbox{Res }}{Resi}{was without i!}
\newcommandmacro{\Indi}{\mbox{Ind }}{Indi}{was without i!}
\def\RR(#1,#2){R^{#1}_{#2}}   %projection after restriction (after projection)
\def\TT(#1,#2){T^{#1}_{#2}}   %translation

\def\Chi{\chi}

\newcommandmacro{\Ann}{\mbox{Ann}}{Ann}{}
\newcommandmacro{\Cen}{\mbox{Cen}}{Cen}{}
\newcommandmacro{\End}{\mbox{End}}{End}{}
\newcommandabbreviation{\semisimple}{semisimple}{semisimple}
\newcommandabbreviation{\Bratteli}{Bratteli}{Bratteli}
\newcommandabbreviation{\JBC}{Jones Basic Construction}{JBC}
\newcommandabbreviation{\pa}{partition algebra}{pa}
\newcommandabbreviation{\TM}{transfer matrix}{TM}
\newcommandabbreviation{\PM}{Potts model}{PM}
\newcommandabbreviation{\QSC}{quantum spin chain}{QSC}
\newcommandabbreviation{\Hamiltonian}{Hamiltonian}{Hamiltonian}
\newcommandabbreviation{\YS}{Young symmetrizer}{YS}
\theoremstyle{pu}   %%% defined in martin07.tex
\newtheorem{mmlem}[minidef]{Lemma}
\newcommand{\mlem}[1]{\begin{mmlem} #1 \end{mmlem}}
\newtheorem{mmcor}[minidef]{Corollary}
\newcommand{\mcor}[1]{\begin{mmcor} #1 \end{mmcor}}
\newtheorem{mmth}[minidef]{Theorem}
\newcommand{\mth}[1]{\begin{mmth} #1 \end{mmth}}

\newcommand{\printglossary}{} % needed if glossary not defined

\newcommand{\dd}{{\sf d}}  %%% hook for line distances
\newcommand{\hht}{{\rm ht}}  %%% hook for heights
\newcommand{\Nm}{\N_-}   %%%   {-1,0,1,2,...}

\newcommand{\ppicture}{picture} %% polygonal case

\newcommand{\pvert}{arc-vert}

\newcommand{\catP}{{\mathcal P}}
\newcommand{\catB}{{\mathcal B}}

\newcommand{\kk}{{\mathsf k}}  %% comm ring
\newcommand{\units}{\times}  %% k^\units is gp of units

\begin{document} \maketitle
%{{{ draft level and draft header 
 \newcommand{\ignoreifnotdraft}[1]{\ignore{#1}}
\ignoreifnotdraft{
%%% PAGEHEADINGS (REMOVED FOR JOURNAL SUBMISSION):
\pagestyle{myheadings} \markboth{Draft}{\today}
}

\ignoreifnotdraft{
\noindent
\begin{tabular}{l} 
  {\tiny \filename (Draft)} \hspace{4.4in} Jan 1997 \\ \hline 
\end{tabular}
}
%}}}
%{{{ %abstract
%
\begin{abstract}
We define an infinite chain of subcategories of the partition 
category by introducing the left-height ($l$) of a partition. 
For the Brauer case, the chain starts with the Temperley-Lieb 
($l=-1$) and ends with the Brauer ($l=\infty$) category. The End 
sets are algebras, i.e., an infinite tower thereof for each $l$, 
whose representation theory is studied in the paper. 
\end{abstract}
\ignoreifnotdraft{ \newpage } %%% Journal submission only %%%
\printglossary
%}}}
\section{Introduction}
%{{{ intro
%{{{ INTRO
%{{{ intro
\newcommand{\TL}{Temperley--Lieb} 
\newcommand{\BMW}{Birman--Murakami--Wenzl}

The partition algebra and its Brauer and \TL\ (TL) subalgebras 
\cite{Brauer37,TemperleyLieb71,Martin94} 
%all 
have many applications and 
a rich representation theory 
(see e.g. 
\cite{Brauer37,Brown55,Wenzl88b,Martin96,Mazorchuk95}
and references therein).
In particular  each 
representation theory 
has  %a distinct and 
an intriguing geometrical characterisation 
\cite{Jantzen87,CoxDevisscherMartin0609,Martin91,MartinWoodcock2000}
%embedded in the wider representation theory of the 
%Brauer algebra \cite{Brauer37,Brown55,Wenzl88b} and
(in the Brauer case also embracing 
the \BMW\ (BMW) algebra 
\cite{Murakami87,BirmanWenzl89,Kauffman91,Morton10,Yu05}).
The TL case can be understood in terms of Lie Theoretic notions of
alcove geometry and geometric linkage, via generalised Schur-Weyl
duality \cite{Martin92,DuParshallScott98,Jantzen87}, 
but the Brauer case is much richer
\cite{Martin09a}
and, 
although its complex 
representation theory is now {\em intrinsically} well-understood, 
the paradigm for the corresponding notions
%in the Brauer case is much richer and 
is more mysterious. 
Here we introduce a sequence of (towers of) algebras $J_{l,n}$ 
which  interpolate between the TL algebra and the Brauer algebra
 as algebras.  %%The
A particular aim is to  study the geometry in their representation theory
by
lifting this new connection 
to the representation theory level.
%to the level of, and hence study, the geometry in their 
%representation theory. 
To this end 
we  %%and 
investigate the representation theory of the new algebras using 
their amenability to tower of recollement (ToR)
\cite{CoxMartinParkerXi06}
and monoid methods 
\cite{Putcha98}. 
The representation theory for large $l,n$ eventually becomes very
hard, but we are able to prove a number of useful general results,
and results over the complex field. 
%In particular
For example we obtain the
`generic' semisimple structure in the sense of \cite{ClineParshallScott99}.

%}}}
%{{{ further motive

By way of further motivation (although we will not
develop the point here) we note that 
both the Brauer and TL algebras 
%have many applications. For example they
provide solutions to the 
Yang--Baxter (YB) equations
\cite{Baxter,ReshetikhinTuraev90}. 
In addition to their interest from a representation theory
perspective, 
our new algebras can be seen as ways to address the problem 
of contruction of natural solutions to
the boundary YB equations in the TL setting (generalising the blob
approach and so on --- see e.g. \cite{DoikouMartin03}
and references therein). 
A paradigm here is the XXZ spin chain --- a `toy' model of quantum
mechanical interacting spins on a 1-dimensional spacial lattice
derived from the Heisenberg model \cite{Bethe31};
see e.g. \cite[Ch.6]{LiebMattis66}. 

%}}}
%{{{ ignore Here having contructed

%}}}
%{{{ outline

\medskip

An outline of the paper is as follows. The partition category has a
basis of set partitions, and the Brauer and \TL\ categories are
subcategories with bases of certain restricted partitions. In
particular the \TL\ category has a basis of {\em non-crossing}
partitions. 
Here we provide a classification of partitions generalising the 
geometrical notion
of non-crossing. Many such games are possible in principle
(see e.g. \cite{CautisJackson10}), but we
show that our classification (like non-crossing) is preserved under
the partition category composition. 
This closure theorem allows us to define a
sequence of new `geometrical' subcategories.
Next we turn to our motiviating aim:  investigation of geometric features
in the algebraic representation theory contained in these categories. 
We focus in particular here on the extensions of the \TL\ category in
the Brauer category. 
In this paper we establish a framework for
modular representation theory of the corresponding towers of algebras.
In the case that is modular over $\C$ in the sense of \cite{Brauer39}
we prove that the algebras are generically semisimple. 
We observe
an intriguing subset of %the  %%cases in 
parameter values for 
which they are not semisimple,
distinct from both TL and Brauer cases. 
We conclude by determining branching rules, and hence give
combinatorial constructions for the ranks of these algebras.

%}}}
%{{{ contour etc

The TL algebra has %, in turn, 
a sequence of known generalisations 
using its characterisation via an embedding of pair partitions in the plane
--- the blob algebras and the contour algebras 
\cite{MartinGreenParker07};
as well as various beautiful generalisations due to 
R.~Green et al \cite{Green98,FanGreen99}, 
tom~Dieck \cite{tomDieck98} and others.
The blob algebra also has a rich geometrically-characterised
representation theory \cite{MartinWoodcock2000}. 
However  %%neither the contour algebras nor any of the other 
none of the 
previously known cases 
%do not 
serve to interpolate between the TL algebra and the Brauer algebra.
%(either as algebras, or representation theoretically). 

%}}}
%{{{ Jobs
%***
%}}}
%}}}
\subsection{Notations and pictures for set partitions}
%{{{ macros

\newcommand{\vac}{f}    %% hook for vacuum function on partitions
\newcommand{\vloop}{\Pi^{\vac}}  %% hook for vacuum loop number
\newcommand{\Part}{P}   %% hook for set of set partitions
\newcommand{\gr}{{\mathcal G}}  %% class of graphs
\newcommand{\hashp}{\#^{{\mathsf p}}}  %% propagating no.

\newcommand{\gloss}[2]{}

%}}}
%{{{ 1

We  need to recall 
a pictorial (and so {\em en passant} geometric) realisation of the
partition algebra (i.e. of set partitions). This realisation is 
in common use  %standard, as an informal device
(see e.g. \cite{Martin96}), but we will need to develop it more formally.

\gloss{$P(T)$}{ set of partitions of set $T$}

\mdef \label{de:PiS}
If $T$ is a set then $\Part(T)$ denotes the set of partitions of $T$.
Noting the standard bijection between $\Part(T)$ and 
the set of equivalence relations on $T$ 
we write $a \sim^p b$ when $a,b$ in the same part in $p \in \Part(T)$.

Suppose that $p \in \Part(T)$ and 
% is a set partition of some set $T$. 
 $S \subset T$.
Write $p|_S$ for the restriction of $p$ to $S$. 
Write 
$\vac_S(p) \; := \; \# \{ \pi \in p \; | \; \pi \cap S = \emptyset \}$,
the number of parts of $p$ that do not intersect $S$. 
(Here we follow \cite[Def.20]{Martin94}. 
See also e.g. \cite{Martin08a}.)

%}}}
%{{{ graph1

\mdef \label{de:PiS2}
Let $\gr$ denote the class of graphs;  
$\gr(V)$ the subclass of graphs on vertex set $V$;
and $\gr[S]$ the subclass of graphs whose vertex set {\em contains}
set $S$.
Define
$$
\Pi : \gr(V) \rightarrow \Part(V)
$$
by $v \sim^{\Pi(g)} v'$ if $v,v'$ are in the same  connected component in the
graph $g$. 
Define
$$
\Pi_S : \gr[S] \rightarrow \Part(S)
$$
by $\Pi_S(g) = \Pi(g)|_S$.
%
%}}}
%{{{ drawings (ignore)

%}}}
%{{{ fig

\begin{figure}
\[
\includegraphics[width=1.62in]{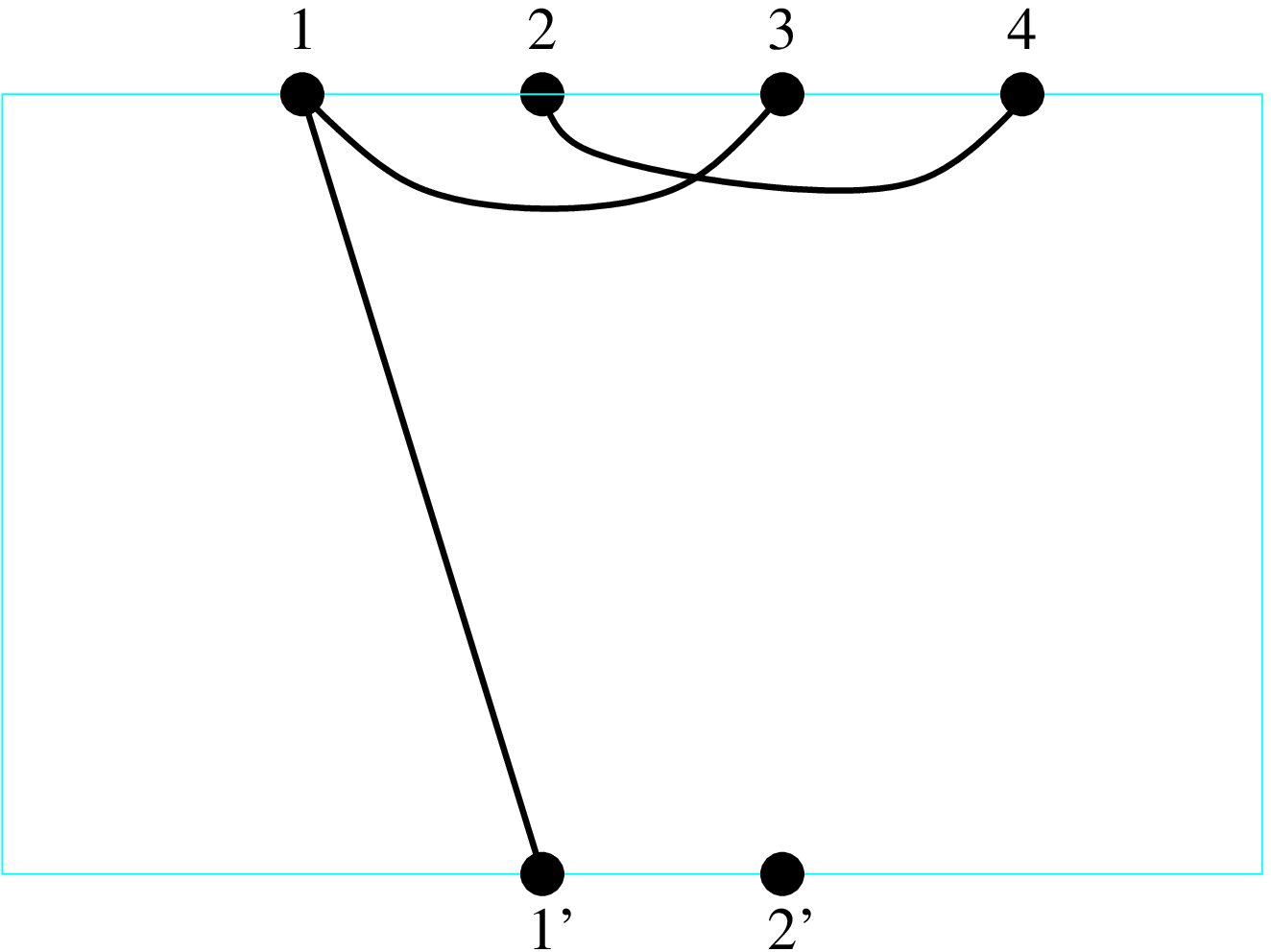}
\qquad
\includegraphics[width=1.62in]{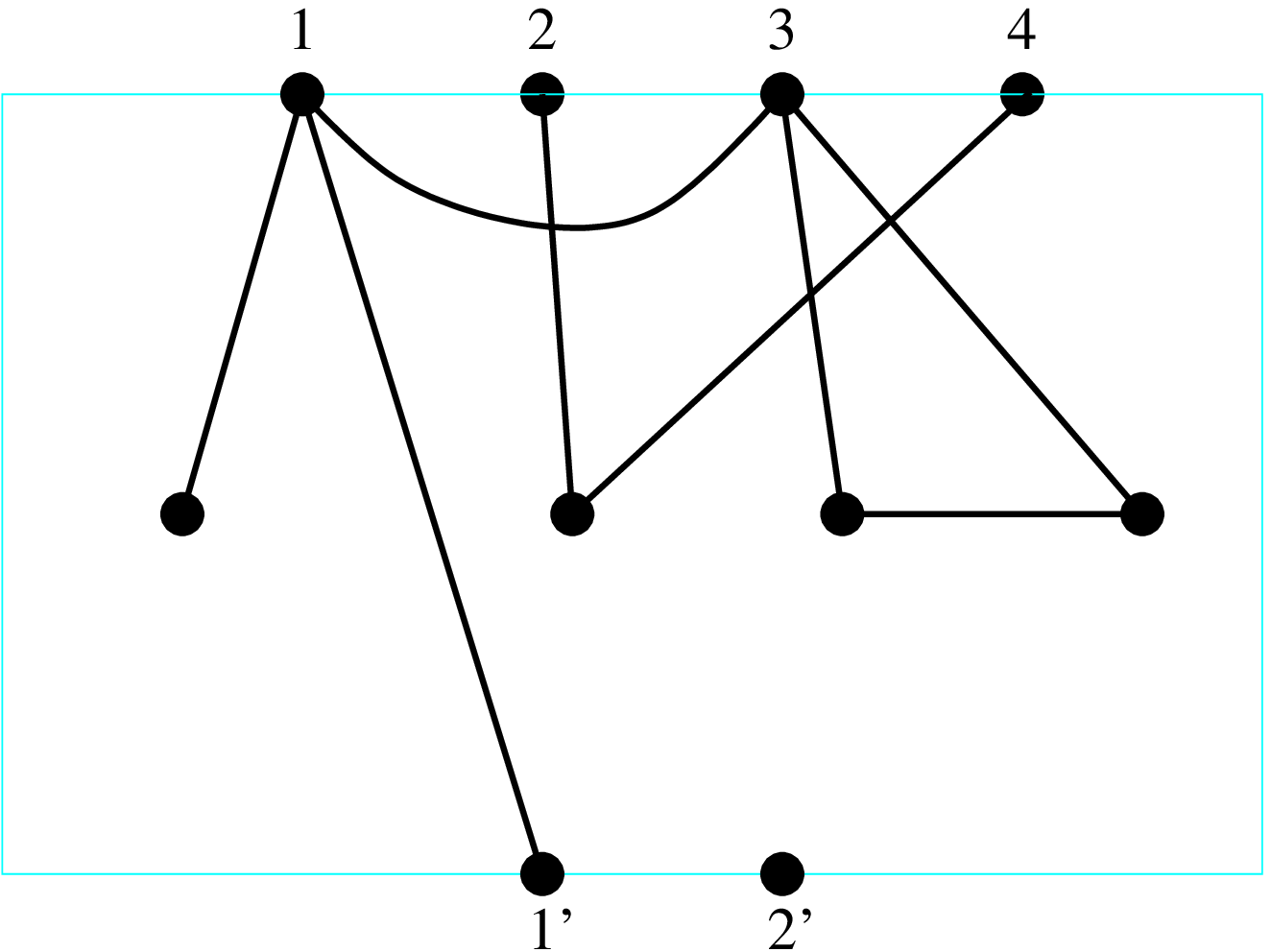}
\qquad
\includegraphics[width=1.62in]{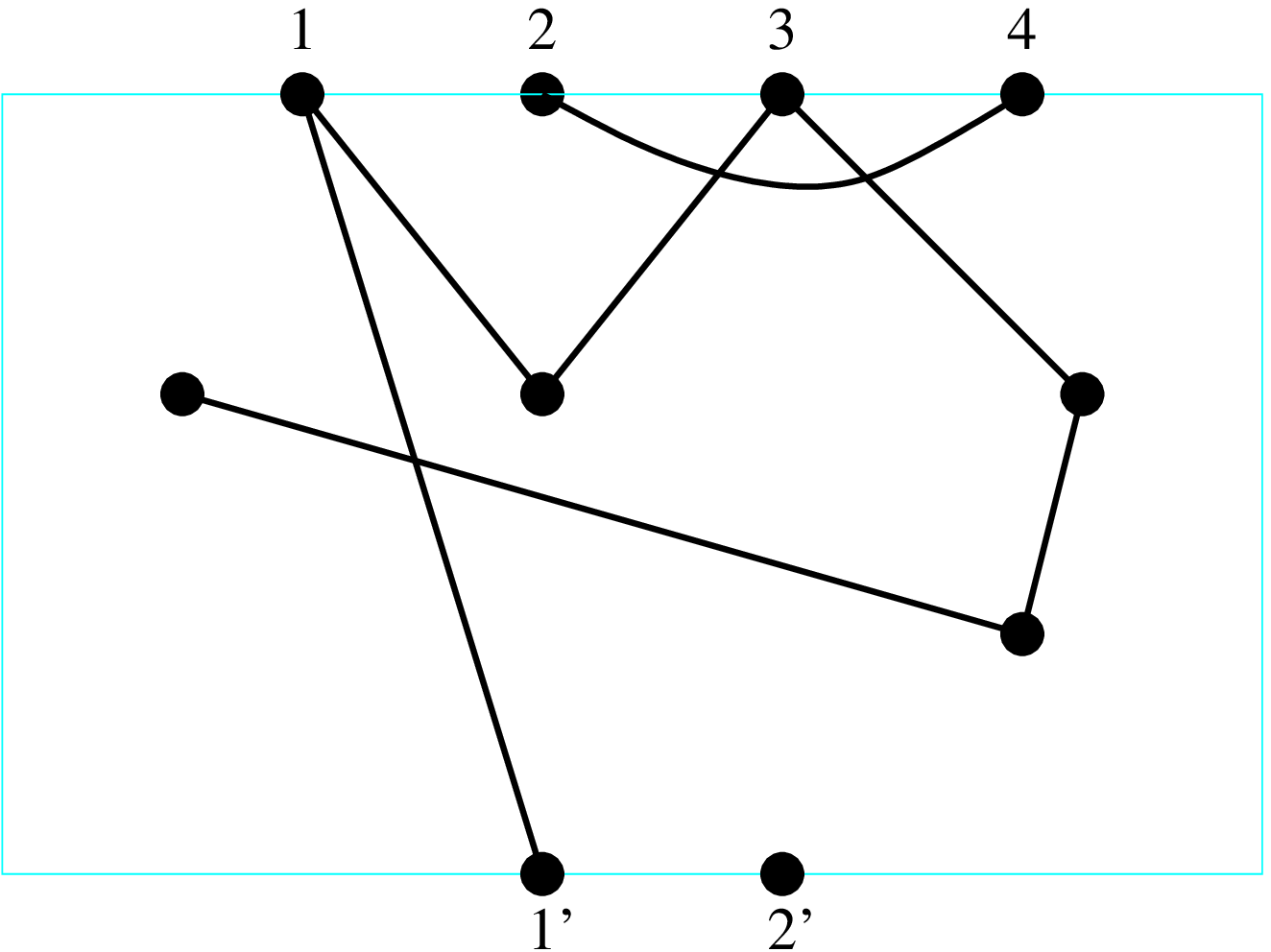}
\]
\caption{Graphs for the partition 
$\{ \{ 1,3,1' \}, \{ 2,4 \} , \{ 2' \}\}$
of the set $\{1,2,3,4,1',2' \}$.
\label{fig:p121}}
\end{figure}

%}}}
%{{{ drawings

\begin{figure}
\[
\includegraphics[width=3.4in]{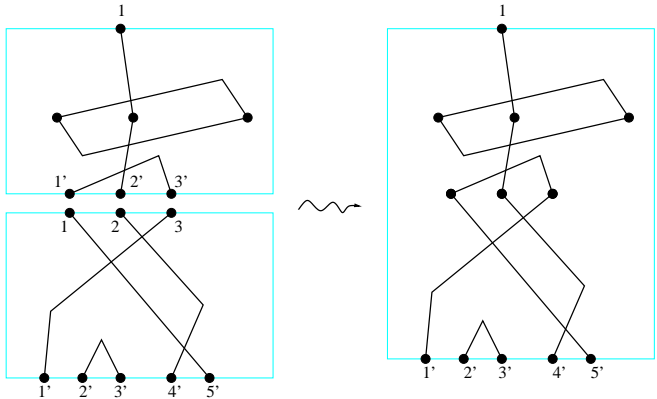}
\]
\caption{Picture stacking composition. \label{fig:p121pex}}
\end{figure}

\mdef \label{de:draw1}
We shall use drawings to represent graphs in a conventional way: 
vertices by points and edges by lines 
(polygonal arcs %, say) 
between vertex points),  %vertices (lines), 
as in Fig.\ref{fig:p121} or~\ref{fig:p121pex}.

A {\em \ppicture} 
$d$ of a graph  %$g$   %$g \in \gr[S]$ 
is thus
(i) a rectangular region $R$ of the plane; 
(ii) an injective map from a finite set into $R$
(hence a finite subset of points identified with vertices); 
and (iii) a subset of $R$ that is the union of lines. 
Line (polygonal arc) crossings are not generally avoidable
(in representing a given graph in this way), 
but we stipulate
{\em `line regularity'}:
that, endpoints apart, lines touch  % crossings occur 
only at points in the interior of straight segments;
and that a line does not touch any vertex point except its endpoints,
or the boundary of $R$ except possibly at its endpoints.

Note: (I)
Regularity ensures that no two graphs have the same picture, and hence
gives us a map `back' from pictures to graphs.
(II) Any finite graph can be represented this way 
(indeed with the vertices in any position, see e.g. \cite{CrowellFox}). 

%}}}
%{{{ External/interior vertices

For $g \in \gr[S]$ one thinks of $S $ as a set of 
{\em `external'} vertices,
and draws them on  the horizontal part of the rectangle boundary.
{\em Interior} vertices ($v \not\in S$)
will generally not need to be explicitly labelled here
(the choice of label 
will be unimportant). 
%%

%}}}
%{{{ partition pics

By (I) and (II),  %Therefore,
via $\Pi_S$, we can use a picture of $g \in \gr[S]$ to represent a partition.
The drawings   in Fig.\ref{fig:p121}  
all represent the same partition,
when regarded as pictures of set partitions of %the set 
$S = \{1,2,3,4,1',2' \}$.
Specifically in each case  %they all represent the partition 
$
\Pi_S(g) \; = \; \{ \{ 1,3,1' \}, \{ 2,4 \} , \{ 2' \}\}   .
$

%}}}
%{{{ vacuum bubble

\mdef
A  {\em vacuum bubble } in $g\in \gr[S]$ is a purely
interior connected component 
\cite{BjorkenDrell65} (cf. Fig.\ref{fig:p121pex}).
The  vacuum bubble {\em number} %\cite{Feynman,Marcolli10}
is
$$
\vloop_S(g) \; = \; \vac_S(\Pi_{}(g)) \; = \; 
  \# \{ \pi \in \Pi_{}(g) \; | \; \pi \cap S = \emptyset   \}
$$

%}}}
%{{{ _n_ _n'_
\gloss{$\underline{n}$}{$\{1,2,...,n \}$}
\gloss{$P(n,m)$}{ $P(\underline{n} \cup \underline{m}')$}
\gloss{$P(n,l,m)$}{subset of $P(n,m)$ of partitions with $l$ prop lines}

\mdef
Let $\underline{n} := \{ 1,2,...,n \}$ and 
 $\underline{n}' := \{ 1',2',...,n' \}$, and so on.
Let 
$$
\Part(n,m) \; := \; \Part(\underline{n} \cup \underline{m}')  .
$$
%}}}

\newcommand{\Nnm}{\underline{n}' \cup \underline{m}}
\newcommand{\Nmn}{\underline{n} \cup \underline{m}'}

\mdef \label{de:nmgraph}
An {\em  $(n,m)$-graph} is an element of %$\Gamma(\Nnm)$
$\gr(n,m) \; := \; \gr[\Nmn]$.
We draw them as in Figures~\ref{fig:p121} and~\ref{fig:p121pex}.
We define $\Pi_{n,m} = \Pi_{\Nmn}$, so 
\beq \label{eq:pismap}\label{eq:Pinm}
\Pi_{n,m} : \gr(n,m) \rightarrow P(n,m)
\eq

%}}}
\medskip
%{{{ P cat

\mdef \label{de:catP0}
Next we recall the partition category $\catP$,  as defined in
\cite[\S7]{Martin94}. 
We first fix a commutative ring $k$ say, and $\delta \in k$.
The  object set in $\catP$ is $\N_0$.
%In $\catB^\circ$, $(n,m)$-arrows are  pair-partitions
The arrow set $\catP(n,m)$ is 
%In $\catP$, an $(n,m)$-arrow is an element of 
the free $k$-module with basis $\Part(n,m)$.
Noting (\ref{eq:pismap})
this means that  elements of  $\catP(n,m)$ could be represented as
formal $k$-linear combinations of $(n,m)$-graphs.
In fact one generalises this slightly.
In $\catP$ an $(n,m)$-graph (as in (\ref{de:nmgraph}))
maps to an element of $kP(n,m)$ via:
\[
\Pi_{\catP} : 
g \stackrel{}{\mapsto} \delta^{\vloop_{n,m}(g)} \Pi_{n,m}(g)
\]
\label{de:PiP}

%}}}
%{{{ fig
%}}}
%{{{ composition and pics for P cat

The composition $p*q$ in $\catP$
can be defined and computed in naive set theory \cite{Martin94}.
But it can also be computed by representing
composed partitions as stacks of corresponding pictures of graphs,
as in Fig.\ref{fig:p121pex}.
First  a composition 
$\circ: \gr(n,m) \times \gr(m,l) \rightarrow \gr(n,l)$
is defined:
 $a \circ b$ is given by stacking pictures of $a$ and $b$ so
that the $m$ vertex sets in each meet and are identified 
as in the Figure.
Then  $p*q = \Pi_\catP (a \circ b)$ for suitable $a,b$.
For example, in case $p = \{ \{ 1,2'\},\{ 1',3'\} \} $ 
in $P(1,3)$ and 
$q =  \{ \{ 1,5' \}, \{ 2,4'\}, \{ 3,1' \}, \{ 2',3' \} \} $
in $P(3,5)$ then the composition 
$\delta p * q$, or more explicitly
\[
\delta
\{ \{ 1,2'\},\{ 1',3'\} \} 
  * \{ \{ 1,5' \}, \{ 2,4'\}, \{ 3,1' \}, \{ 2',3' \} \} 
\;  = \;
\delta
\{ \{ 1,4' \}, \{ 1',5' \}, \{ 2',3' \} \}
\]
can be verified via Fig.\ref{fig:p121pex}.
%}}}
%{{{ composition cont
In general, denoting the stack of pictures by $d | d'$, then 
\beq \label{eq:multidef}
p * q  = \; \Pi_{\catP} (d | d' )
\eq
for any $d,d'$ such that $p =  \Pi_{\catP} (d  )$ and 
$q =  \Pi_{\catP} ( d' )$.
(Given the set theoretic definition of $\catP$ \cite{Martin94} the 
identity (\ref{eq:multidef}) would be a Theorem.
Here we can take it as the definition,
and one must check well-definedness.) 

Remark:
From this perspective the
pictures constitute a mild modification of the 
plane projection of arrows in the tangle category, in which arrows are certain collections of non-intersecting polygonal
arcs in a 3D box (see later, or e.g. \cite{Kassel95}).

%}}}

\subsection{Overview of the paper}
%{{{ preamble
We start with a heuristic overview and summary. 
Later, in order to prove the main Theorems, 
we will give  more formal definitions. 

%}}}
%{{{ invariants

Besides the representation of a set partition $p$ by a graph $g$, the task of
 constructing a  picture $d$ of $p$ contains 
another layer --- the embedding and depiction of graph $g$ in the plane. 
%(specifically within a rectangle $R$, with the  vertices from
% $S$ drawn on the frame; and edges drawn as, say, polygonal arcs). 
Both stages of the representation of set partitions are highly
non-unique. 
%\footnote{DEMOTE?:
However, they lead to some remarkable and useful invariants.
To describe these invariants we will need a little 
preparation.  %notation. 

%}}}
%{{{ alcove

Suppose we have a \emph{picture} $d$ of a 
{\em partition} of this  kind. 
Then each %(let us assume WLOG) 
polygonal arc $l$ of $d$ partitions the rectangle $R$ into 
various parts: one or more 
connected components of $R\setminus  l$; 
and $l$ itself.
Overall, a picture $d$ subdivides $R$ into a number of connected
components, called {\em alcoves}, 
of $R \setminus d$ (regarding $d$ as the union of its
polygonal arcs), together with $d$ itself. 
%}

%}}}
%{{{ distance

\mdef
Given a picture $d$, the {\em distance} 
$\dd_d(x,y)$ from point $x$ to $y$ is
the minimum number of polygonal arc segments 
crossed in
any continuous path from $x$ to $y$. 
Examples (the second picture shows distances to $y$ 
from points in various alcoves):
\[
\includegraphics[width=1.62in]{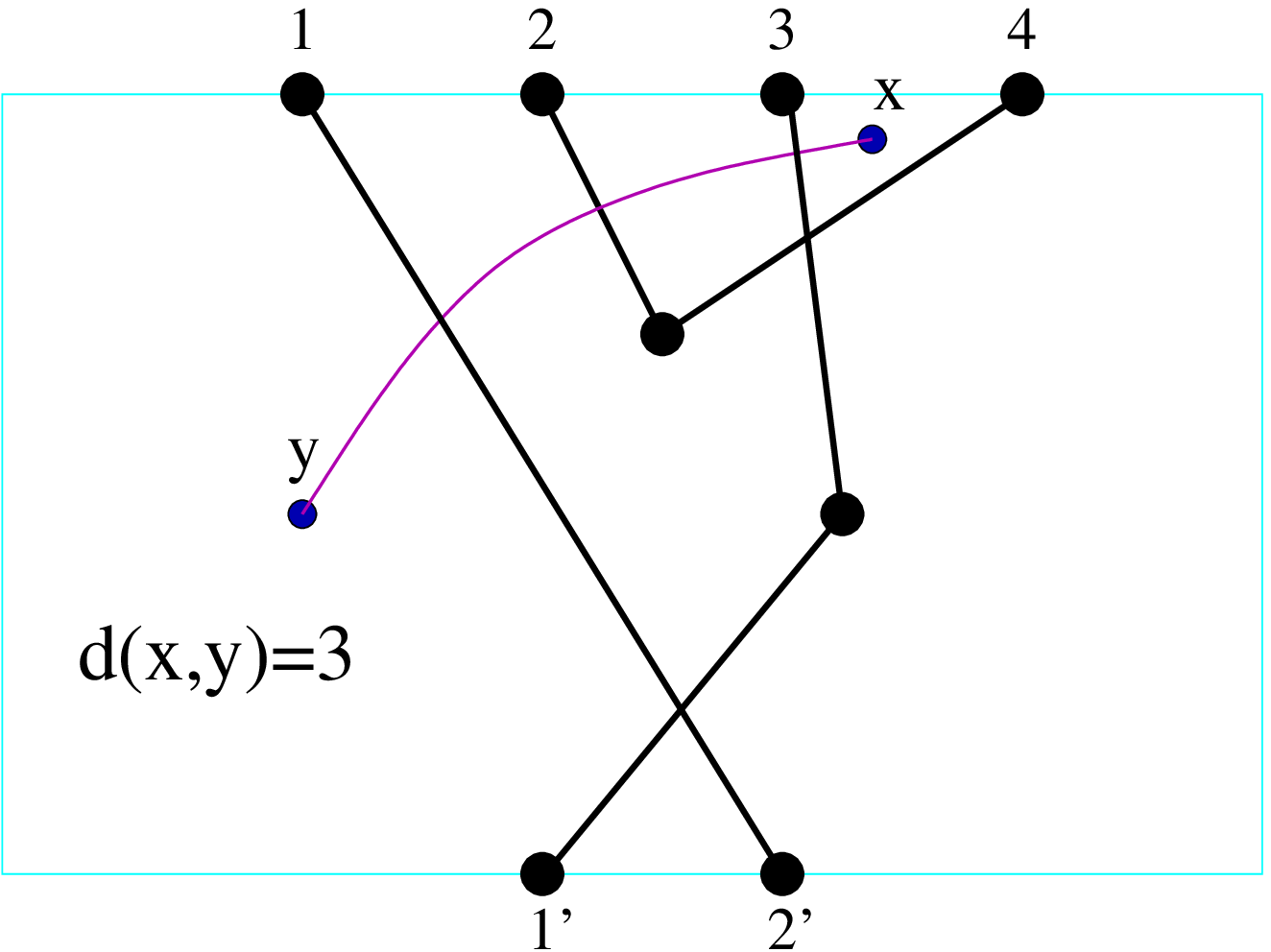}
\qquad
\includegraphics[width=1.62in]{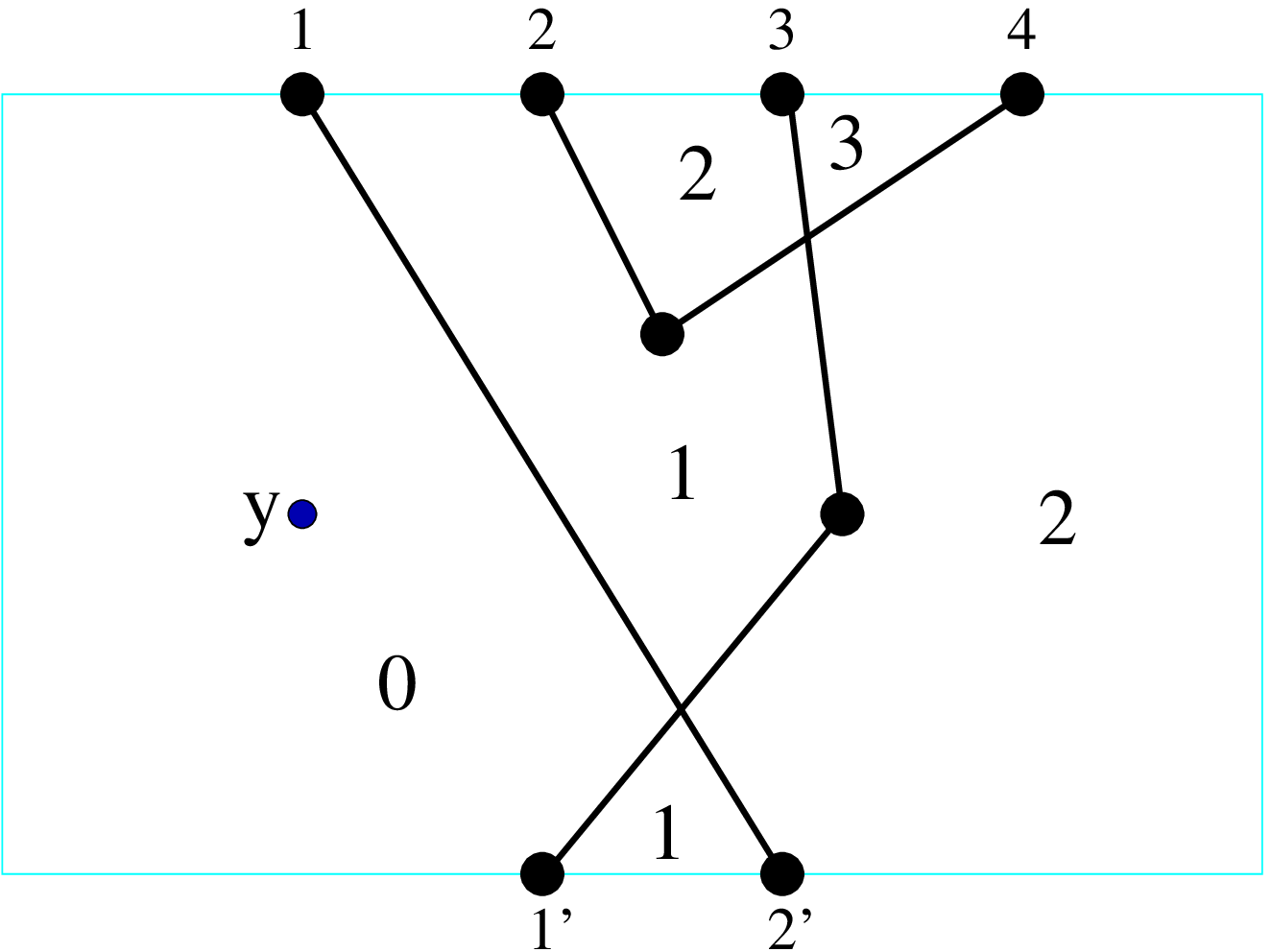}
\]
Note that there is a well-defined distance between a point and an
alcove; or between alcoves.

\mdef
The {\em (left-)height} of a point %/alcove 
in $d$ is defined to be the
distance from the leftmost alcove. 
(By symmetry there is a corresponding notion of right-height.)

Given a picture $d$, a {\em crossing point} is a point 
where two  %not-necessarily distinct 
polygonal arc segments cross.
Note that %Thus 
these points in particular have heights. 
For example  the upper of the two crossing points in the
picture above has height 1, and the other has height 0.

The {\em (left-)height $\hht(d)$ 
of a picture} $d$ with crossing points is defined to be
the maximum (left-)height among the heights of its crossing points.
(We shall say that the left-height of a picture without crossings is
-1.) 

%}}}
%{{{ TL
 
\mdef
Although the picture $d$ of a  %pair 
partition $p$ is non-unique, 
we can ask, for example, if it is possible to draw $p$ without 
arc crossings --- i.e. if among the drawings $d$ of $p$ there is one
without crossings. 
If it is not possible to draw $p$ without crossings, we can ask 
what is the minimum height of  %crossing
picture needed ---
that is, among all the pictures $d$ representing $p$, 
what is the minimum picture height?
We call this minimum the {\em (left-)height of partition $p$}. 

%}}}
%{{{ PLAN
%}}}

%{{{ TL

\newcommand{\catPP}{(\N_0 , k \Part(n,m), *)}

\mdef
Returning to the partition category $\catP  = \catPP$, 
the existence of a Temperley-Lieb %monoid 
subcategory in $\catP$  %$\catB^\circ$ 
(see e.g. \cite[\S6.2]{Martin91}, \cite[\S5.1]{Martin08a})
corresponds to the observation that the product 
$p * p'$ in $\catP$ 
of two partitions of height -1 (i.e. non-crossing) gives rise to a
partition that is again height -1.
Our first main observation is a generalisation of this:
\\
{\em The height of $p*p'$ in %the Brauer monoid category 
$\catP$ 
does not exceed the greater of the heights of $p,p'$.}  Thus:
\\
{\em For each $l \in \Nm \; := \; \{-1,0,1,2,... \}$
there is a subcategory  spanned by  %%consisting of 
the partitions of height at most $l$.}

\medskip

\newcommand{\catBB}{(\N_0 , k J(n,m), *)}

We first prove this result. 
This requires formal definitions of `left-height' and so on, and then
some mildly geometrical arguments.
Write $J(n,m) \subset \Part(n,m)$ for the subset of partitions of 
$\underline{n} \cup \underline{m}'$ into pairs. 
The partitions of this form span the Brauer subcategory:
$\catB  = \catBB$;
and the construction above induces a sequence of subcategories here too.   
The rest of the paper is concerned with the representation theory of
the tower of  
algebras
associated to each of these categories, that is, the  %monoids 
algebras
that are the End-sets in each of these categories.

%}}}

%}}}
\section{Formal definitions and notations}
%{{{ 1
%}}}
%{{{ graph pic

We start with a formal definition of a  {\em \ppicture},  a 
drawing as in (\ref{de:draw1}).
Notation is taken largely from 
Moise \cite{Moise77} and Crowell--Fox \cite{CrowellFox}.

\mdef Given a manifold $M$ we write $\partial M$ for the
manifold-theoretic boundary; and $(M) = M \setminus \partial M$ for the
interior \cite[\S0]{Moise77}. 

\mdef 
A {\em polygonal arc} 
is an embedding $l$ of $[0,1]$ in $\R^3$ consisting of
finitely many straight-line segments. The open arc $(l)$ of $l$ is the
corresponding 
embedding of $(0,1)$.
An 
{\em \pvert ex }
in $l$ is the meeting point of two maximal straight segments.

%}}}
\newcommand{\epsilo}{\epsilon}     %% embedding
\newcommand{\ppic}{(V,\lambda,L)}  %% picture triple (R implicit)
%{{{ take x (polygonal g; regular)

\mdef A {\em polygonal graph} is 
(i) an embedding $\epsilon$ of the vertex
set $V$ of some $g \in \gr(V)$ as points in $\R^3$; and 
(ii) for graph edges $E$ a polygonal embedding
$\epsilon : \sqcup_{e \in E} (0,1) \hookrightarrow \R^3 \setminus
\epsilon(V)$
such that  the closure points of $(0,1)_e$ agree with the endpoints of $e$.

\mdef Note that every $g$ has an embedding for every choice of
$\epsilon:V \hookrightarrow \R^3$. 

\mdef Note that if the edge labels are unimportant then we can recover
the original graph from the map $\epsilon:V$ 
(labeling 
graph 
vertex points) and the image $\epsilon(g)$.
Note well
the distinction between graph vertices $\epsilo(v \in V)$ and
polygonal arc-vertices. 

\mdef \label{de:reg}
Fix a coordinate system on $\R^3$. 
A polygonal graph  $G=\epsilon(g)$ is {\em  regular} (in rectangle
$R \subset \R^2$) if 
\\
(i) the projection $p(x,y,z)=(x,y)$ is injective on vertices; \\
(ii) for $k \in G\setminus \epsilon(V)$ then $|p^{-1}(p(k))| \leq 2$; \\
(iii)    $|p^{-1}(p(k))| =1$ if $k$ an \pvert ex; \\
(iv) for $k \in G$, $p(k) \in p(\epsilon(V))$ implies $k \in \epsilon(V)$; \\
(v) $p(G) \subset R$ and $p(k) \in \partial R$ implies $k \in \epsilon(V)$.

%}}}
%{{{ PICTURE

\mdef \label{de:ppicture}
A {\em \ppicture} is a triple 
$d = (V,\lambda,L)$ 
consisting of
a set $V$, an injective map $\lambda:V \hookrightarrow \R^2$ and a
subset $L \subset \R^2$ such that $\lambda = p \circ \epsilon |_V$ 
for some regular polygonal graph with $g \in \gr(V)$; 
and $L$ is the image $L=p(\epsilo(g))$.
(The datum also includes the containing rectangle $R$, 
but notationally we leave this implicit.)

The point here is that such a $d$, consisting only of labeled points
and a subset of  %%lines in
 $R$, {\em determines} a graph $g$;
and every graph has a picture.
Note that $d$ also determines the set of points where 
$|p^{-1}(p(k))|=2$ 
in (\ref{de:reg})(ii)
--- the set $\Chi(d)$ of {\em crossing points}. 

%}}}
%{{{ ORIENTATION/CONVENTIONS

\mdef  \label{de:framedrawn00}  \label{de:conventions}
Let us consider  \ppicture s  
with  $R$ oriented so that its edges lie in the $x$ and $y$ directions.
If the vertex points on the 
northern (respectively southern) edge of 
$R$ are not labelled explicitly then they may be understood to be
labeled $1,2,...$ (respectively $1',2',...$)
in the natural order from left to right.

In particular such a 
{\em frame-drawn}
picture {\em without any vertex point labels}
determines a graph in some 
$\gr[\underline{n}\cup\underline{m}']$ 
up to labelling of the other `interior' vertices. 

We identify  pictures differing only by an overall vertical shift.
Given our vertex labelling convention above we could also identify under 
horizontal shifts, but the horizontal coordinate will be a useful tool in
proofs later, so we keep it for now.

%}}}

\subsection{Stack composition of pictures}
%{{{ stacking
\newcommand{\danger}{\includegraphics[width=.21in]{dangerousbend.eps}
}

%{{{ category pic conventions

\newcommand{\nn}{{\tt n}}
\newcommand{\sss}{{\tt s}}
\newcommand{\hh}{{\tt h}_0}

\newcommand{\fdp}{frame drawn \ppicture}
\newcommand{\chainpicture}{chain \ppicture}

Here we follow the usual construction of `diagram categories' 
(e.g. as in \cite[\S7]{Martin94}), 
but take care to emphasise specific
geometrical features that we will need later. 
(See also e.g. \cite{ReshetikhinTuraev90,BaezDolan95,AlvarezMartin06}.)

%}}}

\newcommand{\mpt}{marked point}  %% alt name for vertex

%{{{ defns (marked pt)

\mdef \label{de:hh}\label{de:mapt}
Given $d=\ppic$ 
in rectangle $R$ 
write $\nn(d)$ for the subset of $\R$ giving the intersection of $L$
with the northern edge of $R$  % the frame 
(thus 
by (\ref{de:reg})(v)
the collection of $x$-values of `northern'
exterior vertex points, or `\mpt s'). 
Write $\sss(d)$ for the corresponding southern set.  

Write $\hh(a,b)$ for the class of pictures $d$ with  $\nn(d)=a$
and  $\sss(d)=b$,
and $L$ not intersecting the  two vertical edges
of the containing rectangle.

%}}}
%{{{ stacking

\mdef
Note that for $d \in \hh(a,b)$ there is an essentially identical 
picture with $R$ {\em wider}. Thus any two such pictures may be taken to have
the same (unspecified, finite) interval of $\R$ as their
northern edges,  and southern edges. 
The juxtaposition of rectangular
intervals, $R,R'$ say, by
{\em vertical stacking} 
of $R$ over $R'$ thus produces a rectangular interval, 
denoted $R|R'$. 
This is almost a disjoint union, except that  %$\sss(R,S)$ 
the southern edge of $R$ 
is identified with  %$\nn(R',S')$
the northern edge of $R'$.

Given a pair of \ppicture s $d$ and $d'$, 
stack $R|R'$ induces a corresponding pair of subsets 
$\lambda(V) | \lambda'(V')$  %$S_0|S_0'$ 
and $L | L'$ in
the obvious way.
For example see Fig.\ref{fig:p121pex}.
%}}}
%{{{ cat

{\mpr{ \label{pseudopicturecategory}
The stack juxtaposition of a picture $d$ in  $\hh(a,b)$ 
over a picture $d'$ in 
 $\hh(b,c)$ 
defines 
a picture $d|d'$ in  $\hh(a,c)$.
}}

\noindent
Proof: %\danger
As noted, 
the stack $R|R'$ induces a corresponding pair of subsets 
$\lambda(V) | \lambda'(V')$  %$S_0|S_0'$ 
and $L | L'$. % in the obvious way. 
The former is a union of finite point sets which 
clearly agrees with $a$ and $c$ on the relevant edges
of $R|R'$. 
The latter is a union of lines, and the only new meetings are at the
marked points in $b$ (as it were). These are now interior marked
points.
Conditions (\ref{de:reg})(i-v) hold by construction.
\Qed

%}}}
%{{{ cat2

\newcommand{\hhcat}{\hh}  %% picture category

\mdef \label{pr:whax}
Allowing rectangles of zero vertical extent in any
$\hh(a,a)$ 
allows for an 
identity element $1_a$ of stack composition in $\hh(a,a)$. % in general).
Write $\hhcat$ for the `picture category'.

\newcommand{\pisse}{\pi_{e}} %% projection to partitions from fdpics
\newcommand{\pissep}{\pi_{p}} %% projection to k-partitions from fdpics

{\mpr{
For finite %discrete 
subsets $a$ and $b$ of $\mathbb{R}$, there is a surjection
\eql(eq:surjp)
\pisse : \hh(a,b) \rightarrow P(|a|, |b| )
\eq
given by counting the elements of $a$ (resp. $b$) from left to right
and using $\Pi_{n,m}$ from (\ref{eq:Pinm}).
\Qed
}}

%}}}
%{{{ catP

{\mpr{\label{pr:whay}
Fix a commutative ring $k$ and $\delta \in k$, and let 
 $\pissep$ denote the generalisation of $\pisse$ corresponding to
$\Pi_\catP$ from (\ref{de:PiP}).
Let $d \in \hh(a,b)$, $d' \in \hh(b,c)$. 
Then 
$\pissep (d|d') =  \; \pissep (d) * \pissep (d')$,
where 
the product on the  %%RHS
right is in the partition category $\catP = (\N_0 , k P(n,m), *)$
\cite{Martin94}. 
}}

We can take this as a {\em definition} of $\catP$ cf. (\ref{de:catP0}). 
For a proof of equivalence to the original definition see e.g.
\cite{Martin94}. 
In outline, one checks that the stack composition implements the
transitive closure condition \cite[\S 7]{Martin94}. 
\Qed

%}}}
%{{{ flip pic, flop partition

\newcommand{\flip}{*} % flip for pictures
\newcommand{\flop}{*} % flip for partitions

\medskip

\mdef 
If $d \in \hh(a,b)$ is a \ppicture\ as above, 
let $d^\flip \in \hh(b,a)$ denote the \ppicture\ obtained
by flipping $d$ top-to-bottom. 

Let $p \in \Part(n,m)$. 
Write $p^\flop$ for the element of $\Part(m,n)$ obtained by swapping primed
and unprimed elements of the underlying set.

Note that if $d$ is a \ppicture\ of $p \in P(n,m)$ then 
$d^\flip$ is a \ppicture\ of $p^\flop \in P(m,n)$. Furthermore, 
this $\flop$ 
is a contravariant functor between the corresponding partition categories.

\mdef \label{de:otimespic} \label{de:catmonoid}
Note that for any \ppicture\ in the category 
$\hhcat$
with distinct northern and southern
edge we can vertically rescale to arbitrary finite separation of these edges. 
Thus we can make any two \ppicture s have the same separation. 
For two such pictures $d,d'$ there is a \ppicture\ $d \otimes d'$
obtained by side-by-side juxtaposition. 

%}}}
%}}}

%{{{ chain

%{{{ frame drawn

\mdef
We call  a \ppicture\ a {\em \chainpicture} if every exterior marked point
(as in (\ref{de:mapt}))
is an endpoint of precisely one line, and every interior marked point
is an endpoint of precisely two lines
(e.g. as in Fig.\ref{fig:p121pex}). 
Write $\hh^2(a,b) \subset \hh(a,b)$ for the subset of \chainpicture s. 
Note that every  $p \in J(n,m)$ has a 
\chainpicture.
We have the following.

%}}}
%{{{ Brauer cat / TL cat

{\mlem{
\label{pr:whaz}
The stack composition (\ref{pseudopicturecategory}) 
closes on \chainpicture s.  
This gives a subcategory of the category in (\ref{pr:whax}).
The corresponding $\pissep -$ quotient category 
(as in (\ref{pr:whay})) is the Brauer category 
$\catB = (\N_0 , k J(n,m), *)$. 
\Qed
}}

\newcommand{\catT}{{\mathcal T}} %% TL cat

\mdef \label{de:catPBTT}
A pair partition is {\em plane} if it has a \fdp\ %frame drawn picture 
(as in (\ref{de:conventions}))
without crossings of lines. We write $T(n,m)$ for the subset of plane
pair partitions  %% (Temperley--Lieb 
(TL partitions)
and $\catT =(\N_0 , k T(n,m), *)$ for the corresponding subcategory
of $\catB$. 

%}}}
%}}}

\subsection{Paths and the height of a picture/a partition}
%{{{ path
%{{{ alcove

{\mrem{
Fix a rectangle $R$. 
Each non-self-crossing
line $l$ with exterior endpoints in  $R$
can be considered to define a  
separation of $R$ into
three parts --- the component of $R \setminus l$ containing the left
edge;
the other component of  $R \setminus l$; and $l$. 
}}

\mdef
An {\em alcove} of picture $\ppic$ is a connected component of $R \setminus L$. 

\medskip

%}}}
%{{{ left-height 1
\newcommand{\simple}{simple}  %% simple point

\mdef \label{de:path}
A {\em \simple\ point} of a subset $U$ of $\R^2$ is a point having a
neighbourhood in $U$ that is an open straight segment.
(For example in some picture $d=\ppic$, with $L=p(\epsilo(g))$, 
the {\em non}-\simple\ points of $L\setminus \lambda(V)$
are the
\pvert ices and crossing points.) 

Given a picture $\ppic$, a {\em path in $\ppic$} 
is a %smooth path 
line $l$
in $R$ such that 
every $k \in (l) \cap L$ is a \simple\ point of $(l)$ and also 
a \simple\ point of 
$L \setminus \lambda(V)$.

Thus a path $l$ in $\ppic$ has a well-defined number of line
crossings,
$\dd_L(l)$.

%}}}
%{{{ left-height 2

{\mlem{\label{avoidl}
Given a path $l$ in picture $d$ connecting points $x,y\in R$ 
and a distinct point $z\in (l)$, there is a path $l'$ connecting $x,y$ 
that does not contain $z$.}}

\proof Since $z\in (l)$ it  has a neigbourhood either 
containing only a segment of $l$; or only a crossing  of $l$ 
with a straight segment of $L$.
If we modify the path inside the 
neighbourhood by a small polygonal detour  %half-circle 
then the modification is a path and  does not contain $z$. \Qed     

{\mlem{
For each picture $d$ and $x,y \in R$ there is a path in $d$
from $x$ to $y$.
}}

\proof
Draw a small straight line $l_1$  from $x$ to a point $x'$ in 
an  adjacent alcove,
choosing $x'$ so  
 that the tangent of the straight line $x'-y$ is
not in the finite set of tangents of  
segments of lines of $L$; 
and the line does not contain
any element of the finite set of crossing points of $d$. 
Then $x-x'-y$ is a path.
\Qed

%}}}
%{{{ x-ht

\mdef
Given a picture $d=\ppic$, and  points $r,x$ in $R$,
the {\em $x$-height} 
$\dd_d(r,x) = \dd_L(r,x)$
of $r$  %% a point $r$ in $R$ %an alcove
 is the minimum  %number of line crossings
value of  $\dd_L(l)$ 
over paths $l$ in $d$ from $r$ to $x$.

We suppose that
$L$ does not intersect the left edge $R_L$ of $R$. 
The {\em (left)-height} 
$\; \hht_L(r) = \dd_L(r,x)$
in case $x$ is any point on $R_L$. 
(Note that this is well-defined.) 

The {\em (left)-height} of an alcove $A$ 
is the left-height of a point in $A$. %
See Fig.\ref{fig:wild1} for examples. 

%}}}
%{{{ Fig fig:wild1

\newcommand{\smoothy}{(Remark: By \cite[\S6]{Moise77} 
piecewise linear and piecewise smooth lines are effectively
  indistinguishable as far as physically drawn figures are concerned.)
}

\begin{figure}
\[
d_1 = \includegraphics[width=2in]{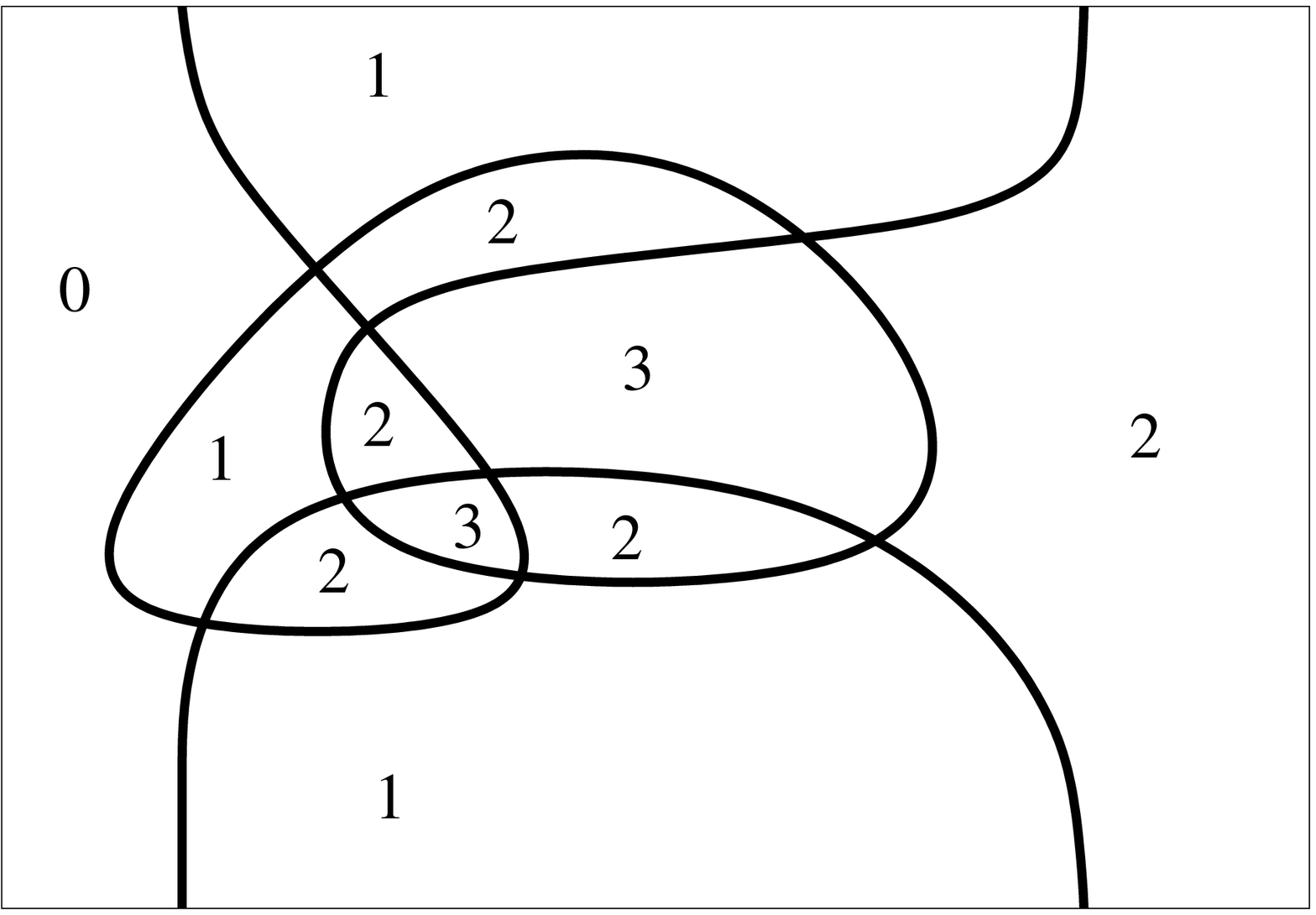}
\qquad
d_2 = \includegraphics[width=2in]{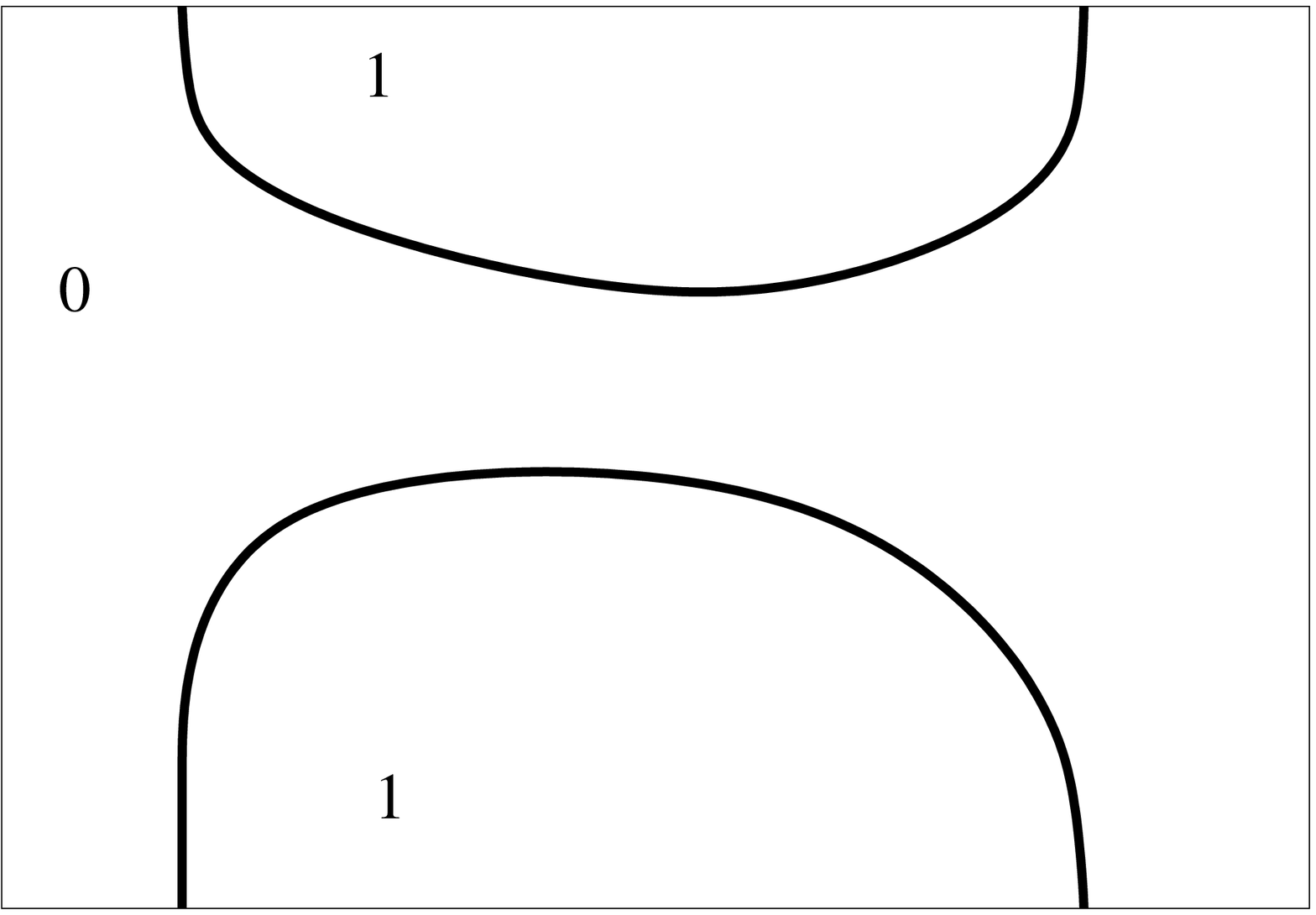}
\]
\caption{\label{fig:wild1} Example pictures with left-heights of
  alcoves. \smoothy
\label{fig:squigg}
}
\end{figure}

%}}}
%{{{ left-height 3

\mdef \label{de:picht}
Given a picture $d=\ppic$, 
recall     %%write 
$\Chi(d)$  %%for 
is the set of crossing points of lines in $d$
(recall from (\ref{de:draw1}) that, vertex points aside,
lines in $d$ only meet at crossing points).
The {\em left-height } $\hht(d)$ is the
greatest of the left-heights of the points $x \in \Chi(d)$;
or is 
defined to be 
-1 if there are no crossings.

For example, $d_1$ in Fig.\ref{fig:wild1} has left-height 2;
and $d_2$ has left-height -1.

\mdef \label{de:ppht}
Finally we say that a  partition  $p \in P(n,m)$ has left-height 
$\hht(p) = l$ if 
it has a \fdp\ %frame drawn picture 
of left-height $l$, and no such \ppicture\ 
of lower left-height.
For example, both \ppicture s in Fig.\ref{fig:wild1} give the same
partition $p$, %. The left-height of this partition is -1.
so $\hht(p)=-1$.

Since every $p$ has a picture, 
it will now be clear that $\hht(p)$ defines a function 
$$
\hht : P(n,m) \rightarrow \{ -1,0,1,2,... \}  .
$$ 

A path realising the left-height of a point 
in  a \ppicture\ is called a {\em low-height  path}. 
A \ppicture\ realising the left-height of a partition is called a 
{\em low-height \ppicture}. 

%}}}
%{{{ notation

\medskip

\newcommand{\Pel}[1]{P_{ #1}} 
\newcommand{\Jel}[1]{J_{ #1}} 
\newcommand{\Pl}[1]{P_{\leq #1}} 
\newcommand{\Jl}[1]{J_{\leq #1}} 

\mdef
Define $P_l(n,m)$ as the subset of partitions in $P(n,m)$ of
left-height $l$, and 
$$
\Pl{ l} (n,m) \; = \bigcup_{j \leq l} P_j(n,m)
.$$
Define $J_l(n,m)$ as the 
corresponding subset of $J(n,m)$, and $\Jl{l}(n,m)$ analogously.
%subset of pair partitions in $J(n,m)$ of left-height $l$, and 

%}}}

%{{{ new

\mdef \label{rem:maxheight} Remark: 
Observe that for $p \in J(n,n)$, $\hht(p) \leq n-2$.
Hence,  
in particular, $\Jl{r}(n,n) = \Jl{n-2}(n,n)$ for any $n-2 < r < \infty$.

%}}}

%{{{ ex

{\mex{
Here we give the  $\Jel{l}(3,3)$ subsets of $J(3,3)$.
Each element is represented by 
a low-height picture  %which happen to  
(of course,  other pictures could have been chosen instead). 
Note that it is a Proposition that a given picture is low-height.
One should keep in mind that the elements 
of $\Jel{l}(3,3)$ are pair partitions, not pictures!
\[ J_{-1}(3,3) = \{
\includegraphics[width=.5in]{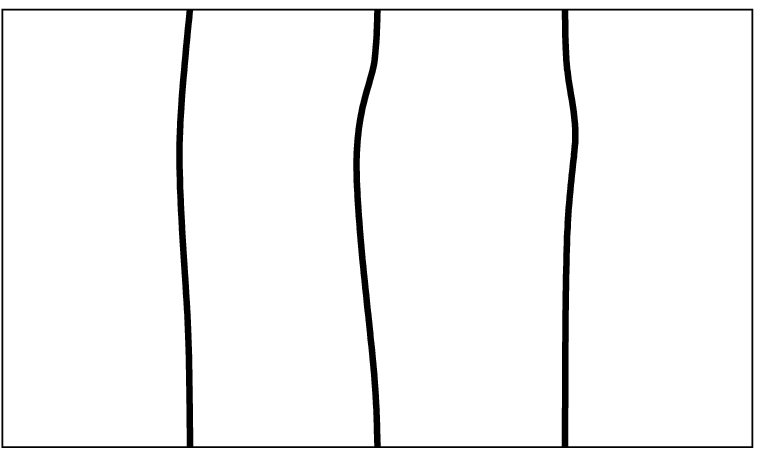},
\includegraphics[width=.5in]{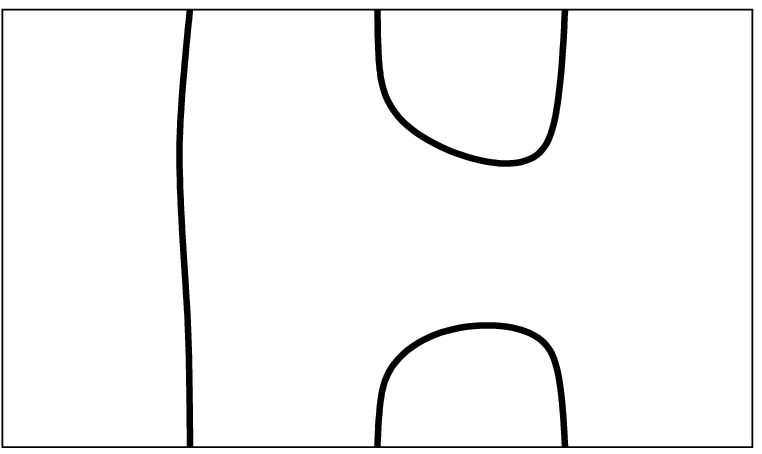},
\includegraphics[width=.5in]{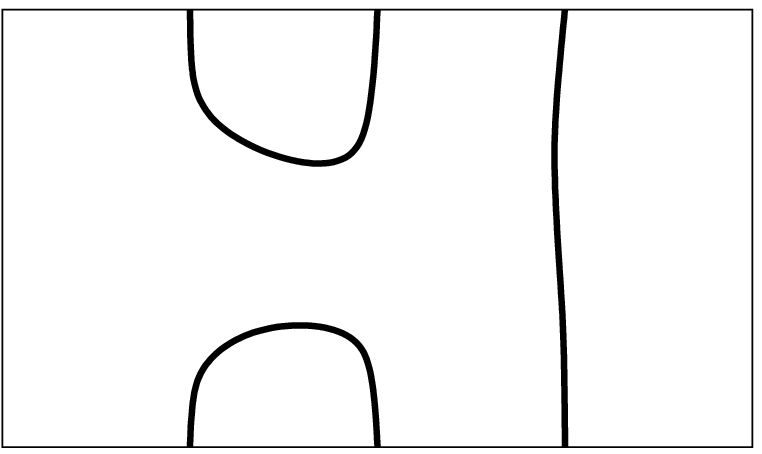},
\includegraphics[width=.5in]{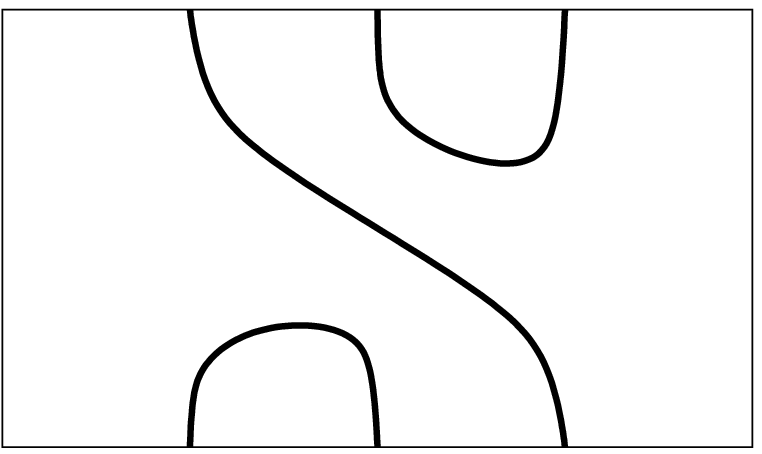},
\includegraphics[width=.5in]{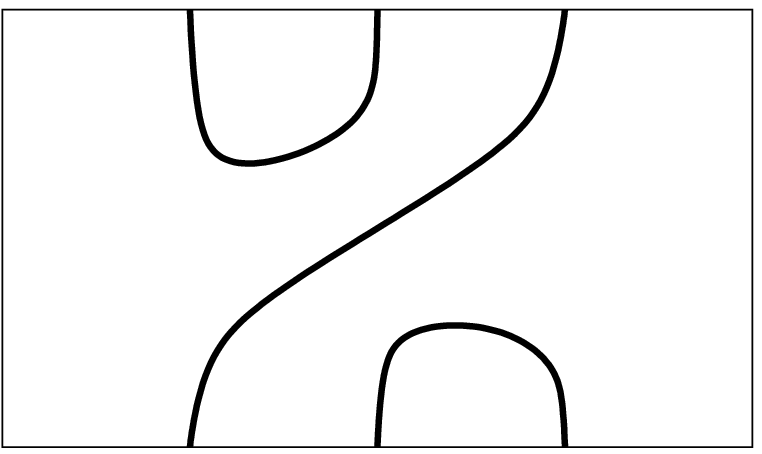}
\}
\]
\[ J_{0}(3,3) = \{
\includegraphics[width=.5in]{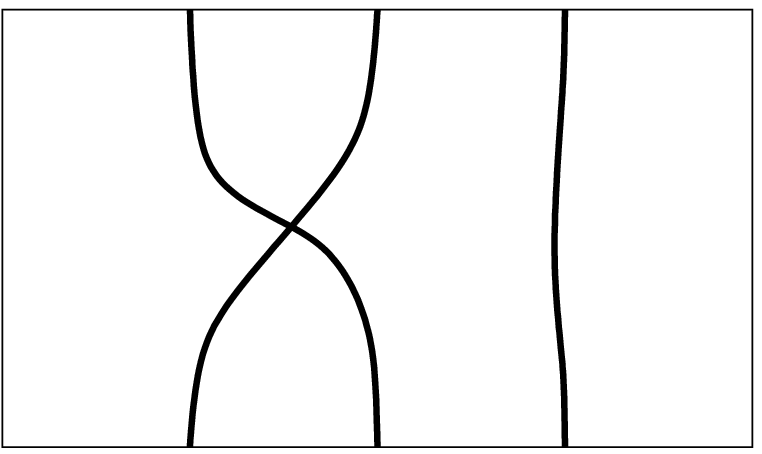},
\includegraphics[width=.5in]{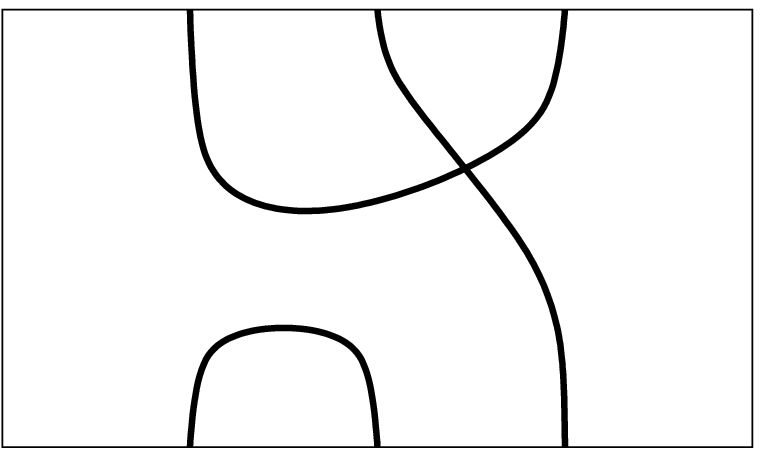},
\includegraphics[width=.5in]{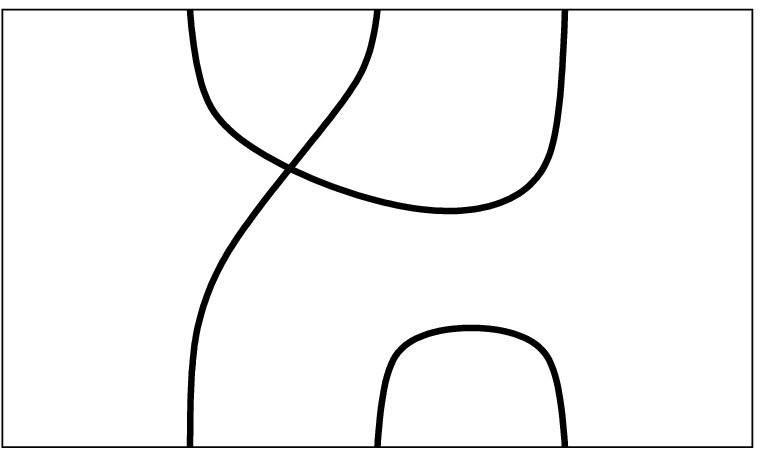},
\includegraphics[width=.5in]{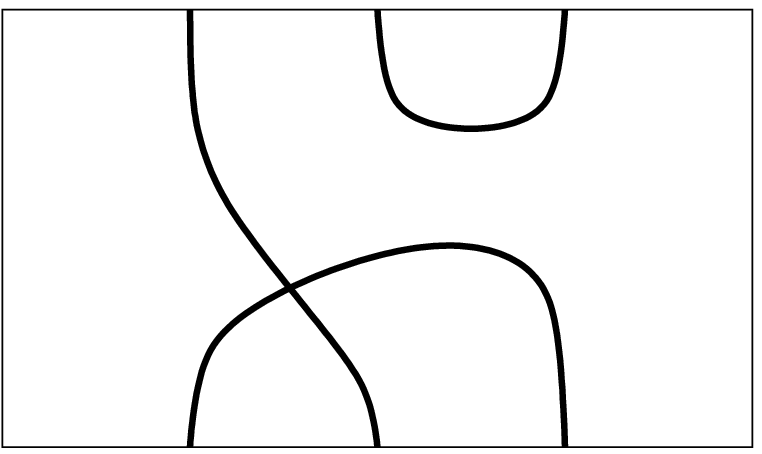},
\includegraphics[width=.5in]{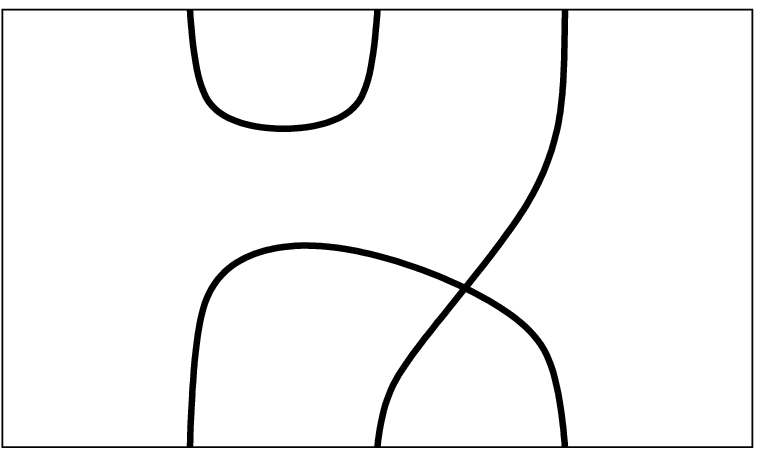},
\includegraphics[width=.5in]{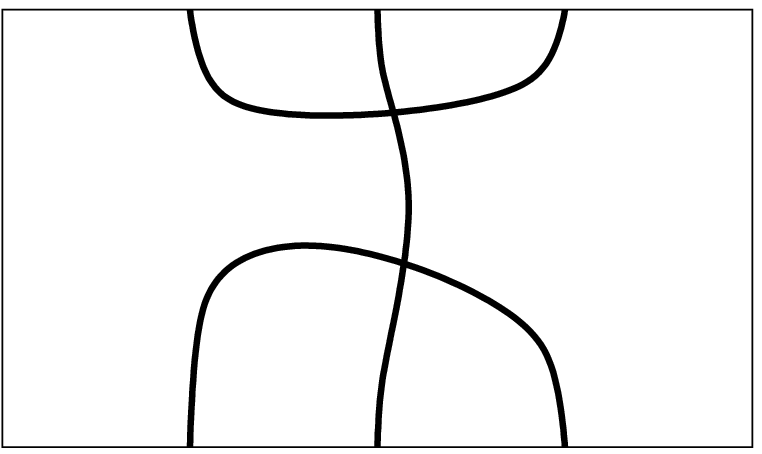}
\}
\]
\[ J_{1}(3,3) = \{
\includegraphics[width=.5in]{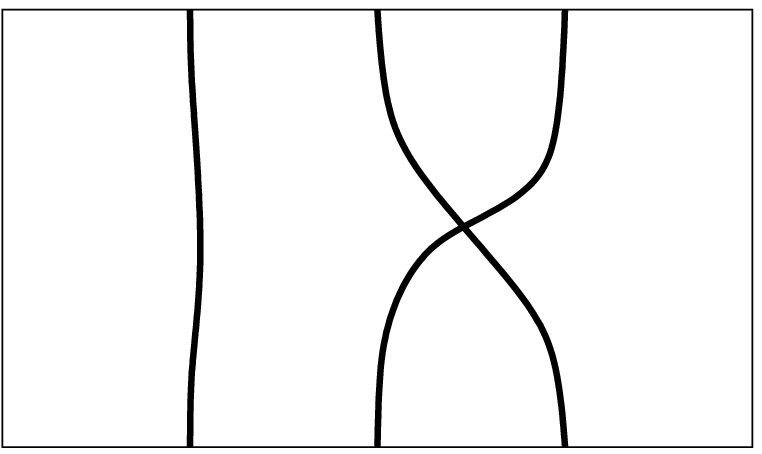},
\includegraphics[width=.5in]{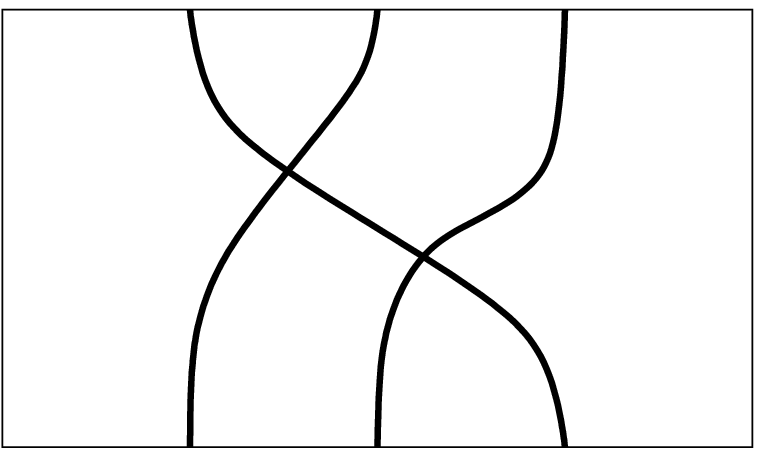},
\includegraphics[width=.5in]{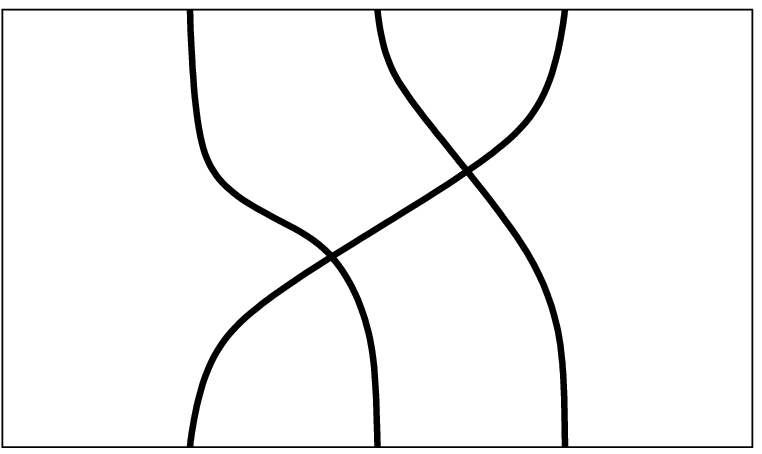},
\includegraphics[width=.5in]{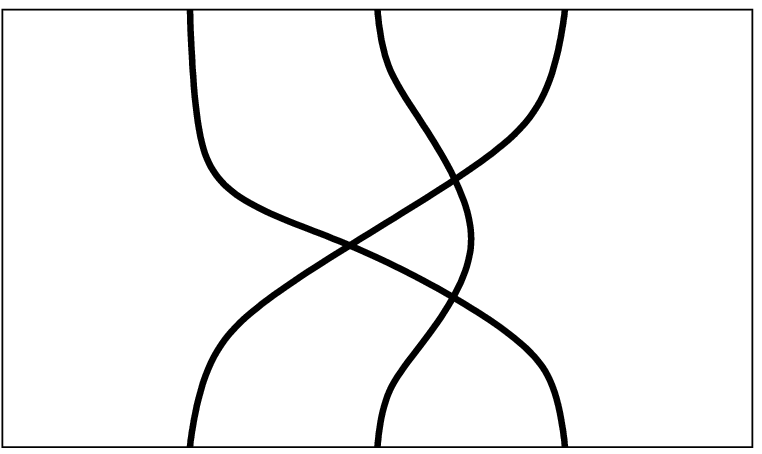}
\}
\]
More generally the Brauer algebra identity element $1_n \in J_{-1}(n,n)$;
the symmetric group Coxeter generator $\sigma_1 \in J_0(n,n)$;
$\sigma_i \in J_{i-1}(n,n)$;  
\[
J_{1}(9,9) \ni 
\raisebox{-.2in}{
\includegraphics[width=1.5in]{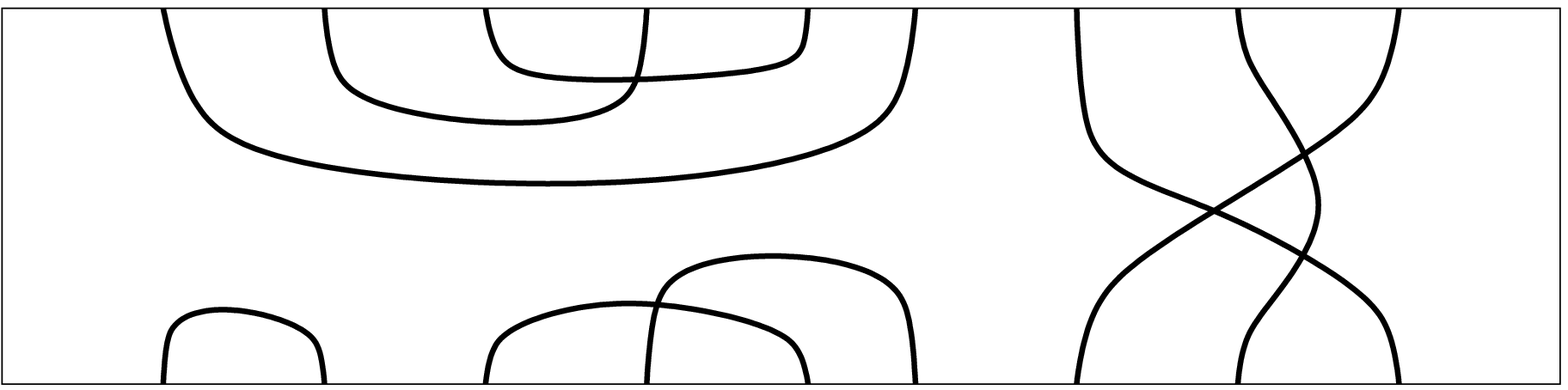}}
;\qquad 
\qquad %\mbox{ and }
J_{2}(10,10) \ni 
\raisebox{-.2in}{
\includegraphics[width=1.8in]{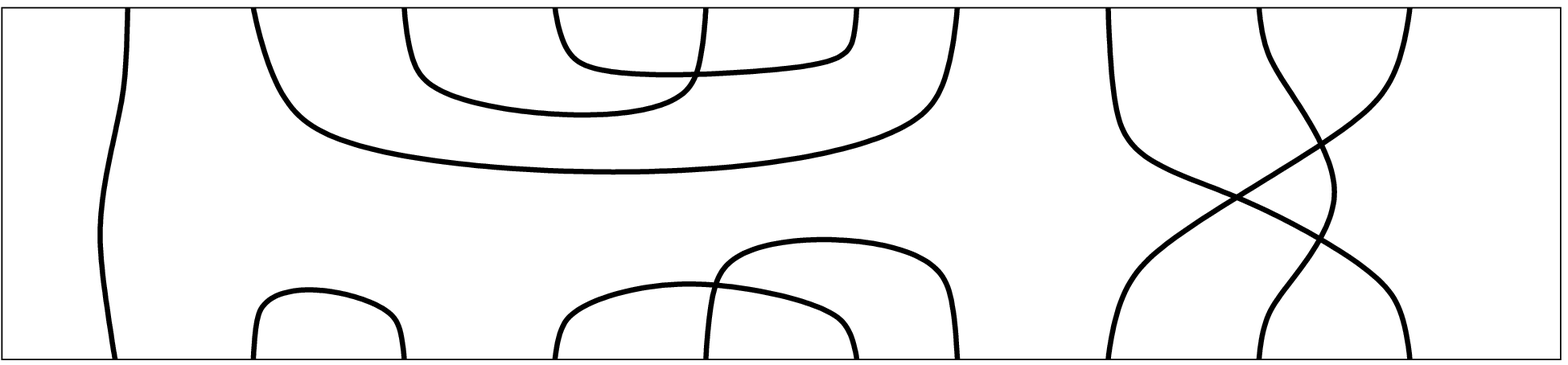}}
\]
\[
J_{2}(10,10) \ni 
\raisebox{-.2in}{
\includegraphics[width=1.8in]{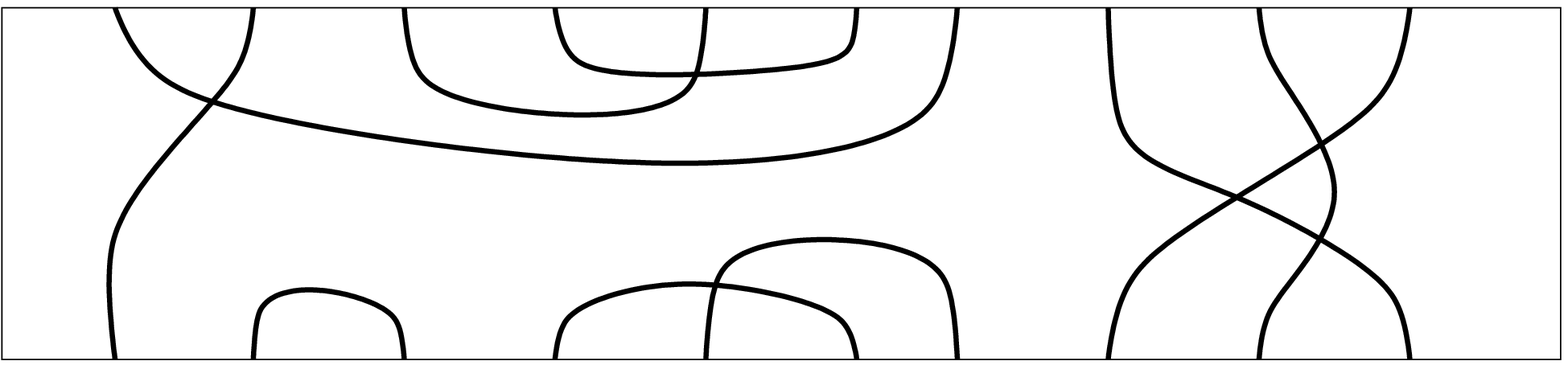}}
\]
}}

%}}}
%{{{ remove

{\mlem{ \label{lem:deletion} 
Removing part or all of a line from a \ppicture\ 
cannot %make the picture height higher.
produce a \ppicture\ with higher height.
}}

\proof
The number of crossings of a path cannot be increased by removing a
line. 
\Qed

\mdef In particular, if a line has a self-crossing then we can 
`short-circuit' the path without increasing the height,
or changing the partition. The self-crossing point becomes an \pvert ex 
with a regular neighbourhood, so the regularity of the picture remains 
to hold.  Thus for each low-height picture there is a low-height picture without
line self-crossings. 

%}}}
%}}}

%}}}
\section{Algebraic structures over $\Jl{l}(n,n)$}
%{{{ CONSTRUCTION
\newcommand{\Jln}{J_{l,n}}  %% our main algebra
\gloss{$\Jln$}{algebra}
\newcommand{\Jlnn}[1]{J_{l,#1}}  %% our main algebra - n variable
\newcommand{\Jlnkk}{J_{l,n}^\kk}  %% our main algebra over \kk

%{{{ alg
Define the commutative ring $\kk = \Z[\delta]$
and $\Bbbk = \Z[\delta, \delta^{-1}]$.

\mdef
Recall (e.g. from \cite{Brauer37} or (\ref{pr:whaz}))
that the %algebra 
multiplication in the Brauer  
$\kk$-algebra
$B_n^\kk = \kk J(n,n)$
may be defined via
vertical juxtaposition of representative diagrams.

\mdef \label{de:Jln}
Define $\Jlnkk$ as the $\kk$-subalgebra of 
$B_n^\kk$
generated by $\Jl{ l}(n,n)$. 
For $k$ a fixed commutative ring and $\delta_c \in k$ we write
$J_{l,n} = J_{l,n}(\delta_c)$ for the base-change:
\beq \label{eq:base}
J_{l,n}(\delta_c) \; := \; k \otimes_{\kk} \Jlnkk
\eq
(i.e., regarding $k$ as a $\kk$-algebra in which $\delta$ acts as
$\delta_c$). 
For example $J_{l,n}(1)$ is the monoid $k$-algebra for the
submonoid of the Brauer monoid associated to $B_n(1)$
\cite{Mazorchuk95}. 

%}}}
\newcommand{\caJ}{{\mathcal J}}
%{{{ cat

\mdef \label{de:caJ}
Define $\caJ_{l}$ as the subcategory of 
%monoid $k$-algebra $J_n = B_n(1) \subseteq P_n(1)$
$\catB$
generated by the sets $\Jl{l}(n,m)$, over all 
$n,m \in \N_0$  %%positive integers $n$ and $m$
(that is, the smallest $\kk$-linear subcategory such that the collection of arrows 
$\caJ_l(n,m)$ contains $\Jl{l}(n,m)$ for each $n$, $m$). 

%}}}
%{{{ comments and main
\medskip

Note that the smallest possible height for a crossing is 0, so
$J_{-1}(n,m)$ is the subset with no crossings, that is $J_{-1}(n,m)=
T(n,m)$. It will be clear then that $J^\kk_{-1,n}$ is the Temperley--Lieb
algebra. 
Complementarily, $J^\kk_{n-2,n}$ ($=J^\kk_{\infty,n}$, by Remark
\ref{rem:maxheight}) is the 
Brauer %monoid 
algebra. 
Next we will show (our first main Theorem)
that the `interpolation' is proper in the following sense: 

%}}}
%{{{ main theorem

\newcommand{\caJJ}{(\N_0, \kk \Jl{l}(n,m), *)}

{\mth{ \label{pr:main} \label{th:main}
(I) The sets $\Jl{l}(n,m)$ form a basis for the $\kk$-linear category
    $\caJ_{l}$. 
That is, $\caJ_l = \caJJ$.   $\;$ 
(II) The set $\Jl{l}(n,n)$ is a $\kk$-basis for $\Jlnkk$.
 }}

\proof{ 
(I) Recall from definition (\ref{de:caJ}) 
that the $n,m$-arrow set in $\caJ_{l}$ is generated by $\Jl{l}(n,m)$.
Note that  $\Jl{l}(n,m)$ is linearly independent over $\kk$
in $\caJ_l$, 
as it is linearly independent in the Brauer category $\catB$.
It is thus enough to show that 
$\Jl{l}(n,m) \times \Jl{l}(m,j)$   %%is closed under 
maps into $\kk\Jl{l}(n,j)$ under the Brauer category product,
and that the sets of arrows of form $\Jl{l}(n,m)$
are therefore suitably spanning.
%the algebra product in $J_n$ and hence spans $J_{l,n}$. 
For this, we need to
show that for every pair of pair-partitions 
$(p_1,p_2) \in \Jl{l}(n,m) \times \Jl{l}(m,j)$, 
determining a 
partition $p_3 \in J(n,j)$, %then 
we have $p_3 \in \Jl{l}(n,j)$.
By (\ref{de:ppht}) 
the pair-partitions $p_1, p_2$ have 
composable minimum-height pictures, denoted by $d_1,d_2$, respectively. 
By (\ref{pr:whay}) 
their vertical juxtaposition $d_1|d_2$ 
gives $p_3$, so  %we are done if we can
it is sufficient to show that  $\hht(d_1|d_2)  \leq l$.

Observe that by construction
the set of crossing points of $d_1|d_2$ is precisely the
disjoint union of those of $d_1$ and $d_2$.
Now, also by construction a low-height path from any point $x$
in $d_1$ remains a (not necessarily of low height) path in
$d_1|d_2$. See the path from $x$ in the figure below for example. 
\[
\raisebox{.64in}{\includegraphics[width=3.25cm]{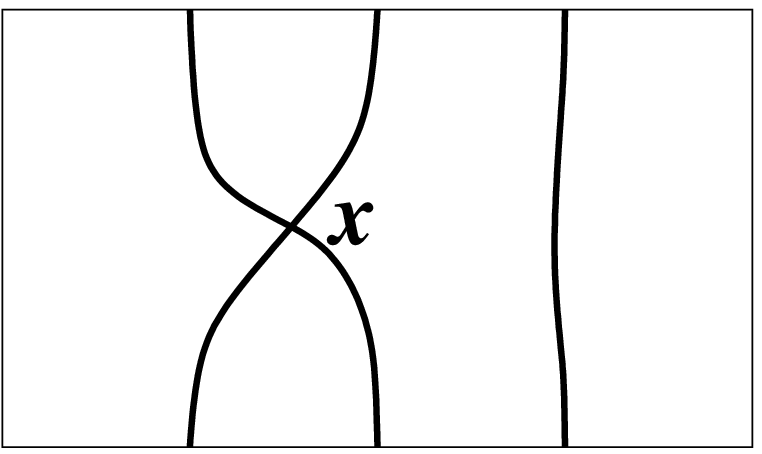}} 
\hspace{.25in}
\raisebox{.64in}{\includegraphics[width=3.25cm]{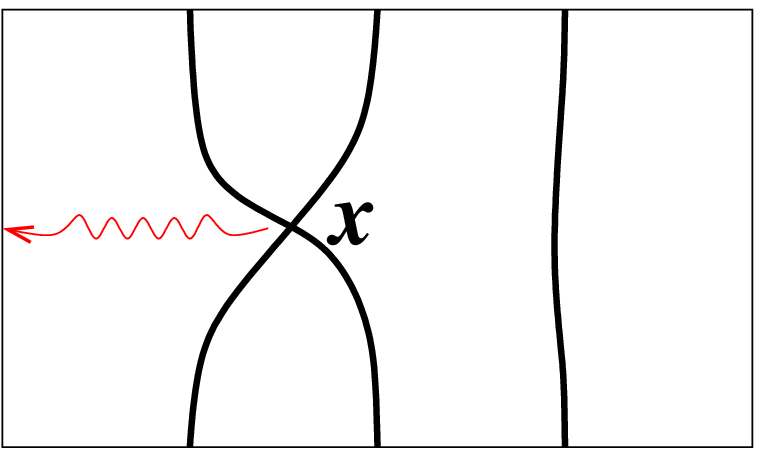}} 
\hspace{1.25in}
\includegraphics[width=3.25cm]{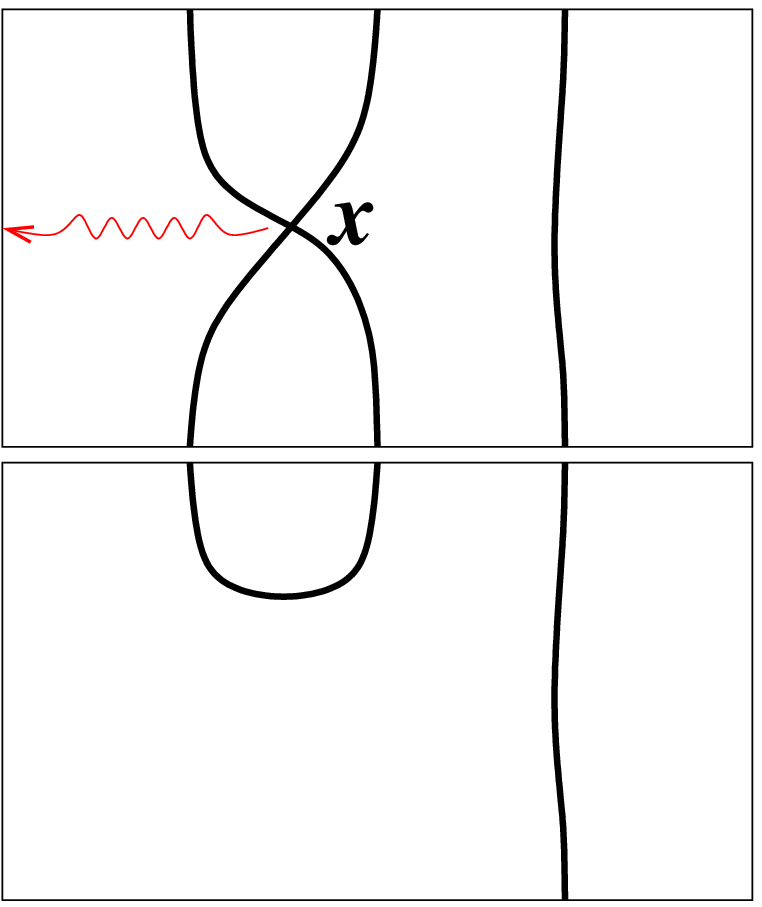} 
\label{fig:prmain}
\] 
Thus, the left-height of a point (in $d_1$ or $d_2$) cannot increase
after concatenation. 
Hence   %the left-height of the picture 
$\hht(d_1|d_2) \leq l$. (II) follows immediately.
\Qed
}

%}}}
%{{{ remarks

\begin{figure}
\[
\includegraphics[width=9cm]{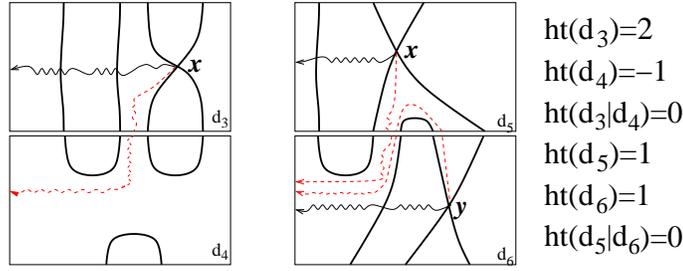} \label{fig:prmainb}
\]
\caption{Examples with heights shown. \label{fig:bellow}}
\end{figure}

{\mrem{
The left-height $ht(d|d')$ may be smaller than  that of $d$ or
$d'$, 
due to paths (in red in Fig.\ref{fig:bellow}) in $d|d'$, which are not 
paths in either $d$ or $d'$. 
}}

%}}}
%{{{ catP_l theorem and remarks

{\mth{
There is a subcategory $\catP_l = (\N_0 , \kk \Pl{l}(n,m) , *)$ of $\catP$.
}}
\proof{The proof of (\ref{pr:main}) works {\em mutatis mutandis}. \Qed}

\medskip

We will discuss connections with known constructions in
    \S\ref{ss:remarkies}.
%}} 

%}}}

%}}}
\section{Representation theory of $J_{l,n}$}
%{{{ rep
%{{{ rep

We now begin to examine the representation theory of $J_{l,n}$.
We follow a tower of recollement (ToR)
approach \cite{CoxMartinParkerXi06}. 

%}}}
%{{{ footnote: J(n,m), P(n,l,m) etc
A part in $p \in \Part(n,m)$ is {\em propagating} if it contains both
primed and unprimed elements of $\underline{n} \cup \underline{m}'$ 
(see e.g. \cite{Martin91,Martin94}). 
Write $\Part(n,l,m)$ for the subset of $\Part(n,m)$ of partitions with $l$
propagating parts; 
and similarly $J(n,l,m)$. 
Define 
$$
J_l(n,r,m)  = J_l(n,m) \cap J(n,r,m) .
$$ 
We write $\hashp (p)$ for the {\em propagating number}
--- the number of propagating parts. 
\label{def:hash}
%}

%}}}
%{{{ elements (\ee ...)
\newcommand{\ee}{{\mathsf e}}  
\newcommand{\botimes}{\otimes}  

Define $u$ as the unique element in $J(2,0)$.
Note that $\catP$ and $\catB$ are isomorphic to their respective
opposite categories via the opposite mapping $c \mapsto c^\flop$. 
Thus $u^\flop$ is the unique element in $J(0,2)$.
Define U to be the pair partition in $J(2,2)$ determined by the
following picture. 
\[
\qquad 
U := u u^\flop
= \raisebox{-.21in}{\includegraphics[width=.32in]{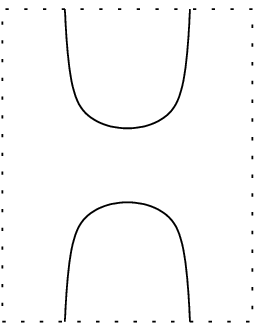}}
\]
We use $\botimes$ to denote the monoidal/tensor category composition in 
$\catP$  %the Brauer category 
(that is, the image of the 
side-by-side concatenation of pictures from (\ref{de:otimespic}),
extended $k$-linearly). 
For any given $n \geq 2$, set $\ee \; = \; e_n$ where
\[
e_n \; = \; 1_{n-2} \otimes U \; = \;
\raisebox{-.21in}{\includegraphics[width=1.1in]{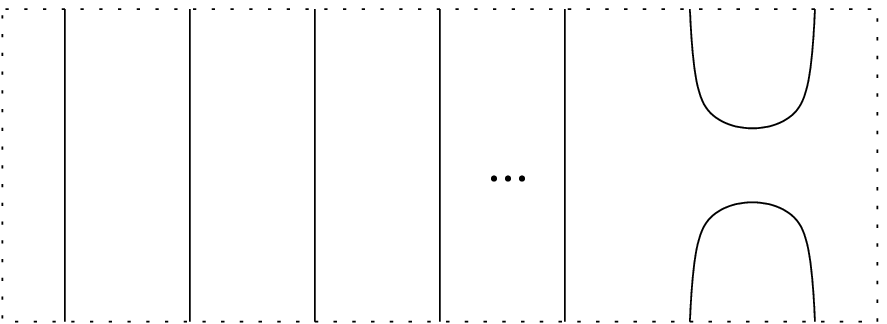}}
\; \in J(n,n) 
\]
Given a partition $p$ in $P(n,n)$ we write $p|_{n-2}$ for the natural
restriction to a partition in $P(n-2,n-2)$. (Note that this
restriction does not take $J(n,n)$ to $J(n-2,n-2)$.) 

%}}}
\subsection{Index sets for simple $J_{l,n}$-modules}
%{{{ INDEX
%{{{ index

Here we assume we have base-changed  as in (\ref{eq:base}) to a
field $k$. 
We write $\Jln$ for $\Jln(\delta_c)$ if we do not need to
emphasise $\delta_c$.
Write $k^\units$ for the group of units.
{\em For simplicity only}
we generally 
assume $\delta=\delta_c$ is invertible in $k$ (e.g. as in $\Bbbk$). 

%}}}
%{{{ eJe

\mpr{ \label{pr:eJe}
Suppose $\delta \in k^\units$ or $n>2$.    %%is invertible in $k$.
%For $n \geq 2$  % and $n > l$ 
There is a $k$-module map
\[
\Psi : \ee J_{l,n} \ee \cong J_{l,n-2}
\]
given by
$
\Psi : \ee d \ee \mapsto \ee d \ee |_{n-2}  . 
$
For $\delta \in k^\units$ (and $n \geq 2$) this $\Psi$ is an algebra isomorphism.
}

\proof{ Note $\ee d \ee = d' \botimes U$ for some $d' \in J_{l,n-2}$, so 
$\Psi(\ee d \ee) = d'$ so the map is injective (just as in the ordinary
  Brauer case). 
To show surjectivity in case $\delta \in k^\units$ 
consider $d'$ in $J_{l,n}$ by the natural
inclusion of $J_{l,n-2} \hookrightarrow J_{l,n}$
(the key point here is that the natural inclusion $J_{n-2}
\hookrightarrow J_n$ takes  $J_{l,n-2} \hookrightarrow J_{l,n}$
since the embedding does not change the height of crossings in the
$d'$ part, as it were, and does not introduce further crossings), so 
$\Psi ( \ee d' \ee) = \delta d'$ for any $d' \in J_{l,n-2}$. 
Other cases are similar.
\Qed
}

%}}}
%{{{ Corol

{\mco{
Suppose that $\Lambda(J_{l,n})$ denotes an index set for classes of
simple modules of $J_{l,n}$, for any $n$.
Then the set of classes of simple modules $S$ of $J_{l,n}$ such that 
$\ee S \neq 0$ may
be indexed by  $\Lambda(J_{l,n-2})$.
}}

\proof{Note that $\delta^{-1} \ee$ is idempotent and apply Green's Theorem in 
\cite[\S6.2]{Green80}. 
\Qed}

\medskip

{\mco{
\label{pr:green1}
The index set $\Lambda(J_{l,n})$ may be chosen so that 
\[
\Lambda(J_{l,n}) \setminus \Lambda(J_{l,n-2}) = 
  \Lambda(J_{l,n} /(J_{l,n} \ee J_{l,n} ))
\]
where  $   \Lambda(J_{l,n} /(J_{l,n} \ee J_{l,n} ))$ is an index set for
simple modules of the quotient algebra by the relation $\ee = 0$. 
In other words   
$\Lambda(J_{l,n}) \cong 
    \Lambda(J_{l,n-2}) \sqcup   \Lambda(J_{l,n} /(J_{l,n} \ee J_{l,n} ))$
and the sets  $\{ \Lambda(J_{l,n}) \}_n$ are determined iteratively
by the sets $ \{  \Lambda(J_{l,n} /(J_{l,n} \ee J_{l,n} )) \}_n$. 
\Qed
}}

%}}}
%{{{ JeJ prop

Define 
$J_l(n,<\!\!m,n) \; = \;  \cup_{r<m} J_l(n,r,n)$ 
and so on (e.g. $\Jl{l}(n,<\!n,n)$ includes every pair partition in
$\Jl{l}(n,n)$  
with submaximal number of propagating lines). 

{\mpr{ \label{pr:JeJ} %CLAIM. 
The ideal $J_{l,n} \ee J_{l,n} \; 
= \; k \Jl{l}(n,<\!\!n,n) 
$    %%%\; := \; k \cup_{r<n} J_l(n,r,n)$ 
as a $k$-space. 
}}

%}}}
%{{{ proof

\proof 
Let $p \in \Jl{l}(n,n)$ have submaximal number of propagating
lines, that is $\hashp(p) \leq n-2$ (see page \pageref{def:hash}),
so that it has at least one northern and one southern pair.
Let $d$ be a low-height \ppicture\ of $p$. 
We will use $d$ to show that $p \in J_{l,n} \ee J_{l,n}$.

For $X$ some subset of $\{N,S\}$, let $d[-X]$ denote the picture
obtained from $d$ by deleting the lines 
from north to north ($N$), south to south ($S$), or both.
Thus $d[-S] $ is a picture of some $p_t \in J(n,j,j)$
where $j = \hashp(p)$.
Similarly, abusing notation slightly by writing $u^\flop$ for some low-height
\ppicture\ of $u^\flop$, then 
$
d[-S] \botimes u^* \;
$ 
%%\danger \\
is a picture of some 
$
p' \in \; J(n,j,j+2) .
$

\begin{figure}
\eql(figeq1)
\includegraphics[width=13.6cm]{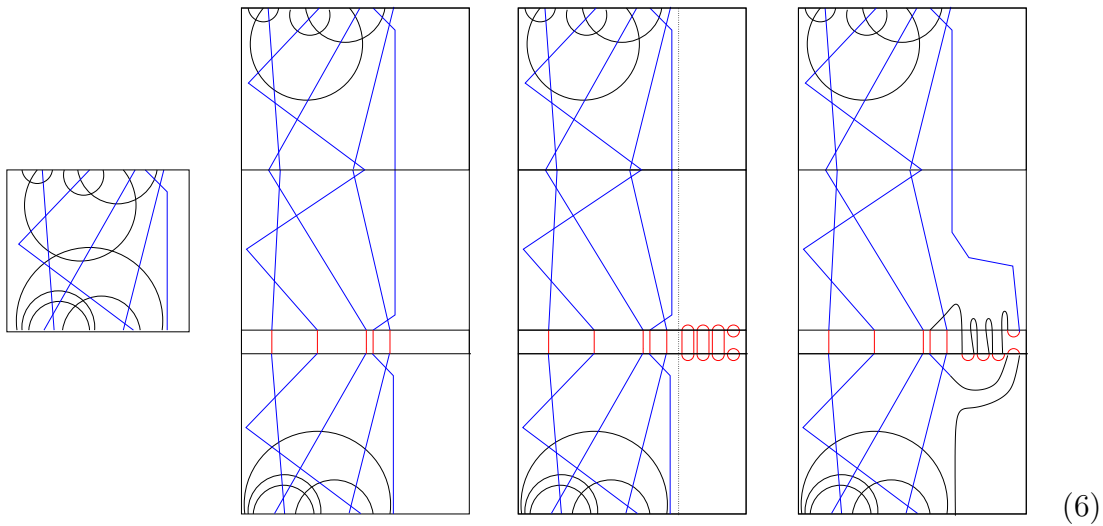}
\eq
\caption{Writing $p \propto  p' e p''$. \label{fig:figeq1}}
\end{figure}

By the Deletion Lemma~\ref{lem:deletion} the height of $d[-N]$ does
not exceed that of $d$, and similarly for $d[-NS]$ and $d[-S]$.  
Note that, since $d[-NS]$ is a picture of a permutation (of the
propagating lines); and $d[-NS]^\flip$ is a picture of the inverse, 
we have that in the picture category (Prop.\ref{pseudopicturecategory}),
\begin{equation} \label{eq:del1}
d' \;\; = \;\;\;  d[-S] \; | \;  d[-NS]^\flip  \; | \; d[-N]\label{deco}
\end{equation} 
is another picture for $p$. 
(An example is provided by 
Fig.\ref{fig:figeq1}:
the original picture of $p$ is on the left, whereas  
$d[-S] \; | \;  d[-NS]^\flip  \; | \; d[-N]$ is on the 
center-left.) 

%}}}
%{{{ proof contd

Next, observe that one can add some loops on the right of this picture, 
(red circles in the example 
Fig.\ref{fig:figeq1})  %%figure in (\ref{figeq1}))
with no change to the height of the picture, or the resulting partition $p$.
Thus,
up to an overall factor of a power of $\delta$, 
$p$ can be expressed in the form $p' \ee p''$ where 
a picture for $p'$ is 
$ ( d[-S] \; | \;  d[-NS]^\flip ) \otimes (u^*)^{\otimes (n-j)/2}$
(a picture of height $\leq l$ by construction), 
and $p''$ is   %%similar, 
$d[-N]\otimes u^{\otimes (n-j)/2}$,
and hence
$p',p'' \in J_{l,n}$.
Finally note that the red loops can be replaced by 
suitable non-crossing deformations of lines from above and below
(cf. the rightmost picture in example Fig.\ref{fig:figeq1}). 
\Qed

%}}}
%{{{ IGNORE \bittwo

\newcommand{\bittwo}{{
\subsection{Lemmas for (\ref{pr:JeJ})}
%{{{ lems

Let $J(n,1_m,m) \subset J(n,m,m)$ be the subset in which propagating
lines do not cross. Evidently (see e.g. \cite{}) we can resolve
$p \in J(n,m,n)$ by 
\[
J(n,m,n) \stackrel{\sim}{\rightarrow} 
    J(n,1_m,m) \times J(m,m,m) \times J(m,1_m,n)
\]
which we will write as $p \mapsto (|p\rangle , c(p) , \langle p |)$.

\mdef The above decomposition of partition $p$ is called a {\em polar
  decomposition} of the partition. 

\mdef
A picture is called {\em north polar} if no crossing involves two
propagating lines; and every line that touches the southern edge is
propagating. 
Note that every $| p \rangle$ has a north polar picture.

A polar decomposition of a {\em picture} is a picture which can be separated
into three parts by two horizontal lines, such that the top part is
north polar; 
the middle part is a permutation of propagating lines; 
and the bottom part is south polar. 

Note that every  $p$ has a picture with a 
polar decomposition obtained by stacking
pictures from the polar decomposition of $p$. 

\medskip

First we aim to show the following.

{\mlem{ 
The resolution of an element of 
$J_{\leq l}(n,n)$ obeys: $|p\rangle , c(p) , \langle p | $ have height
$\leq l$. 
}}

\proof 
The idea is to suppose that we start with a low-height picture of $p$
and show that we can move crossings of propagating parts out of the
`polar regions' without making the {\em picture} height higher. 

In order to do {\em this} we will need some more Lemmas.

...

%}}}
\subsection{More Lemmas: Pictures again}
%{{{ pics

\mdef Now consider a line $l$ that is a cup line (north-to-north) without
self-crossings in a picture. 
Suppose there is a crossing of propagating lines in the alcove of $l$
(in the obvious sense). 

...then we claim this crossing can be `moved out of the alcove' without 
raising the picture height.

An example of this is:
\[
\includegraphics[width=2.6in]{reidem01.eps}
\]
This case agrees with the claim. However this case is not generic. 
There may be other lines present, such as here:
\[
\includegraphics[width=2.96in]{reidem02.eps}
\]
%\end{figure}

%}}}
%{{{ braid move

\mdef In order to address this we consider the effect on height of an
arbitrary braid move. 

The consequences for height of the `braid move' depends on
whether the lines involved are propating or not. 
Here are typical cases: in Fig.\ref{fig:YB1x};
Fig.\ref{fig:YB2x}.

\begin{figure}
\caption{One checks that the `local' height is $h+1$ in both of these
  renderings (assuming the leftmost alcove of height $h$, as marked).
The high crossings are marked; and the directions of increasing height
across a line are marked.
\label{fig:YB1x}}
\[
\includegraphics[width=3.96in]{reidem044.eps}
\]
\end{figure}

\begin{figure}
\caption{One checks that the `local' height is $h$ in both of these
  renderings (assuming the leftmost alcove of height $h$, as marked)
--- the high crossings are marked.
\label{fig:YB2x}}
\[
\includegraphics[width=3.96in]{reidem045.eps}
\]
\end{figure}

\mdef Note (e.g. from \cite{Martin12w})
\[
\includegraphics[width=2.6in]{reidem03.eps}
\]
that this move does not increase height locally or globally.

%}}}
}}

%}}}
%\subsection{Back to the representation theory}

\newcommand{\pl}{p_l}  %% the min(p,l+2)
\newcommand{\pll}[1]{\min( #1 , l+2 )}  %% ...explicitly

%{{{ back - simple index theorem

%{\em Proof of (\ref{pr:JeJ}) continued.}

\mdef \label{pa:minnow} %\muco %In this case
By (\ref{pr:JeJ}) the quotient algebra 
$J_{l,n} /(J_{l,n} \ee J_{l,n} )$ has a basis 
which is the image of  %%corresponding to 
$\Jl{l}(n,n,n)$.
Note that these elements of $J_{l,n}$ form a subgroup as well
as their image spanning a quotient algebra. 

For $n<l+2$ this group %quotient 
is isomorphic to $ S_n$; 
otherwise it is isomorphic to $ S_{l+2}$
(since there can be no crossings after the first $l+2$ lines). 
The quotient itself is then the corresponding group algebra.
\ignore{{
For given $l$ define 
$
\pl \; = \; \min(p,l+2)   .
$
Then we have 
}}
That is, 
$\Lambda(J_{l,n} /(J_{l,n} \ee J_{l,n} )) \cong \Lambda (k S_{\pll{n}} )$.
Combining with Prop.\ref{pr:green1}, we thus have the following. 

%}}}
%{{{ the theorem

\newcommand{\LambdaJln}{\Lambda(J_{l,n})}  %% index set for simples

{\mth{ \label{th:indexsimples}
Let  $\Lambda_n $ denote the set of integer partitions of $n$
(so that $\Lambda(\C S_n) = \Lambda_n$ for any $n$). 
Let $\Lambda_n^{(p)} \; := \; \{ p \} \times \Lambda_n$ 
(so that $\Lambda_n^{(1)}$ and  $\Lambda_n^{(2)}$ are disjoint copies
of $\Lambda_n$). 
Define $\Gamma_n = \{ n, n-2, ... , 1/0 \}$.
\ignore{{
The index sets $\Lambda(J_{l,n}/(J_{l,n} \ee J_{l,n} ))$
 for simple modules %of $J_{l,n}$ 
over $\C$, say, are: % then known:
\[
\Lambda(J_{l,n}/(J_{l,n} \ee J_{l,n} )) 
 = \left\{ \begin{array}{ll} 
     \Lambda(\C S_n) & n < l+2
\\   \Lambda(\C S_{l+2}) & n \geq l+2   \end{array} \right.
\] 
%%%%  FIX ME!
and hence, over $\C$,
}}
Then over $k=\C$,
\[
\Lambda(J_{l,n}) \;
 =  \; \bigcup_{p\in \Gamma_n}  \Lambda^{(p)}_{\pll{p}} \; 
 \; \cong \; \left( {\bigcup_{p=0}^{l+1}}{}^{'} \Lambda_p \right) 
      \cup \left( {\bigsqcup_{p=l+2}^{n}}\!\!{}^{'} \Lambda_{l+2} \right)
\]
where $\bigcup'_p$ denotes a range %in steps of 2, 
including only $p$ congruent to $n$ mod.2.
\Qed
}}

%}}}
%{{{ polar decomp

\newcommand{\Jid}{J^{||}}  %% hook for no crossing propagating lines set
\newcommand{\Jidl}[1]{\Jl{#1}^{||}}  %% hook for no crossing propagating lines set
\newcommand{\nux}{\nu}   %%  polar mult map
\newcommand{\nui}{\nux^{-1}} %% inverse polar mult map

\mdef \label{de:polard}
Here $\Jid(n,l,m)$ denotes the subset of $J(n,l,m)$ of elements $p$
having a picture $d$ for which $d[-NS]$ 
has no crossings.  
Recall (e.g. from \cite{Martin09a})
the {\em polar decomposition} of an element of $J(n,m,n)$:
the inverse of the map
$
\nux : \Jid(n,m,m) \times J(m,m,m) \times \Jid(m,m,n) 
          \stackrel{\sim}{\rightarrow} J(n,m,n)
$, 
given by the category composition.
Note that if $p \in \Jl{l}(n,m,n)$ then $l$ bounds the height of all
three factors in the polar decomposition. 
(The argument is analogous to the argument at (\ref{eq:del1}).
Firstly note that if $d$ is a low-height picture of $p$ then 
$d[-S] \; | \; d[-NS]^\flip$ is a picture of the northern polar factor of no
higher height. 
The other factors are similar.) 
That is, the restriction
$\nui : \Jl{l}(n,m,n) \hookrightarrow  
   \Jidl{l}(n,m,m) \times \Jl{l}(m,m,m) \times \Jidl{l}(m,m,n)
$. 
On the other hand the image of $\nux$ on this codomain lies in 
$ \Jl{l}(n,m,n)$ by Th.\ref{th:main}, so the restriction as given is a
bijection. 

%}}}
%}}}
\subsection{Quasiheredity of $J_{l,n}$ }
%{{{ JeeeeJ

The proof of the main result of this section (Thm.\ref{pr:qher1})
follows closely the Brauer algebra case, as for example in 
\cite{CoxDevisscherMartin0609}.
We focus mainly on the new features required for the present case.

\mdef For $n \geq 2t$ define 
$ \; 
e_{n,t} = 1_{n-2t} \otimes U^{\otimes t}
$
and (when $\delta\in k^\units$) $e'_{n,t} =\delta^{-t} e_{n,t}$.
Note that $e_{n,t} \in J_{-1}(n,n-2t, n)$. 
For example $\ee =e_{n,1}$ and 
\[
e_{n,2} = 
\raisebox{-.21in}{\includegraphics[width=1.1in]{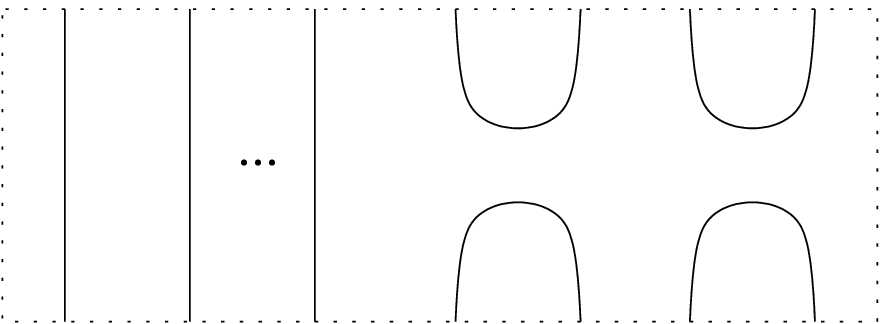}}
\]
We have the following useful corollary to the proof of Prop.\ref{pr:JeJ}. 

{\mcor{  \label{pr:JeeeJ} %CLAIM. 
Provided that $\delta \in k^\units$,   %%is a unit in $k$, 
the ideal $J_{l,n} e_{n,t} J_{l,n} = k \Jl{l}(n,\leq n-2t,n)$. 
}}

\proof{
The proof is the same as for Prop.\ref{pr:JeJ} except that we use more
pairs of 
loops, instead of single loops, 
in the final stage of the construction of $d \in J_l(n, \leq n-2t , n)$
(the red loops in (\ref{figeq1})).
Note that there is room for enough loops because of the bound on the
number of propagating lines in $d$.
\Qed
}

%}}}
%{{{ quotient algebra

\mdef Define the quotient algebra
\[
J_{l,n,t} \; = \;  J_{l,n} / J_{l,n} e_{n,t+1} J_{l,n}
\]
By Prop.\ref{pr:JeeeJ} this algebra has basis 
$\Jl{l}(n,\geq n-2t,n)$. 

{\mpr{ \label{pr:heredid}
For each triple $n,l,t$ the following hold when $\delta \in k^\units$. 
\\
(i) The algebra $A = e_{n,t} J_{l,n,t} e_{n,t}$ is semisimple over $\C$. 
\\
(ii) The multiplication map 
$
J_{l,n,t} e_{n,t} \otimes_{A} e_{n,t} J_{l,n,t} 
 \stackrel{\mu}{\rightarrow} J_{l,n,t} e_{n,t} J_{l,n,t}
$
is a bijection of $J_{l,n,t} , J_{l,n,t}$-bimodules.
}}

\proof{
(i) The number of propagating lines must be at least $n-2t$, but with
  the $e_{n,t}$'s   %%idempotents 
present this is also the most it can be, so every
  propagating line in $e_{n,t}$
is propagating in $A$, and indeed $ A$ is isomorphic to the group
  algebra of a symmetric group.
\\
(ii) The map is clearly surjective. We construct an inverse
using the polar decomposition  %%(\ref{MISSING:polard}) 
(\ref{de:polard}).
Note that $J_{l,n,t} e_{n,t}$ has a basis in bijection with 
$\Jl{l}(n,n-2t,n-2t) \cong 
    \Jl{l}^{||}(n,n-2t,n-2t) \times \Jl{l}(n-2t,n-2t,n-2t)$.
Subset $ e_{n,t} J_{l,n,t}$ can be treated similarly.
It then follows from the definition of the tensor product that the
left-hand side has a spanning set whose image is independent on the right.
\Qed
}

{\mth{ \label{pr:qher1} \label{th:qher1}
If $\delta$ invertible in $k=\C$ then $J_{l,n}$ is quasihereditary,
with heredity chain $(1,e'_{n,1}, e'_{n,2},...)$.
}}

\proof This follows from (\ref{pr:heredid}). Cf. e.g. 
\cite{DlabRingel,ClineParshallScott88,Donkin98,CoxDevisscherMartin0609}.
\Qed

%}}}
\subsection{Aside: Slick proof of quasiheredity in the monoid case}
%{{{ qh

Note that for $a \in J(n,n)$ we have $aa^* a = a$ in $J_{l,n}$.
Since the height of $a^*$ is the same as $a$ we have the following.

\mpr{
The algebra $J_{l,n}(1)$ is a regular-monoid $k$-algebra
(i.e. $a \in aJ_{l,n}a$ for all $a \in J_{l}(n,n)$).
\Qed
}

{\mcor{
The algebra
$J_{l,n}(1)$ is quasihereditary when $k = \C$.
}}

\proof{Use Putcha's Theorem \cite{Putcha98}. \Qed}

%}}}
\subsection{%Quasiheredity and s
Standard modules of $J_{l,n}$}
%{{{ STANDARDS
%{{{ standard modules

\newcommand{\Jak}{{\mathfrak J}} %% hook for section modules

We may construct a complete set of `standard modules' for each $J_{l,n}$
as follows. 
The modules we construct are `standard' with respect to 
a number of different compatible axiomatisations
(the general idea of standard modules, when such exist, 
is that they interpolate between
simple and indecomposable projective modules). 
For example (I) we can construct quasihereditary standard modules by
enhancing the heredity chain in Th.\ref{th:qher1} to a maximal chain
cf. \cite{Donkin98}; 
(II) we can construct the modular reductions of lifts of
`generic' irreducible modules in a modular system cf. \cite{Benson95};
(III) we can construct %inflations and 
globalisations of suitable
modules from lower ranks cf. \cite{Martin96,KoenigXi99,Martin09a}.
We are mainly
interested in a useful upper-unitriangular property of decomposition
matrices that we establish in (\ref{pr:uutri}).

%}}}
%{{{ bimodules

\mdef
Note that $k\Jl{l}(n,m)$, with $m\leq n$ say,  is a left $J_{l,n}$ right
$J_{l,m}$-bimodule, by the category composition.
By the bottleneck principle
it has a sequence of submodules:
\[
k\Jl{l}(n,m) = k\Jl{l}(n,\leq m,m) \supset k\Jl{l}(n,\leq m-2,m) \supset ...
\]
For given $l$, each section
$$
\Jak_{n,m}^p \; : = \; k\Jl{l}(n,\leq p,m) / k\Jl{l}(n,\leq p-2,m)
$$
thus has basis $\Jl{l}(n,p,m)$. 
In particular the top section has basis $\Jl{l}(n,m,m)$. 

%}}}
%{{{ in partic (quotient II)

\mdef \label{de:cX} \label{de:propideals}
The above holds in particular for the case $m=n$, where our sequence is an
ideal filtration of the algebra, cf. Corol.(\ref{pr:JeeeJ}). 
Define quotient algebra
\[
J_{l,n}^{/p} = J_{l,n} / k J_l (n , \leq p , n)
\]
Note that this is the same as $J_{l,n,t}$ with $p+2 =n-2t$
(but now without restriction on $\delta$).
The index $p$ tells us that partitions with $p$ or fewer propagating
lines are congruent to zero in the quotient. 

\mdef \label{de:cY}
In particular, as noted in (\ref{pa:minnow}), 
$
J_{l,n}^{/n-2} \cong k S_{\pll{n}}  .
$
\ignore{{
\[
J_{l,n}^{/n-2} \cong 
\left\{ \begin{array}{ll}
 k S_{l+2}  \qquad  & n \geq l+2
\\ k S_n   \qquad  & n < l+2   \end{array} \right.
\]
}}
Specifically
$\Jln^{/n-2}$ has basis $\Jl{l}(n,n,n)$, which is of form 
$J(l+2,l+2,l+2) \otimes 1_{n-(l+2)}$, 
i.e. $S_{l+2} \otimes 1_{n-(l+2)}$,
when $n>l+2$
(since there can be no crossing lines after the first $l+2$ in this case). 

%}}}
%{{{ functor

\mdef \label{de:cZ}
Note that $\Jak_{n,m}^p$ is also a left $J_{l,n}^{/p-2}$
right $J_{l,m}^{/p-2}$ bimodule. Thus we have a functor 
\[
\Jak_{n,m}^{p} \otimes_{J_{l,m}^{/p-2}} - : J_{l,m}^{/p-2} - mod
\; \rightarrow J_{l,n}^{/p-2} - mod
\]
and in particular a functor
\[
\Jak_{n,p}^{p} \otimes_{k S_{p}} - : kS_{p} - mod
\; \rightarrow J_{l,n}^{/p-2} - mod
\qquad (p \leq l+2)%\label{lr}
\]
\[
\Jak_{n,p}^{p} \otimes_{k S_{l+2}} - : kS_{l+2} - mod
\; \rightarrow J_{l,n}^{/p-2} - mod
\qquad (p > l+2)%\label{hr}
\]
(in case $p>l+2$, the right action of 
$kS_{\pl} =  kS_{l+2}$ in the form $ kS_{l+2}\otimes 1_{p-l-2}$
is understood).

%}}}

%{{{ specht

\newcommand{\Spec}{{\mathcal S}}  %% Specht hook

\mdef \label{de:cT}
With this functor in mind we recall some facts about the symmetric
groups. 
For any symmetric group $S_m$ 
and a partition $\lambda \vdash m$,
let $\Spec_\lambda$ denote the corresponding
Specht module of  %the symmetric group 
$S_m$ 
--- see e.g. \cite{James,JamesKerber81}.  % (any $m$). 
Recall %%for later use 
that there is an element
$\epsilon_\lambda$ 
in $k S_m$ such that $\Spec_\lambda = k S_m \epsilon_\lambda$.
If $k \supset \Q$ then $\epsilon_\lambda$ may be chosen  idempotent.
Example: 
The element $\epsilon_{(2)} \in kS_2$ is unique up to scalars:
$\epsilon_{(2)} = 1_2 + \sigma_1$ (in the obvious notation).
Thus a basis for %the one-dimensional Specht module
$\Spec_{(2)}$ is $b_{(2)} = \{ \epsilon_{(2)} \}$. 

%}}}
%{{{ global-standard

\mdef \label{de:stan1}
For given $l$ define $ \pl \; := \; min(p,l+2)$.
For any $p$ and $\lambda \vdash \pl $,   %%\; := \; min(p,l+2)$ 
that is for $ (p,\lambda) \in \LambdaJln$,
we define a global-standard $J_{l,n}$-module
\[
\Delta^n_{p,\lambda} \;  =  \;  
\Delta_{p,\lambda} \;  =  \;  % kJ_l(n,p,p) 
      \Jak_{n,p}^p \otimes_{k S_{\pl}} \Spec_\lambda
\]

%}}}
%{{{ d-module

\mdef
For given $l,p$ let $\lambda \vdash \pl$. 
We can consider 
$kS_{\pl} \hookrightarrow kJ(\pl,\pl)$, indeed
$kS_{\pl} \stackrel{\sim}{\rightarrow} k\Jl{l}(\pl,\pl,\pl)$
by the inclusion in (\ref{pa:minnow}), 
and hence 
$\epsilon_\lambda \in  k\Jl{l}(\pl,\pl,\pl)$. Thus
$$ 
E_n(p,\lambda) 
  \; := \epsilon_\lambda \otimes 1_{p-\pl} \otimes U^{\otimes (n-p)/2}
$$
lies in $J_{l,n}$. 
%}}}
%{{{ pic
For any $\epsilon_\lambda$ we may write the image in $\catP$ 
{\em schematically} as 
$\framebox{$\; \lambda \;$}$ with $2\times |\lambda |$ legs; 
e.g. $\framebox{$\; (2) \;$}$ with two pairs of legs. 
Thus for example with $l=1$ we have:
\[
E_{10}(6,\lambda\vdash 3) \;\; = \; 
\raisebox{-.21in}{
\includegraphics[width=2.3in]{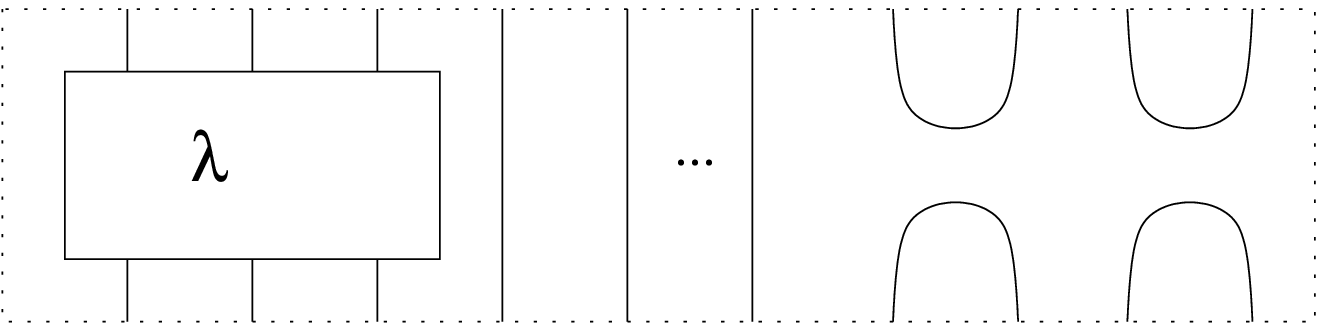}}
\]
Define  another $\Jln$-module by 
\[
D_{p,\lambda} \; 
 = J_{l,n}^{/p-2} E_n(p,\lambda)   
          %( \epsilon_\lambda \otimes 1_{x} \otimes UUU)
\]
\mupr
If $k \supset \Q$ and $\delta \in k^\times$ then 
$E_n(p,\lambda)$ is an (unnormalised) idempotent and 
$\{ D_{p,\lambda} \; | \; (p,\lambda) \in \LambdaJln \}$
are the
quasihereditary-standard modules \cite[\S A1]{Donkin98}
of the quasihereditary structure in
Theor.\ref{th:qher1}. \Qed $\;$
(We will shortly explain this explicitly.)

\medskip

%}}}
%{{{ Example

Example: 
\ignore{{
The element $\epsilon_{(2)} \in kS_2$ is unique up to scalars:
$\epsilon_{(2)} = 1_2 + \sigma_1$ (in the obvious notation).
Thus a basis for the one-dimensional Specht module
$\Spec_{(2)}$ is $b_{(2)} = \{ \epsilon_{(2)} \}$. 
}}
For any $l \geq 1$ a basis of $D_{2,(2)}$  %% $\Delta_{2,(2)}$ 
for $n=4$ is
in bijection with $J_l(4,1_2,2)$, and 
(omitting an irrelevant arc in the bottom-right) may be depicted: 
\eql(eg:42)
\includegraphics[width=5in]{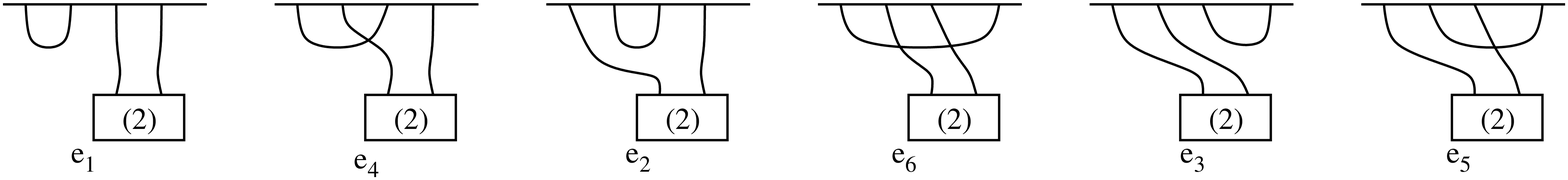}
\eq
Note that the action of the algebra in this depiction is from above.
Note that, if $l=0$, then  $J_0(4,1_2,2)$ is smaller. 

%}}}
%{{{ basis

{\mpr{\label{ebas}
Let $\lambda \vdash \pl$ and let 
$b_\lambda$ be a basis for $\Spec_\lambda = k S_{\pl} \epsilon_\lambda$. 
Then a basis for the $J_{l,n}$-module $D_{p,\lambda}$ 
is, up to (irrelevant arc ignoring) isomorphism, 
\[
B_{p,\lambda} \; 
  = \; \{  x(y \otimes 1_{p-\pl}) 
            \; | \; x \in \Jl{l}(n,1_p,p); \; y \in b_\lambda  \}
\]
regarded as a submodule of $\Jak_{n,p}^p$,
where $\Jl{l}(n,1_p,p) = \Jl{l}^{||}(n,p,p)$ as in (\ref{de:polard}). 
\ignore{{
\hspace{-0.6cm}In case $p < l+2$ the quotient algebra is not restricted by the left-height: 
$L_{l,p}^{/p-2}\cong kS_p$, so 
we have $\lambda\vdash p$ and we have no factor of $1_{p-l-2}$ in the tensor product expression
of neither $E_n(p,\lambda)$ or that for its basis.} 
}}
}

%}}}
%{{{ proof

{\proof 
By definition $D_{p,\lambda} = J_{l,n}^{/p-2} E_n(p,\lambda)   $ 
is spanned by elements 
$x(y \otimes 1_{p-\pl})$ as in
$B_{p,\lambda}$ except with $x \in \Jl{l}(n,p,p)$
(again ignoring irrelevant arcs in the bottom-right).
Thus we  need to show that we can omit $x$'s with crossing propagating
lines without breaking the spanning property. 
Some
example elements of $D_{5,\lambda}$ with $p>l+2$ (in case $p=5$,
$l=1$)
provide a useful visualisation here:
\begin{center}\includegraphics[width=12cm]{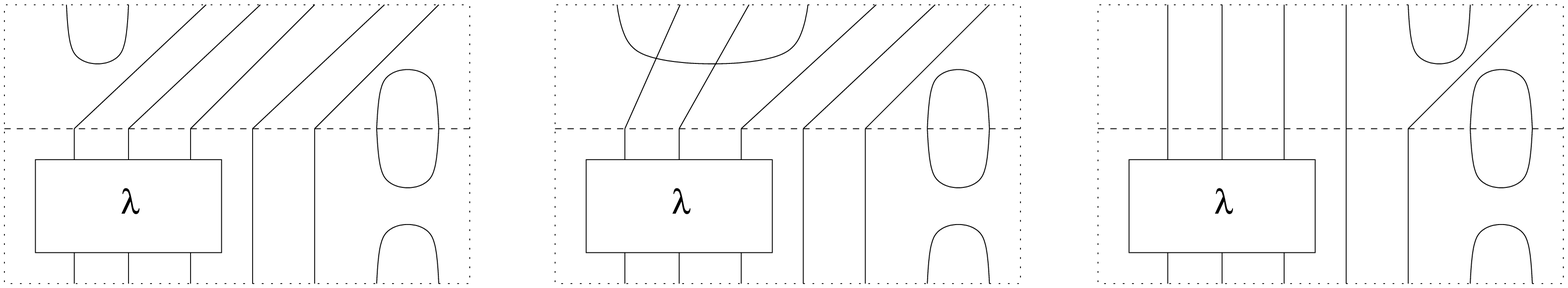}
\end{center}
First note that 
a basis element must have $p$ propagating lines by the quotient.
Any crossing in the first $\min(p,l+2)$ of these can be `absorbed' by the 
$b_\lambda$ part of the basis. There cannot be a crossing in any
remaining propagating lines by the height restriction, cf. (\ref{de:polard}).
\Qed

}

{\mcor{ \label{co:DDel}
$D_{p,\lambda} \cong \Delta_{p,\lambda}$
}}

\proof Compare our basis above with the construction for 
$\Delta_{p,\lambda}$. The main difference is  combination via
$\otimes$ rather than multiplication (up to some subtleties when $\delta=0$). 
This gives us a surjective map right to left. % to right. 
One then compares dimensions. 
\Qed

%}}}
%{{{ id
{\mpr{ \label{pr:EJEeE}
For given $l$,  $\;$
$E_n(p,\lambda) \; J_{l,n}^{/p-2} \; E_n(p,\lambda) \; = \; k E_n(p,\lambda)$.
}}

\proof{ %We claim that 
Pictorially/schematically we have:
\[
\includegraphics[width=1.5in]{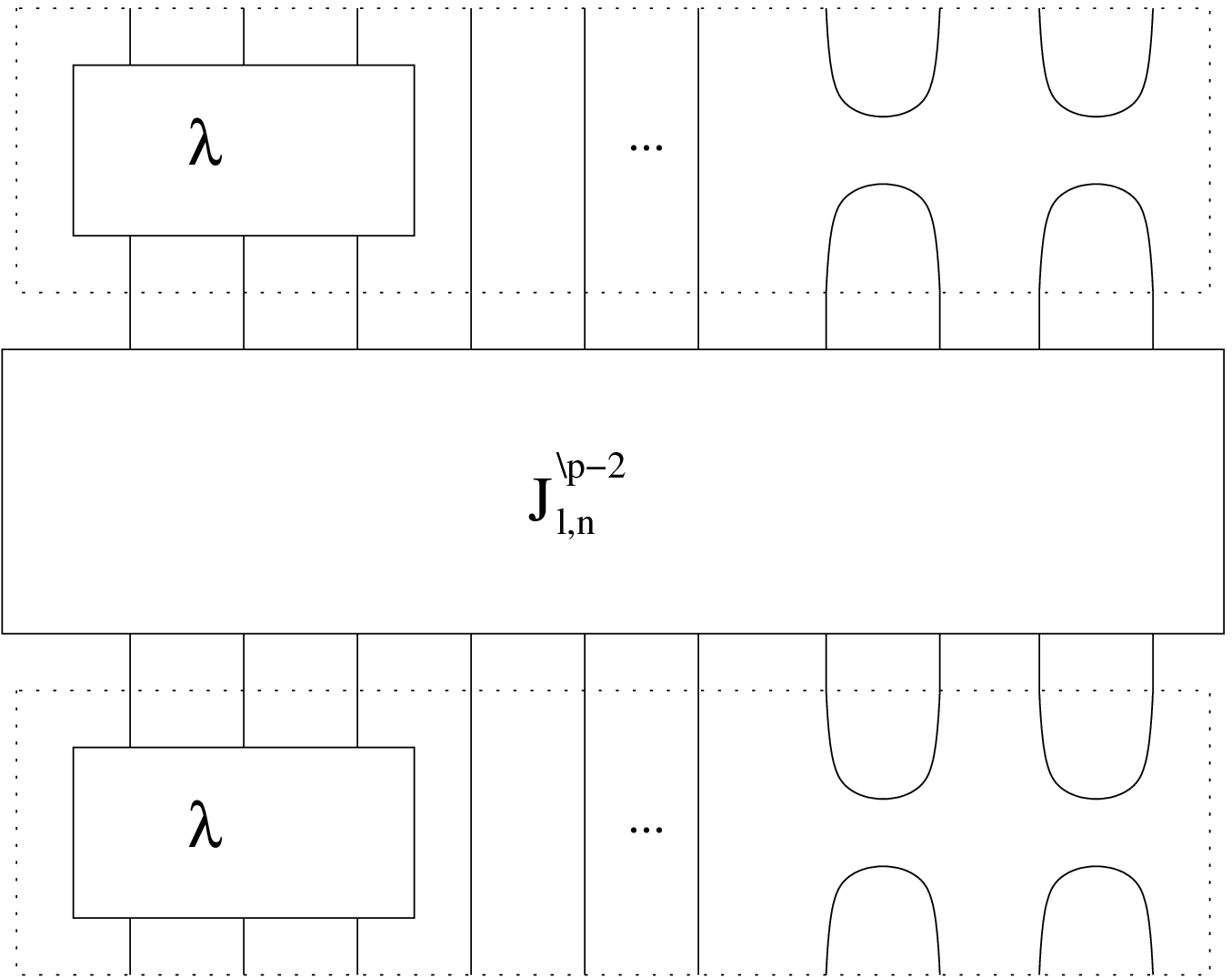}
\; %\qquad 
\stackrel{(\ref{pr:eJe}),(\ref{de:cX})}{=}
\; %\qquad
\includegraphics[width=1.5in]{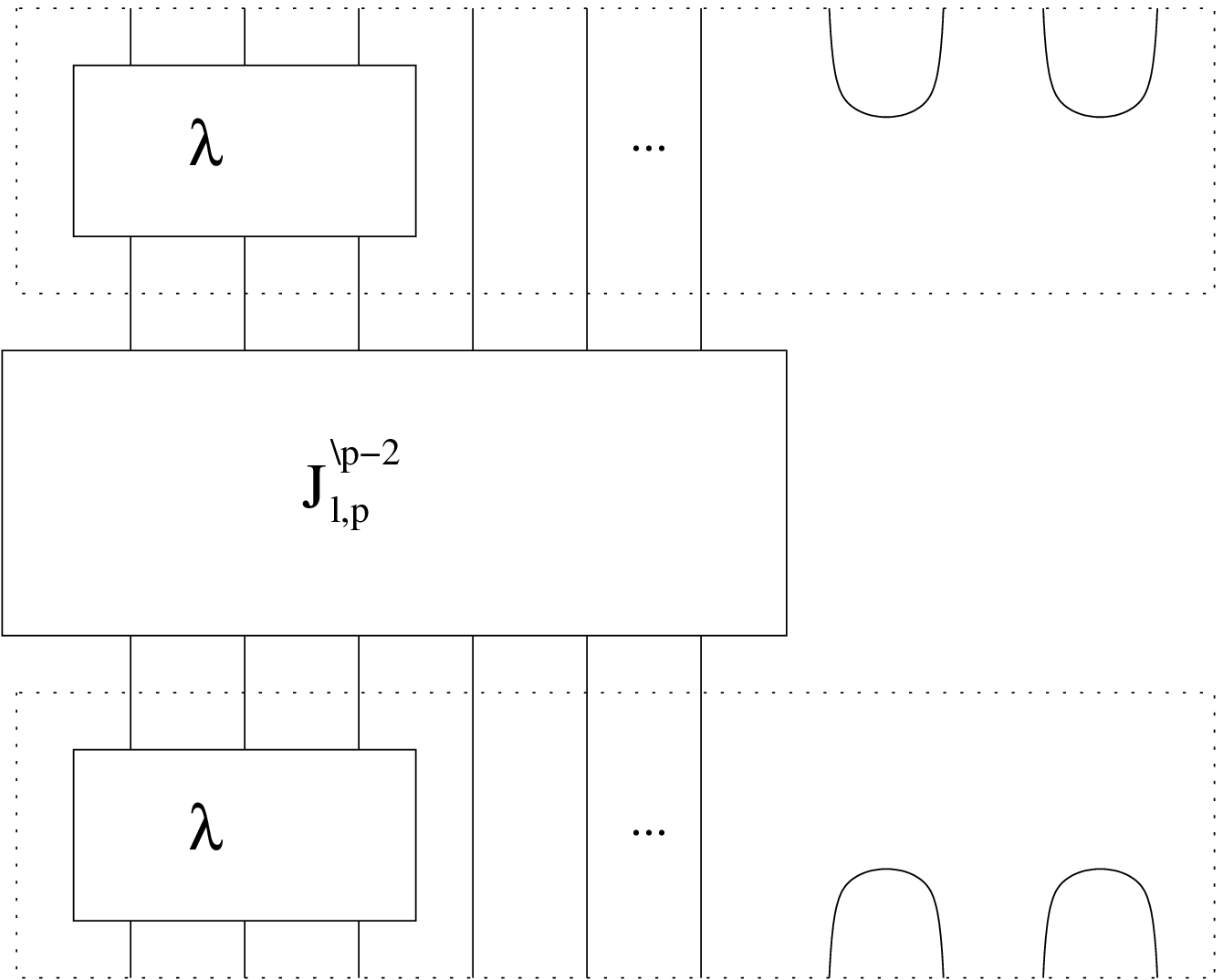}
\; %\qquad 
\stackrel{(\ref{de:cY})}{\cong}
\; %\qquad
\includegraphics[width=.95in]{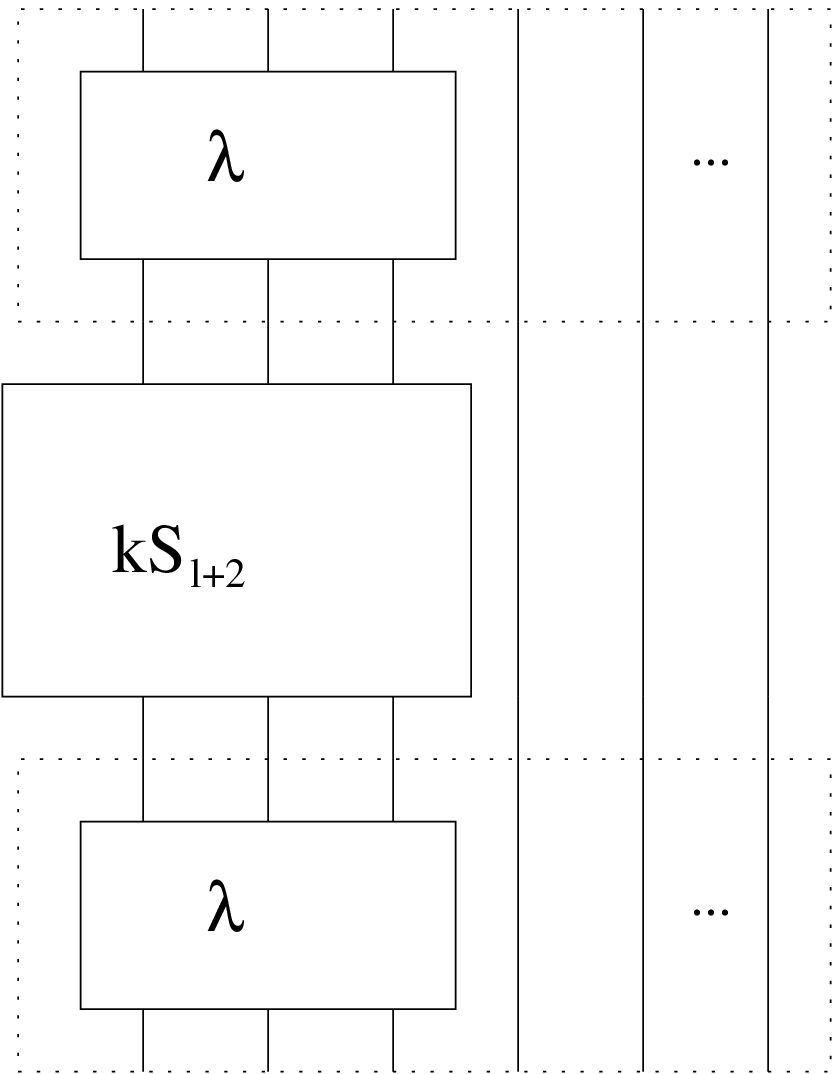}
\]
Finally 
$\epsilon_\lambda kS_{l+2} \epsilon_\lambda = k \epsilon_\lambda$
by the Specht property (\ref{de:cT}) \cite{JamesKerber81}.
\Qed
}

{\mcor{
For $(p,\lambda) \in \LambdaJln$,
$D_{p,\lambda}$ is indecomposable projective as a
    $J_{l,n}^{/p-2}$-module; and hence indecomposable with simple head
    as a  $J_{l,n}^{}$-module.
\Qed
}}

%}}}
%{{{ G_e

In light of (\ref{pr:eJe}) $ J_{l,n+2} e_{n+2}$ is a right
$\Jln$-module. Thus
we have also the adjoint pair of functors 

$$
\xymatrix{  J_{l,n}-mod \;\;\; \ar@/^/[r]^{G_e}  
& \;\;\; J_{l,n+2}-mod \ar[l]^{F_e}
}
$$
given by $G_e M \; = \; J_{l,n+2} e_{n+2} \otimes M $ 
$\;$
and $\;$
$F_e N = e_{n+2}N$. 

%}}}
%{{{ G \Delta

{\mpr{ \label{pr:GeD}
For $\delta \in k^\times$, 
$$
G_e \Delta^{n}_{p,\lambda} = \Delta^{n+2}_{p,\lambda}\ ,
$$
and
$$
F_e \Delta^{n+2}_{p,\lambda} = \left\{ \begin{array}{ll}
    \Delta^n_{p,\lambda} &  p\leq n \\
    0 & p=n+2 \end{array} \right.
$$
}}

\proof See e.g. \cite{Martin09b}. In this proposition our new geometrical
constraints do not affect the argument. 
\Qed

%}}}
%{{{ unitri

Write $L_{p,\lambda} = \head( \Delta_{p,\lambda})$ for the simple head
(see e.g. \cite[\S1.2]{Benson95}).
Write $[\Delta_{p,\lambda} : L_{p',\lambda'}]$ for the multiplicity of
$L_{p',\lambda'}$ as a composition factor in $\Delta_{p,\lambda}$;
and
 $C^\Delta = ([\Delta_{p,\lambda} : L_{p',\lambda'}])_{(p,\lambda),(p',\lambda')}$ 
for the corresponding decomposition matrix.
It is routine %straightforward 
to check the following. 

{\mpr{ \label{pr:uutri}
(I) The modules $\{ \Delta_{p,\lambda} \}$ 
are a complete set of standard modules in the
quasihereditary algebra cases ($\delta$ invertible in $k=\C$). 
(II) The simple decomposition matrix 
$C^\Delta $  
for this set of modules is upper
unitriangular (when written out in any order so that 
$(p,\lambda) > (p',\lambda')$ when $p> p'$). 
}}

\proof (I) Follows from Corol.\ref{co:DDel}. 
(II) There are several ways to prove this. Prop.\ref{pr:GeD}
implies that the composition factors with the same label in 
$ \Delta^{n}_{p,\lambda} $ and $ \Delta^{n+2}_{p,\lambda} $ have the
same multiplicity. The only possible new factors in  
$ \Delta^{n+2}_{p,\lambda} $ have the property $e_{n+2} S=0$.
Applying this iteratively on $n$ gives the claimed result.
\Qed

%}}}
\subsection{Towards the Cartan decomposition matrices}
\newcommand{\gram}{Gram}  %% gram matrix capital letter
%{{{ preamble

The Cartan decomposition matrix encodes the fundamental invariants of
an algebra \cite{Benson95}. 
Given the difficulty experienced in computing them in case $l=-1$ and
particularly $l=\infty$ \cite{Martin09a}
we can anticipate that they will not be easy to determine
in general. 
However, the main tools used in cases  $l=-1$ and
$l=\infty$ can be developed in general, as we show next. We start with
a corollary to Prop.\ref{pr:uutri}.

%}}}
%{{{ contravar form

{\mcor{
Module $\Delta_{p,\lambda}$ has a contravariant form  %n inner product 
(with respect to $\flip$)
that is unique up to
    scalars. The rank of this form determines the dimension of the
    simple module $L_{p,\lambda}$. 
}}

\proof The space of contravariant forms is in bijection with 
the space of module maps from $\Delta_{p,\lambda}$ to its contravariant dual
(the analogous right-module $ E_n(p,\lambda) J_{l,n}^{/p-2}$
treated as a left-module via 
ordinary duality).
But by the upper-unitriangular property (\ref{pr:uutri})
this space is spanned by any
single map from the head to the socle.
\Qed

\medskip

The contravariant form here is the analogue of the usual form for the
Brauer algebra \cite{HanlonWales94,Martin09a}. That is, 
a suitable inner product $\langle e_i, e_j\rangle$   is given by 
$e_i^\flip e_j \; = \;  \langle e_i, e_j\rangle E_n(p,\lambda)$
(this is well-defined by (\ref{pr:EJEeE})).
The form rank is given by the matrix rank of the gram matrix over a basis.

Example: A basis for $\Delta^4_{2,(2)}$ is given 
in (\ref{eg:42}).
In particular then (with $\lambda=(2) $):

\begin{center}
\includegraphics[width=9cm]{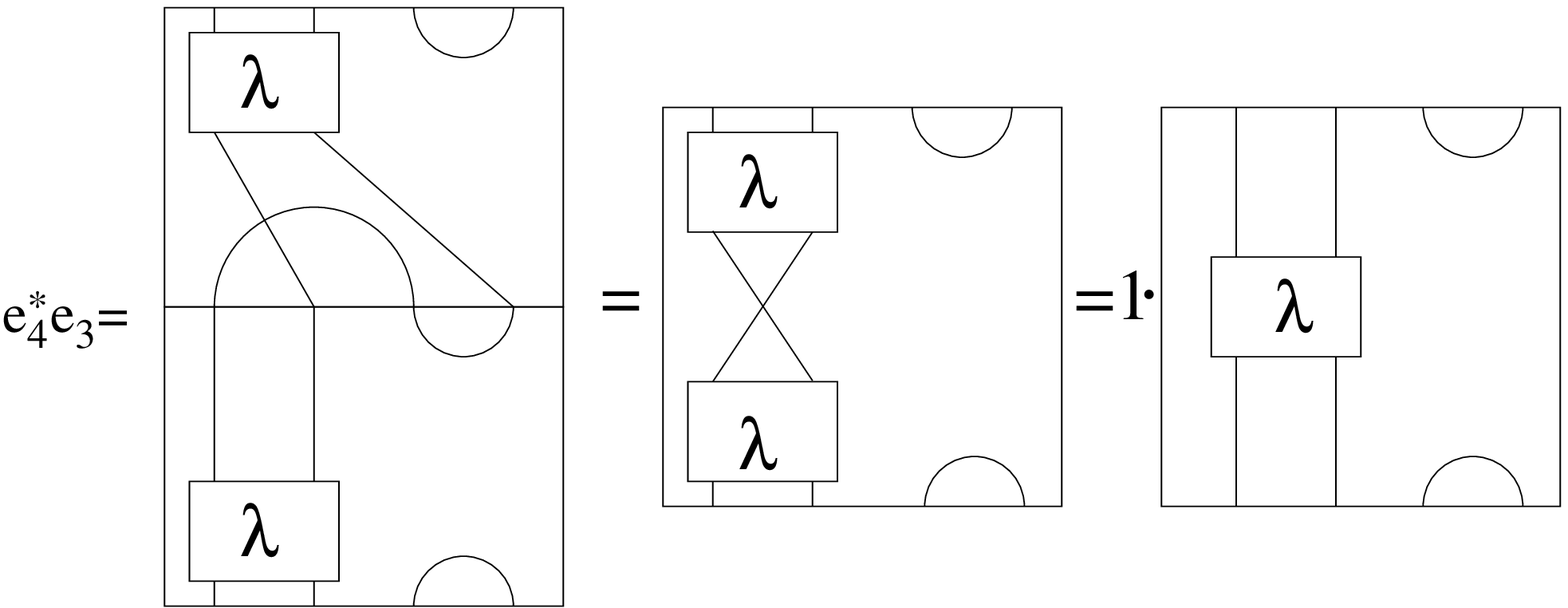}
\end{center}
%}}}
%{{{ Analysis
Here $f_{3,4}\; := \; \langle e_3,e_4 \rangle =1$ since 
$\sigma_1 \epsilon_{(2)} = \epsilon_{(2)}$.
One computes the \gram\ matrix
\beq \label{eq:fgram}
f=\left(\begin{array}{lll|l|ll}
\delta&1&0&1&1&1\\
1&\delta&1&1&1&0\\
0&1&\delta&1&1&1\\\cline{1-3}
1&1&\multicolumn{1}{l}{1}&\delta&0&1\\\cline{1-4}
1&1&\multicolumn{2}{c}{1\;\;\,0}&\delta&1\\
1&0&\multicolumn{2}{c}{1\;\;\,1}&1&\delta
\end{array}\right)
\eq
The lines in (\ref{eq:fgram}) 
indicates the restrictions to a basis $\{e_i, i=1,2,3\}$ 
for $\Delta^4_{2,(1)}$ for $l=-1$
(by a mild abuse of notation); and the basis 
$\{e_i, i=1,2,3,4\}$ for  $l=0$. 

As noted, when the determinant vanishes (e.g. as a polynomial in 
$\delta$) then  the indecomposable module $\Delta$ contains a proper  
submodule $L$. Thus this gives the cases where the algebra 
is non-semisimple.
The roots in the example are given by
\[\delta=\left\{\begin{array}{ll}0,\pm\sqrt{2}&\mbox{if}\quad l=-1\\0,1,\frac{\pm\sqrt{17}-1}{2}&\mbox{if}\quad l=0\\
0,0,0,-4,2,2&\mbox{if}\quad l=1\end{array}\right.
\]

%}}}
%{{{ finale

The `classical' first and last cases are known (they may be 
determined from \cite{Martin91} and \cite{Martin09a} respectively),
but the middle case shows that new features appear for intermediate $l$. 
The {\em details} of this example are therefore most intriguing. 
Here, however, we
note only that it will be clear that every such \gram\ matrix is a
finite polynomial in $\delta$. It follows that, with $k=\C$ say,
the rank of the form is submaximal only on a Zariski-closed subset of
$\delta$-parameter space. In other words:

{\mth{
For each $l$, 
$\Jln(\delta)$ is generically semisimple over $\C$, and non-semisimple
    in finitely many $\delta$-values. \Qed
}} 

The example shows that non-semisimple cases exist.
But computing gram determinants is not easy in general. 
In the classical
cases they are most efficiently calculated using {\em translation
  functors}, which combine use of the $F$ and $G$ functors with 
a geometrical structure obtained via 
induction and restriction. 
We begin to address this in \S\ref{ss:bratt}
by looking at branching rules. 

%}}}
%}}}

\section{Standard Bratteli diagrams} \label{ss:bratt}
%{{{ BRATTELI

\newcommand{\res}{{\mbox{Res}}} %% bare restriction functor
\newcommand{\ind}{{\mbox{Ind}}} %% bare induction functor

%{{{ 1

\mdef
Given a pair of semisimple algebras $A \supset B$, with simple modules
labelled by sets $\Lambda(A)$ and $\Lambda(B)$ respectively, 
then the {\em Bratelli diagram} is the graph with vertices 
$\Lambda(A) \sqcup \Lambda(B)$, and an edge 
(of multiplicity $m$) $\mu \rightarrow \nu$ whenever 
simple $B$-module 
$L_\mu$ is a composition factor of the restriction  $\res^A_B L_\nu$
of simple
$A$-module $L_\nu$, of composition multiplicity $m$
(see e.g. \cite{MarshMartin06} for a review and references). 

If $A_0 \supset A_1 \supset A_2 \supset ...$ is a sequence of subalgebras
then the pairwise
Bratelli diagrams may be chained together in the obvious way. This is
the Bratelli diagram for the sequence.

\mdef More generally, suppose that  $A \supset B$ 
(or a sequence as above) are quasihereditary algebras,
and that the restrictions to $B$ of the 
standard modules of $A$ have filtrations by the standard
modules of $B$ (NB such filtrations are unique if they exist, 
since standard modules are a basis for the Grothendieck group
\cite{ClineParshallScott88}). 
Then the {\em standard Bratelli diagram} encodes the 
standard filtration multiplicities in the same way as the ordinary
Bratelli diagram.

%}}}
%{{{ prop pr:srr

{\mpr{ \label{pr:srr}
Fix $l$. There is a standard Bratelli diagram
for the tower $ \{ J_{l,n} \subset J_{l,n+1} \}_{n\in \N}$ 
(inclusion by adding a line on the right).
The standard restriction rules for 
the $J_{l,n+1}$-modules $\Delta_{p,\lambda}$ 
are given as follows. 
For $p < l+2$:
\[
0 \rightarrow  \oplus_{i} \Delta^{(n)}_{p-1,\lambda-e_i} 
  \rightarrow \res^{J_{l,n+1}}_{J_{l,n}} \Delta^{(n+1)}_{p,\lambda} 
  \rightarrow   \oplus_{i} \Delta^{(n)}_{p+1,\lambda+e_i} 
  \rightarrow 0 
\]
where the sums are over addable/removable boxes of the Young diagram 
(we have noted the value of $n$ explicitly
when different values appear in the same sequence). 
For $p = l+2$:
\[
0 \rightarrow  \oplus_{i} \Delta^{(n)}_{p-1,\lambda-e_i} 
  \rightarrow \res^{J_{l,n+1}}_{J_{l,n}} \Delta^{(n+1)}_{p,\lambda} 
  \rightarrow    \Delta^{(n)}_{p+1,\lambda} 
  \rightarrow 0 
\]
For $p > l+2$:
\[
0 \rightarrow  \Delta^{(n)}_{p-1,\lambda}
   \rightarrow \res^{n+1}_n \Delta^{(n+1)}_{p,\lambda} 
  \rightarrow  \Delta^{(n)}_{p+1,\lambda} 
  \rightarrow 0 
\]
}}

%}}}
%{{{ proof

\proof{ The proof closely follows \cite[Prop.13]{Martin94} or \cite{Martin96}. 
The main differences are due to the fact that here we are working,
for $p>l+2$, with
`inflations' of a subgroup $S_{l+2}$ of $S_p$, rather than inflations
of $S_p$ itself.
 
First note that, combinatorially, we may separate the
basis $B_{p,\lambda}$ 
into two subsets: (I) elements in which the component $x$ from 
$J_l(n,1_p,p)$ has the last line propagating
(i.e. northern marked-point $n$ is connected to southern marked-point $p'$);
and (II) elements in which is does not
(i.e. northern marked-point $n$ is connected to some earlier northern
marked-point). 
For example, in 
(\ref{eg:42}) the first three elements have the last line propagating.   

Next note that the subset in (I) is indeed a 
basis of a submodule with regard to
the $J_{l,n-1}$ action (which acts trivially on the last line). 
For $p> l+2$ it is isomorphic, as it were, to the basis of 
$\Delta^{(n)}_{p-1,\lambda}$. 
Furthermore (II) spans a submodule {\em modulo} (I)
--- i.e. it is a basis for the quotient. 
It is easy to check that this quotient module is  
$\Delta^{(n)}_{p+1,\lambda}$, noting again that the last line acts
trivially, so that there is an isomorphism obtained by `deforming' the
last marked-point from the top to the bottom of the picture
(effectively adding another propagating line). 

The $p \leq l+2$ cases are similar. One uses the induction and
restriction rules for symmetric group Specht modules. 
Cf. e.g. \cite[Prop.13]{Martin94}, \cite{DoranWalesHanlon99}.
\Qed } 

%}}}
%{{{ more (Rollet)

\medskip

\newcommand{\Roll}{{\mathsf R}}  %% Rollet graph hook

\begin{figure}
\[
\includegraphics[width=3in]{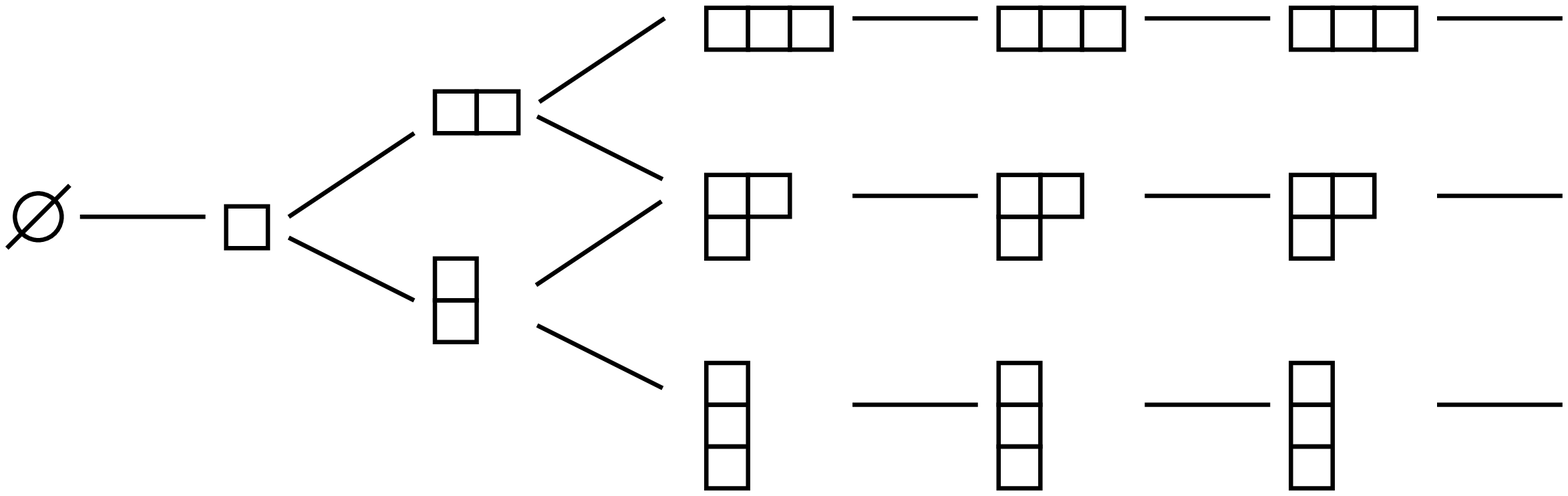}
\]
\medskip
\[
\includegraphics[width=3in]{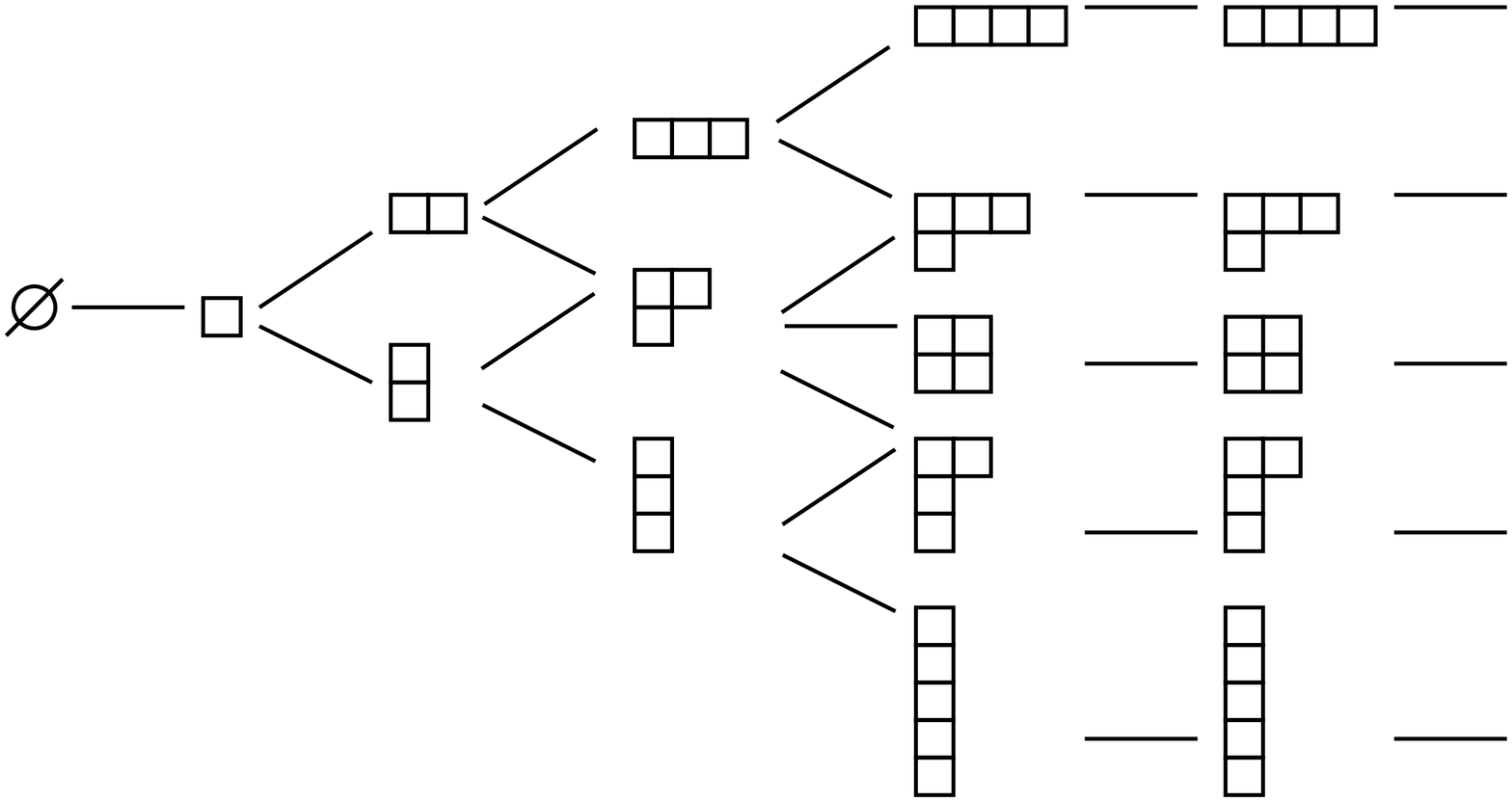}
\]
\caption{The $J_{l,n} \subset J_{l,n+1}$ standard
Rollet
diagram in case $l=1$
and case $l=2$. 
\label{fig:rolletJ0}}
\end{figure}

\mdef \label{de:rollet}
Note in particular that the restriction is multiplicity free here;
and is essentially independent of $n$.
The standard Bratteli diagram in the $l=1$ case 
and in the $l=2$ case 
is encoded in Fig.\ref{fig:rolletJ0}, in the form of the 
corresponding {\em Rollet diagram} $\Roll_l$.
Rollet diagrams are described, for example, in  \cite{MarshMartin06}
(and cf. \cite{Goodman89}).
%(in the sense, for example, of \cite{MarshMartin06}).
In brief, a Rollet diagram arises when there is a global
(large $n$)
limit for the index set for standard or simple modules of a tower of algebras
(such as $J_{l,n}$, $n=0,1,2,3,...$),
induced in the manner of Prop.\ref{pr:eJe};
and {\em also} for the corresponding restriction rules. 
Often, as here, there is a
global limit for the index sets for 
odd and even $n$ sequences of algebras separately
(localisation and globalisation change $n$ in steps of 2).
The Rollet graph is an encoding of all the data in the Bratteli
diagram on a vertex set consisting just of the union of the odd and
even global limits. One takes this vertex set, and an edge of
multiplicity $m$ from $\lambda$ to $\mu$ whenever restriction of 
$\Delta_\lambda$ contains $m$ copies of $\Delta_\mu$
--- this being well-defined by the independence-of-$n$ property.

In particular $\Roll_l$ is an undirected bipartite graph in our
case (cf. the directed Bratteli graph). That is, it is multiplicity
free, and whenever there is an edge there is an edge in both
directions. 

%}}}
%{{{ fig

\mdef
Note that this gives a beautiful combinatorial approach to 
the dimensions of standard modules, and of the various algebras 
(cf. \cite{LeducRam97,MarshMartin06}). 

{\mth{ Let $W^n_{p,\lambda}(l)$ denote the set of walks of length $n$ 
on $\Roll_l$ from
vertex $\emptyset$ to vertex $(p,\lambda)$ (the vertex with label
$\lambda$ at distance $p$ from $\emptyset$). 
Then 
$\dim(\Delta^{n}_{p,\lambda} ) = |W^n_{p,\lambda}(l)|$
and
$\dim (\Jln ) =  |W^{2n}_{0,\emptyset}(l)|$.
\Qed
}}

%}}}

%}}}
%}}}
\section{Discussion}
%{{{ Next steps (Demote to Discussion section)

The construction of $\caJ_{l}$ and $\catP_l$ 
raises  %a number of 
many interesting collateral questions. 
In this section we assemble %%collect together 
some brief general observations on our construction, 
on further developments and on open problems 
(we   
defer full details to a separate note \cite{KadarMartinYu12}). 

\subsection{Next steps in reductive representation theory of $\Jln$}
%{{{ cf Brauer

The next steps  %seek to 
parallel the program for the Brauer algebra used in
\cite{Martin09a}, but now for each $l$ in turn. 
Essentially we should compute the blocks, and
construct `translation functors'; and then construct corresponding analogues of
Kazhdan--Lusztig polynomials. 

We write $\ind -$ for the induction functor adjoint to $\res -$
as in (\ref{pr:srr}) above. 
The precursor of translation functors in the Brauer case is 
the natural isomorphism of functors 
expressed as 
$\ind \cong  \res . G$.
The general setup here is as follows. 

{\mlem{
Suppose we have a sequence of algebras $A \supset B \supset C$ and an
idempotent $e$ in $A$ such that $eAe \cong C$. 
If $B$ and $Ae$ are isomorphic as left-$B$ right-$C$-modules then 
functors
$\ind \; \cong \res\; G $ ($G=G_e$ defined as in (\ref{pr:GeD})).
}}

\proof
We have $\ind_{C}^{B} -  = B_{C} \otimes_{C} - $ and 
$G_{C}^{A} - \; = Ae \otimes_{C} - $.
\qed

\medskip
%}}}
%{{{ pr indresG

In the Brauer case we have the `disk lemma':
Consider a partition in $J(n+1,n+1)$. Replacing the vertex $n+1'$ with
a vertex $n+2$ defines a map 
$\eta: J(n+1,n+1) \rightarrow J(n+2,n)$.
One easily checks that this is a bijection; and an isomorphism of 
$n+1,n$-bimodules.
On the other hand $B_{n+2} e$ has a basis of partitions in 
$J(n+2,n+2)$ with a pair $\{ n+1',n+2'\}$. There is 
a natural bijection of this basis with $J(n+2,n)$ (simply omit the
indicated pair). Altogether then, 
$B_{n+1}$ and $B_{n+2} e$ are isomorphic as bimodules.
In our case we have the following.

{\mlem{
There is a well-defined restriction of $\eta$ to a map
$\eta: \Jlnn{n+1} \rightarrow \Jlnn{n+2} e$; and this 
is an isomorphism of  left-$\Jlnn{n+1}$ right-$\Jlnn{n}$-modules.
}}

\proof  One checks that $\eta$ does not change height. \Qed

Thus we have the following. 

{\mpr{ \label{pr:indresG}
Fix any $l$. Then for all $n$ we have 
$\ind \; \cong \res\; G $.
\qed
}}

This is a  %%very 
powerful result since, for example, the $\ind-$ functor
takes projectives to projectives, while Prop.\ref{pr:srr}
tells us what $\res-$ does to standard modules.
Thus we have an iterative scheme for computing the $\Delta$-content
of projective modules, cf. \cite{Martin09a}.
%
%}}}
%{{{ generalised JBC more
\medskip

\mdef Closely related to Prop.\ref{pr:indresG},
a generalised Jones Basic Construction 
\cite{MarshMartin06} applies here 
(cf. the original Jones Basic Construction \cite{Goodman89}). 
It is analogous to the case in \cite{Martin2000}. 

\mdef We can use the graphs $\Roll_l$ 
from (\ref{de:rollet})
to give an explicit contruction for the
basis states of the standard modules; and indeed of the entire
algebra --- a generalised Robinson-Schensted correspondence \cite{Knuth98}. 
See \cite{KadarMartinYu12} for details.

\newcommand{\appRS}{{  \label{app:RS}
The generalised Robinson-Schensted correspondence
can be illustrated as follows.
In general one needs the {\em usual}
Robinson-Schensted correspondence (see e.g. 
\cite{Knuth98}) to `unwrap' the symmetric group part of the construction
into irreducible components. 
But for low $l$ ($l<2$) we can represent basis for $S_{l+2}$
irreducibles in a simplified way 
(the bases for the two $S_2$ irreducibles, labelled by $(2)$ and $(1^2)$, 
may be represented by the diagrams for $1$ and $\sigma_1$ respectively;
the bases for the three $S_3$ irreducibles, 
labelled by $(3)$, $(2,1)$ (2 basis elements) and $(1^3)$, 
may be represented by the diagrams for $1$,
for $\sigma_1$  and $\sigma_2$, and for $\sigma_1 \sigma_2 \sigma_1$ 
respectively).
See \cite{Martin12w}  %% \cite{PPM-web-page} 
or \cite{KadarMartinYu12}
for an illustration. 

}}

%}}}

%}}}
%{{{ ODDS AND ENDS

\subsection{Remarks on $\Jln$ construction} \label{ss:remarkies}

%{{{ remarks

The $\caJ_{l}$ construction is amenable to several intriguing generalisations. 
Here we briefly mention just one particular such
generalisation, which case makes a contact with existing studies.

\mdef %For $\Jln$ t
The {\em first} case with crossings, $J_{l=0,n}$, is  %closely 
connected to the blob algebra \cite{MartinSaleur94a}:

We say a \ppicture\ is {\em left-simple} if the intersection of the 0-alcove
(as for example in Fig.\ref{fig:squigg})
with the frame of $R$ is connected. 
A partition is left-simple if every picture of it is left-simple. 
(For instance the identity in $B_n$ is left-simple,
while the example in Fig.\ref{fig:squigg} is not.)
Define $\Jl{l}^1(n,m)$ 
as the subset of $\Jl{l}(n,m)$ of left-simple partitions. 

%}}}
%{{{ prop

\mpr{ \label{th:blobx}
The %algebra
subspace $k\Jl{0}^1(n,n)$ is a subalgebra of $J_{0,n}$.
This subalgebra is isomorphic to the blob algebra
$b_{n-1}(q,q')$, with $q,q'$ determined by $\delta$ as follows. 
Parameterising (as in \cite{MartinSaleur94a})
with $x=q+q^{-1}$ as the undecorated loop parameter;
and $y=q'+{q'}^{-1}$ as the decorated loop parameter, 
we have $x=\delta$ and $y=\frac{\delta +1}{2}$.
}

\proof{ (Outline) It is easy to check that the subspace is a subalgebra.
It is also easy to show a bijection between  $\Jl{0}^1(n,m)$
and the set of $(n-1,m-1)$-blob diagrams. This does not lift to an
algebra map, but shows the dimensions are the same.
A {\em heuristic} 
for the algebra isomorphism is to note that the intersection of the
  propagating number zero ideal (\ref{de:propideals}) with the subset of
  left-simple partitions is empty, so there is no
  $\Delta_{0,\emptyset}$ representation. With this node removed, the
  Rollet diagram 
from \S\ref{ss:bratt} 
becomes a (`doubly-infinite') chain, which is the
  same as for the blob algebra. 
See  \cite{KadarMartinYu12} for  an explicit proof.  %%(\Qed)
}

\mdef
On the other hand  %%But 
one can check using \S\ref{ss:bratt} 
that higher $l$ cases such as 
the algebra generated by $\Jl{1}^1(n,n)$ 
do
not coincide  %%either directly or indirectly in this way
with the higher contour algebras \cite{MartinGreenParker07}
or the constructions in \cite{FanGreen99,Green98,tomDieck98}. 

%}}}

%{{{ more (grading/monoidal)

\mdef %Remark: 
Diagram bases may be used to do graded representation theory
(in the sense, for example, of \cite{BrundanKleshchev09})
for graded blob and TL algebras (see e.g. \cite{PlazaRyomHansen12}), 
regarded as quotients of
graded cyclotomic Hecke algebras. It would be interesting to try to
generalise \cite{PlazaRyomHansen12} to $J_{l,n}$.

\mdef As noted (cf.  (\ref{de:catP0}),
(\ref{de:catmonoid}) and (\ref{de:catPBTT})), 
$\catP,\catB,\catT$ are monoidal categories with object monoid
$(\N,+)$ and monoid composition visualised by lateral (as opposed to vertical)
juxtaposition of diagrams. 
Note however that the categories 
$\caJ_l, \; \catP_l$ do not directly inherit this structure 
(except in cases $\caJ_{-1} = \catT$ and $\caJ_{\infty}  = \catB$). 

%}}}

%}}}
%
%{{{ THANKS

\medskip

\noindent
{\bf Acknowledgements.} PM would like to thank Robert Marsh for useful 
conversations. 
ZK and PM would like to thank EPSRC for financial support. 
%}}}
%{{{ APPENDIX

%{{{ APPENDIX OLD

%}}}
%}}}

%{{{ biblio new
\newcommand{\url}[1]{#1}
\bibliographystyle{amsplain}
\bibliography{new31,local}

\providecommand{\bysame}{\leavevmode\hbox to3em{\hrulefill}\thinspace}
\providecommand{\MR}{\relax\ifhmode\unskip\space\fi MR }
% \MRhref is called by the amsart/book/proc definition of \MR.
\providecommand{\MRhref}[2]{%
  \href{http://www.ams.org/mathscinet-getitem?mr=#1}{#2}
}
\providecommand{\href}[2]{#2}
\begin{thebibliography}{10}

\bibitem{AlvarezMartin06}
M~Alvarez and P~P Martin, \emph{A {T}emperley-{L}ieb category for 2-manifolds},
  preprint (2006), arxiv 0711.4777.

\bibitem{BaezDolan95}
John Baez and James Dolan, \emph{Higher-dimensional algebra and topological
  quantum field theory}, Jour. Math. Phys. \textbf{36} (1995), 6073--6105,
  q-alg/9503002.

\bibitem{Baxter}
R~J Baxter, \emph{Exactly solved models in statistical mechanics}, Academic
  Press, New York, 1982.

\bibitem{Benson95}
D~J Benson, \emph{Representations and cohomology {I}}, Cambridge, 1995.

\bibitem{Bethe31}
H~Bethe, Z Phys \textbf{71} (1931), 205.

\bibitem{BirmanWenzl89}
J~S Birman and H~Wenzl, \emph{Braids, link polynomials and a new algebra},
  Transactions AMS \textbf{313} (1989), 249--273.

\bibitem{BjorkenDrell65}
J~D Bjorken and S~D Drell, \emph{Relativistic quantum fields}, McGraw-Hill,
  1965.

\bibitem{Brauer37}
R~Brauer, \emph{On algebras which are connected with the semi--simple
  continuous groups}, Annals of Mathematics \textbf{38} (1937), 854--872.

\bibitem{Brauer39}
\bysame, \emph{On modular and p-adic representations of algebras}, Proc Nat
  Acad Sci USA \textbf{25} (1939), 252--258.

\bibitem{Brown55}
W~P Brown, \emph{An algebra related to the orthogonal group}, Michigan Math. J.
  \textbf{3} (1955-56), 1--22.

\bibitem{BrundanKleshchev09}
J~Brundan and A~Kleshchev, \emph{Blocks of cyclotomic {H}ecke algebras and
  {K}hovanov--{L}auda algebras}, Invent Math \textbf{178} (2009), 451--484.

\bibitem{CautisJackson10}
S~Cautis and D~M Jackson, \emph{On {T}utte's chromatic invariant}, Trans AMS
  \textbf{362} (2010), 509--535.

\bibitem{ClineParshallScott88}
E~Cline, B~Parshall, and L~Scott, \emph{Finite-dimensional algebras and highest
  weight categories}, J. reine angew. Math. \textbf{391} (1988), 85--99.

\bibitem{ClineParshallScott99}
\bysame, \emph{Generic and $q$--rational representation theory}, Publ RIMS
  \textbf{35} (1999), 31--90.

\bibitem{CoxMartinParkerXi06}
A~G Cox, P~P Martin, A~E Parker, and C~C Xi, \emph{Representation theory of
  towers of recollement: theory, notes and examples}, J Algebra \textbf{302}
  (2006), 340--360, DOI 10.1016 online (math.RT/0411395).

\bibitem{CoxDevisscherMartin0609}
A~G Cox, M~De Visscher, and P~P Martin, \emph{A geometric characterisation of
  the blocks of the {B}rauer algebra}, JLMS \textbf{80} (2009), 471--494,
  (math.RT/0612584).

\bibitem{CrowellFox}
R~H Crowell and R~H Fox, \emph{Introduction to knot theory}, Ginn, 1963.

\bibitem{DlabRingel}
V~Dlab and C~M Ringel, Compositio Mathematica \textbf{70} (1989), 155--175.

\bibitem{DoikouMartin03}
A~Doikou and P~P Martin, \emph{Hecke algebraic approach to the reflection
  equation for spin chains}, J Phys A \textbf{36} (2003), 2203--2225,
  hep-th/0206076.

\bibitem{Donkin98}
S~Donkin, \emph{{The $q$-Schur} algebra}, LMS Lecture Notes Series, vol. 253,
  Cambridge University Press, 1998.

\bibitem{DoranWalesHanlon99}
W~F Doran, D~B Wales, and P~J Hanlon, \emph{On semisimplicity of the {B}rauer
  centralizer algebras}, J Algebra \textbf{211} (1999), 647--685.

\bibitem{DuParshallScott98}
J~Du, B~Parshall, and L~Scott, \emph{Quantum {W}eyl reciprocity and tilting
  modules}, Commun Math Phys \textbf{195} (1998), 321--352.

\bibitem{FanGreen99}
C~K Fan and R~M Green, \emph{On the affine {T}emperley--{L}ieb algebras}, J.
  LMS \textbf{60} (1999), 366--380.

\bibitem{Goodman89}
F~M Goodman, P~de~la Harpe, and V~F~R Jones, \emph{Coxeter--{D}ynkin diagrams
  and towers of algebras}, Math Sci Research Inst Publ, Springer--Verlag,
  Berlin, 1989.

\bibitem{Green80}
J~A Green, \emph{Polynomial representations of ${GL}_n$}, Springer-Verlag,
  Berlin, 1980.

\bibitem{Green98}
R.~M. Green, \emph{Generalized {Temperley--Lieb} algebras and decorated
  tangles}, J. Knot Theory Ramifications \textbf{7} (1998), 155--171.

\bibitem{HanlonWales94}
P~Hanlon and D~Wales, \emph{A tower construction for the radical in {B}rauer's
  centralizer algebras}, J Algebra \textbf{164} (1994), 773--830.

\bibitem{James}
G.~D. James, \emph{The representation theory of the symmetric groups}, Lecture
  Notes in Mathematics 682, Springer, 1978.

\bibitem{JamesKerber81}
G~D James and A~Kerber, \emph{The representation theory of the symmetric
  group}, Addison-Wesley, London, 1981.

\bibitem{Jantzen87}
J~C Jantzen, \emph{Representations of algebraic groups}, Academic Press, 1987.

\bibitem{KadarMartinYu12}
Z~Kadar, P~Martin, and S~Yu, \emph{Geometric partition categories}, in
  preparation (2014).

\bibitem{Kassel95}
C~Kassel, \emph{Quantum groups}, Springer, 1995.

\bibitem{Kauffman91}
L~H Kauffman, \emph{Knots and physics}, World Scientific, Singapore, 1991.

\bibitem{Knuth98}
D~E Knuth, \emph{Sorting and searching}, 2 ed., The Art of Computer
  Programming, vol.~3, Addison Wesley, 1998.

\bibitem{KoenigXi99}
S~Koenig and C~C Xi, \emph{Cellular algebras: inflations and {M}orita
  equivalences}, Journal of the LMS \textbf{60} (1999), 700--722.

\bibitem{LeducRam97}
R~Leduc and A~Ram, \emph{A ribbon {H}opf algebra approach to the irreducible
  representations of centralizer algebras: the {B}rauer, {Birman-Wenzl}, and
  type-${A}$ {Iwahori--Hecke} algebras}, Adv. Math. \textbf{125} (1997), 1--94.

\bibitem{LiebMattis66}
E~H Lieb and D~C Mattis, \emph{Mathematical physics in one dimension}, Academic
  Press, 1966.

\bibitem{MarshMartin06}
R~Marsh and P~Martin, \emph{Pascal arrays: Counting {C}atalan sets}, preprint
  (2006), arXiv:math/0612572v1 [math.CO].

\bibitem{Martin12w}
P~P Martin, \emph{web page},
  http://www1.maths.leeds.ac.uk/$\sim$ppmartin/GRAD/pcat-ex-sub.html.

\bibitem{Martin91}
\bysame, \emph{Potts models and related problems in statistical mechanics},
  World Scientific, Singapore, 1991.

\bibitem{Martin92}
\bysame, \emph{On {S}chur-{W}eyl duality, {$A_n$} {H}ecke algebras and quantum
  {$sl(N)$}}, Int J Mod Phys A \textbf{7 suppl.1B} (1992), 645--674.

\bibitem{Martin94}
\bysame, \emph{Temperley--{L}ieb algebras for non--planar statistical mechanics
  --- the partition algebra construction}, Journal of Knot Theory and its
  Ramifications \textbf{3} (1994), no.~1, 51--82.

\bibitem{Martin96}
\bysame, \emph{The structure of the partition algebras}, J Algebra \textbf{183}
  (1996), 319--358.

\bibitem{Martin2000}
\bysame, \emph{The partition algebra and the {P}otts model transfer matrix
  spectrum in high dimensions}, J Phys A \textbf{32} (2000), 3669--3695.

\bibitem{Martin08a}
\bysame, \emph{On diagram categories, representation theory and statistical
  mechanics}, AMS Contemp Math \textbf{456} (2008), 99--136.

\bibitem{Martin09b}
\bysame, \emph{Lecture notes in representation theory}, unpublished lecture
  notes, 2009.

\bibitem{Martin09a}
\bysame, \emph{The decomposition matrices of the {B}rauer algebra over the
  complex field}, Trans AMS, to appear (http://arxiv.org/abs/0908.1500, 2009).

\bibitem{MartinGreenParker07}
P~P Martin, R~M Green, and A~E Parker, \emph{Towers of recollement and bases
  for diagram algebras: planar diagrams and a little beyond}, J Algebra
  \textbf{316} (2007), 392--452, (math.RT/0610971).

\bibitem{MartinSaleur94a}
P~P Martin and H~Saleur, \emph{The blob algebra and the periodic
  {T}emperley--{L}ieb algebra}, Lett. Math. Phys. \textbf{30} (1994), 189--206,
  (hep-th/9302094).

\bibitem{MartinWoodcock2000}
P~P Martin and D~Woodcock, \emph{On the structure of the blob algebra}, J
  Algebra \textbf{225} (2000), 957--988.

\bibitem{Mazorchuk95}
V~Mazorchuk, \emph{On the structure of {B}rauer semigroup and its partial
  analogue}, Problems in Algebra \textbf{13} (1998), 29--45.

\bibitem{Moise77}
E~E Moise, \emph{Geometric topology in dimensions 2 and 3}, Graduate Texts in
  Mathematics 47, Springer-Verlag, New York, 1977.

\bibitem{Morton10}
H~R Morton, \emph{A basis for the {Birman-Murakami-Wenzl} algebra}, unpublished
  (2010).

\bibitem{Murakami87}
J~Murakami, \emph{The {K}auffman polynomial of links and representation
  theory}, Osaka J. Math. \textbf{24} (1987), no.~4, 745--758.

\bibitem{PlazaRyomHansen12}
D~Plaza and S~Ryom-Hansen, \emph{Graded cellular bases for {T}emperley--{L}ieb
  algebras of type {A} and {B}}, J Alg Comb (to appear),
  http://arxiv.org/abs/1203.2592.

\bibitem{Putcha98}
M~Putcha, \emph{Complex representations of finite monoids {II}. highest weight
  categories and quivers}, J Algebra \textbf{205} (1998), 53--76.

\bibitem{ReshetikhinTuraev90}
N~Yu Reshetikhin and V~G Turaev, \emph{Ribbon graphs and their invariants
  derived from quantum groups}, Comm Math Phys \textbf{127} (1990), 1--26.

\bibitem{TemperleyLieb71}
H~N~V Temperley and E~H Lieb, \emph{Relations between percolation and colouring
  problems and other graph theoretical problems associated with regular planar
  lattices: some exact results for the percolation problem}, Proceedings of the
  Royal Society A \textbf{322} (1971), 251--280.

\bibitem{tomDieck98}
T~tom Dieck, \emph{Categories of rooted cylinder ribbons and their
  representations}, J reine angew Math \textbf{494} (1998), 35--63.

\bibitem{Wenzl88b}
H~Wenzl, \emph{On the structure of {B}rauer's centralizer algebras}, Ann. Math.
  \textbf{128} (1988), 173--193.

\bibitem{Yu05}
S.~Yu, \emph{The cyclotomic {B}irman-{M}urakami-{W}enzl algebras}, Ph.D.
  thesis, The University of Sydney, December 2007,
  \url{http://arxiv.org/abs/0810.0069}.

\end{thebibliography}
%}}}
\end{document}